\documentclass[11pt,twoside,a4paper]{amsart}
\usepackage[latin1]{inputenc}
\usepackage[T1]{fontenc}
\usepackage{amsmath,amssymb}
\usepackage{array}
\usepackage{color}
\usepackage{graphicx}

\newcommand{\Id}{\mathbf{1}}

\newcommand{\rank}{\operatorname{rank}}

\renewcommand{\>}{\rangle}
\newcommand{\<}{\langle}

\newcommand{\IM}{\operatorname{\mathcal{M}}}
\newcommand{\BB}{\operatorname{\mathcal{B}}}
\newcommand{\LL}{\operatorname{\mathcal{L}}}
\newcommand{\EE}{\operatorname{\mathcal{E}}}
\newcommand{\UU}{\operatorname{\mathcal{U}}}
\newcommand{\FF}{\operatorname{\mathcal{F}}}

\newcommand{\TiM}{\widetilde{\operatorname{\mathcal{M}}}}
\newcommand{\IA}{\operatorname{\mathfrak{A}}}
\newcommand{\IP}{\operatorname{\mathfrak{P}}}

\newcommand{\Ih}{\operatorname{\textbf{h}}}

\newcommand{\IO}{\operatorname{\mathcal{O}}}
\newcommand{\Si}{\dot{S}}

\newcommand{\CP}{\mathbb{C}\mathbb{P}^1}
\newcommand{\eps}{\epsilon}
\newcommand{\hh}{\operatorname{\textbf{h}}}
\newcommand{\CR}{\bar{\partial}}

\newcommand{\uw}{\underline{w}}
\newcommand{\Mph}{M_{\phi}}
\newcommand{\CM}{\overline{\IM}}

\newcommand{\Aut}{\operatorname{Aut}}
\newcommand{\Symp}{\operatorname{Symp}}
\newcommand{\Det}{\operatorname{Det}}

\newcommand{\Pic}{\operatorname{Pic}}
\newcommand{\Sp}{\operatorname{Sp}}

\newcommand{\cst}{\operatorname{const}}
\newcommand{\op}{\operatorname{op}}

\newcommand{\ev}{\operatorname{ev}}

\newcommand{\loc}{\operatorname{loc}}

\newcommand{\Tor}{\operatorname{Tor}}

\newcommand{\hb}{\bar{h}}

\newcommand{\diag}{\operatorname{diag}}
\newcommand{\ind}{\operatorname{ind}}
\newcommand{\coker}{\operatorname{coker}}
\newcommand{\Coker}{\operatorname{Coker}}
\newcommand{\Ker}{\operatorname{Ker}}

\renewcommand{\Im}{\operatorname{Im}}
\newcommand{\del}{\partial}

\newcommand{\virt}{\operatorname{virt}}
\newcommand{\point}{\operatorname{point}}

\newcommand{\RS}{\IR \times S^1}

\renewcommand{\qed}{\square}

\newcommand{\im}{\operatorname{im}}
\newcommand{\IC}{\operatorname{\mathbb{C}}}
\newcommand{\IZ}{\operatorname{\mathbb{Z}}}

\newcommand{\IR}{\operatorname{\mathbb{R}}}
\newcommand{\IN}{\operatorname{\mathbb{N}}}

\author{Oliver Fabert}
\title{Obstruction bundles over moduli spaces\\ with boundary and the action filtration\\ in symplectic field theory} 
\thanks{Research supported by the German Research Foundation (DFG)}
\begin{document}

\abstract 
Branched covers of orbit cylinders are the basic examples of holomorphic curves studied in 
symplectic field theory. Since all curves with Fredholm index one can never be regular for any choice of cylindrical almost complex structure, we generalize the obstruction bundle technique of Taubes for determining multiple cover contributions from Gromov-Witten theory to the case of moduli spaces with boundary. Our result proves that the differential in 
rational symplectic field theory and contact homology is strictly decreasing with respect to the natural action filtration.} 

\maketitle

\markboth{OLIVER FABERT}{OBSTRUCTION BUNDLES AND THE ACTION FILTRATION IN SFT}

\tableofcontents

\section*{Main result and summary}

Symplectic field theory (SFT), introduced by H. Hofer, A. Givental and Y. Eliashberg in 2000 ([EGH]), is a very large project designed to describe in a unified way the theory of holomorphic curves in symplectic and contact manifolds and can be viewed as a topological quantum field theory approach to Gromov-Witten theory. Besides providing a unified view on known results, symplectic field theory leads to numerous new applications and opens new routes yet to be explored. While the theory promises to provide lots of new invariants with rich algebraic structures, which are currently by a large number of researchers, these invariants are in general very difficult to compute. \\

This paper is concerned with the basic examples of punctured holomorphic curves studied in rational symplectic field theory. Although the definition of symplectic field theory suggests that one has to count holomorphic curves in cylindrical manifolds $\IR\times V$ equipped with a cylindrical almost complex structure $J$, it is already known from Gromov-Witten theory that due to the presence of multiply covered curves we in general cannot achieve transversality for all moduli spaces even for generic choices of $J$. \\

While the contribution of the orbit curves in cylindrical contact homology, namely cylinders staying over one orbit, is still immediately clear, observe that the basic examples of punctured holomorphic curves studied in general symplectic field theory are not only these orbit cylinders but also their branched covers. We show that these multiple covers are in fact the reason why transversality for generic $J$ in general fails in symplectic field theory and whose contribution to the theory is hence not a priori clear. \\

Indeed it is easy to show that in every case where these orbit curves would contribute to the algebraic invariants by index reasons, transversality for the Cauchy-Riemann operator can never be satisfied, so that one has to perturb the Cauchy-Riemann operator appropriately and count elements in the resulting regular moduli spaces. Here it is important that the perturbation chosen for different moduli spaces are compatible with compactness and gluing in symplectic field theory. \\

Instead of referring to the polyfold project of Hofer, Wysocki and Zehnder, which promises to prove transversality for all moduli spaces in symplectic field theory, but which is currently not yet completely finished, we generalize the idea of Taubes of computing Euler numbers of obstruction bundles for determining the contribution of multiple covers from Gromov-Witten theory to the case of moduli spaces with codimension one boundary, as appearing in the study of pseudoholomorphic curves with punctures or boundary. Apart from the reason that we can give a rigoros proof of transversality, the main motivation lies in the fact that we really want to compute the contribution of this fundamental class of moduli spaces to the SFT invariants. More precisely, we prove the following \\
\\
{\bf Main Theorem:} {\it We can choose compact perturbations of the Cauchy-Riemann operator, which make all moduli spaces of orbit curves regular in a way compatible with compactness and gluing, such that the algebraic counts of elements in all resulting zero-dimensional regular moduli spaces (modulo $\IR$-shift) are zero.} \\ 
  
Although the result of our computation may suggest that it follows a global symmetry of the resulting regular moduli space, we want to emphasize that our proof crucially relies on the existence of a smooth finite-dimensional obstruction bundle over all moduli spaces of orbit curves. In particular, while we will show that there exists a global $S^1$-action on the underlying nonregular moduli space of branched covers, this $S^1$-action in general does not lift to an action on the obstruction bundle over this space, so that the resulting perturbed moduli space does {\it not} carry a global symmetry. Finally, we will show that we in general cannot expect to get a well-defined count of elements in all moduli spaces provided that the number of elements in one of the moduli spaces would be different zero, which would make the computation almost impossible. \\

For the significance of this result for symplectic field theory we show that our result has the following immediate \\
\\
{\bf Corollary:} {\it The differential in contact homology and rational symplectic field theory is strictly decreasing with respect to the natural action filtration.} \\

In particular, we outline that the statement of the theorem is true for any choice of coherent compact perturbations 
chosen to make the moduli spaces of symplectic field theory regular. \\

We further introduce the rational symplectic field theory of a single closed Reeb orbit and use our result to compute the underlying generating function. Including the even more general picture outlined in [EGH] needed to view Gromov-Witten 
theory as a part of symplectic field theory, we further prove what we get when we additionally introduce a string of closed 
differential forms $\theta_1,...,\theta_N\in\Omega^*(V)$. Here we prove by using our main result that the Hamiltonian only sees the homology class represented by the underlying closed Reeb orbit. It follows that the Hamiltonian is in general no longer equal to zero when a string of differential forms is chosen, which implies that the differential in rational symplectic field theory and contact homology is no longer strictly decreasing with respect to the action filtration. However, we follow [FOOO] in employing the spectral sequence for filtered complexes to prove the following other important consequence of our main theorem. \\
\\
{\bf Corollary:} {\it Consider a contact manifold or a symplectic mapping torus. Then the $E^2$-page of the spectral sequence $(E^r,d^r)$ for the action filtration on contact homology is given by the graded commutative algebra which is freely generated by the formal variables $q_{\gamma}$ with $\int_{\gamma}\theta_i=0$ for all $i=1,...,N$.} \\ 

Note that this in turn provides us with an easy method to show the vanishing of contact homology: Assume that the string of closed differential forms is chosen in such a way that it indeed generates the cohomology of the 
target manifold (and that none of the corresponding formal variables is set to zero). Then the contact homology vanishes if there are no null-homologous Reeb orbits, like in the case of symplectic mapping tori and unit cotangent bundle of tori.\\
  
For the proof of the main theorem we have to show that for every moduli space of orbit curves the cokernels of the linearizations of the 
Cauchy-Riemann operator indeed fit together to give a global vector bundle over the corresponding compactified moduli space, and prove that there 
exists an Euler number for coherent, that is, gluing-compatible sections in the cokernel bundle which is zero. \\

While in Gromov-Witten theory the existence of the Euler number is immediately clear since all moduli spaces are pseudo-cycles, i.e., 
homologically have no boundary, but their computation is hard in general, the opposite is true here: Since the algebraic invariants of 
symplectic field theory rely on the codimension one boundary phenomena of the moduli spaces of punctured curves, i.e., the regular moduli spaces 
define relative rather than absolute virtual moduli cycles, Euler numbers for Fredholm problems in general do not exist since the count of zeroes in 
general depends on the compact perturbations chosen 
for the moduli spaces in the boundary. \\

In this paper we make use of the fact that the moduli spaces in the boundary again consist of 
branched covers of trivial cylinders and prove the existence of the Euler number by induction on the number of punctures. 
For the induction step we do not only use that there exist Euler numbers for the moduli spaces in the boundary, but it is further important 
that all these Euler numbers are in fact zero. The vanishing of the Euler number in turn can be deduced from the different parities 
of the actual and the virtual dimensions of the moduli spaces following the idea for the vanishing of the Euler characteristic for 
odd-dimensional manifolds. \\

From some invariance argument we deduce that, once the analytical foundations of symplectic field theory are 
established, the result about sections in the cokernel bundles suffices to prove 
that the algebraic number of elements in the regular moduli spaces, obtained by adding general compact perturbations to the Cauchy-Riemann operator
are still zero even when the abstract perturbations no longer result from sections in the cokernel bundles. Despite the analytical work in order to 
show that the cokernels fit together to give a nice vector bundle and showing that studying sections in it gives the right result, the strategy of our proof 
indeed only relies on the difference of the parity of the Fredholm index, i.e., the virtual dimension of the moduli space, 
and the actual dimension of the moduli space, including the moduli spaces in the boundary. Hence it should be applicable to a 
wide range of other multiple cover problems in pseudoholomorphic curve theories. \\
\\
{\bf Remark:} Note that in order to prove $d^2=0$ in embedded contact homology and periodic Floer homology the authors of [HT1] and [HT2] 
also study sections in obstruction bundles over moduli spaces of branched covers of trivial cylinders. 
Beside the fact that their papers became available shortly before this project was finished, we emphasize that there is an 
essential difference between their project and ours: While we view the branched covers of orbit cylinders as basic examples of curves counted 
in the differential of rational symplectic field theory and therefore count orbit curves of Fredholm index {\it one}, M. Hutchings and C. Taubes 
developed a generalized gluing theory for symplectic field theory in dimension four where orbit curves of Fredholm index {\it zero} are inserted between the 
curves to be glued. \\

After describing the geometric setup underlying symplectic field theory, we focus on the 
basic facts about orbit curves in symplectic field theory. Since we have to deal with nonregular moduli spaces we introduce coherent abstract perturbations. 
We then rigorously describe the moduli spaces $\IM$ and $\IM^0$ of orbit curves, obtained by quotiening out or not quotiening out the $\IR$-action, 
and their compactifications. 
We show that $\IM$ and $\IM^0$ are given as products involving the moduli space of punctured spheres and use the 
conservation of energy to describe their compactifications $\CM$ and $\overline{\IM^0}$ which are again made up of moduli spaces of orbit curves. 
Introducing the notion of a tree with (based) level structure $(T,\LL)$, $(T,\LL,\ell_0)$ we show that $\CM$ and $\overline{\IM^0}$
carry natural stratifications and prove that $\CM$ and $\overline{\IM^0}$ are smooth manifolds with corners. 
While for this we explicitly describe the compactifications using Fenchel-Nielsen coordinates on the moduli space of punctured 
spheres, we emphasize that the compactifications $\CM$ and $\overline{\IM^0}$ are different from the one obtained using the 
Deligne-Mumford-Knudsen compactification of the moduli space of punctured spheres, in particular, $\CM$ and $\overline{\IM^0}$ have 
codimension one boundary strata. We then introduce the cokernel bundles $\overline{\Coker}\CR_J$ and $\overline{\Coker_0}\CR_J$ over 
the compactified moduli spaces $\CM$ and $\overline{\IM^0}$. After describing the neccessary Banach space bundle setup, we study 
the linearization of the Cauchy-Riemann operator $\CR_J$ and prove that, by energy reasons, the kernel of the linearization of $\CR_J$ 
agrees with the tangent space to the moduli space. This proves in particular that the cokernel of the linearization of 
$\CR_J$ has the same dimension at every point in $\IM$ and $\IM^0$ which is sufficient to prove that $\overline{\Coker}\CR_J$ and 
$\overline{\Coker_0}\CR_J$ naturally carry the structure of a smooth vector bundle over the strata of $\CM$ and $\overline{\IM^0}$, respectively. 
In order to show that these bundles over the strata fit together to a smooth vector bundle over the manifold with corners $\CM$ and 
$\overline{\IM^0}$ we prove a linear gluing result for cokernel bundles. While we show that the construction of coherent 
orientations in [BM] together with the complex orientations of the strata of $\CM$ and $\overline{\IM^0}$ equips the cokernel bundle 
with an orientation over each stratum, it follows from the results in [BM] that these orientations in general do not fit together to give an orientation of the whole cokernel bundles $\overline{\Coker}\CR_J$ and $\overline{\Coker_0}\CR_J$ but differ by a fixed sign due 
to reordering of the punctures. Equipped with the neccessary analytical results about $\overline{\Coker}\CR_J$ and $\overline{\Coker_0}\CR_J$ 
we finally prove the main theorem. After showing that sections in the cokernel bundle indeed provide us with the 
desired compact perturbations for the Cauchy-Riemann operator, we discuss the gluing compatibility for sections in the cokernel bundle 
and define the notion of coherent sections in $\overline{\Coker}\CR_J$. We finally prove by induction that there exists an Euler number 
for coherent sections in $\overline{\Coker}\CR_J$ and show that it is zero. For this we study sections in the cokernel bundle $\overline{\Coker_0}\CR_J$ 
over $\overline{\IM^0}$. We again emphasize that the induction step does not only need 
the existence result of Euler numbers for the moduli spaces in the boundary but also that these numbers are indeed zero. After this we discuss the implications 
of our result on rational symplectic field theory once the analytical foundations are proven. After explaining why the conclusion of the main result should continue 
to hold for all choices of coherent compact perturbations, we introduce the natural action filtration on symplectic field theory. We further introduce the 
rational symplectic field theory of a single closed Reeb orbit. Including the even more general picture outlined in [EGH] needed to view Gromov-Witten 
theory as a part of symplectic field theory, we further prove what we get when we additionally introduce a string of closed differential forms. 
Finally we follow [FOOO] in employing the spectral sequence for filtered complexes, where we use our result to show 
that after passing from the $E^1$-page to the $E^2$-page we only have to consider those formal variables, where the homology class of the underlying 
closed orbit is annihilated by all chosen differential forms. \\

This paper is organized as follows: After an introductory section on orbit curves in symplectic field theory, section one is concerned 
with the nonregular moduli spaces of unperturbed branched covers of trivial cylinders. While section two is devoted to establishing the existence and the 
properties of the cokernel bundle, we prove the main theorem in section three. In section four we finally introduce the symplectic field theory of a single Reeb orbit and discuss the implications of 
our result on rational symplectic field theory. \\
\\
{\bf Acknowledgements:} This research was supported by the German Research Foundation (DFG). 
The author thanks U. Frauenfelder, H. Hofer, K. Mohnke and K. Wehrheim for useful discussions. Special thanks finally go to my advisor Kai Cieliebak and to 
Dietmar Salamon, who gave me the chance to stay at ETH Zurich for the winter term 2006/07, for their support.  

\setcounter{section}{-1}
\section{Introduction}

\subsection{Symplectic field theory}
We start with describing the geometric setup needed to define symplectic field theory.
Following [BEHWZ] and [CM] a Hamiltonian structure on a closed $(2m-1)$-dimensional manifold $V$ is a closed two-form $\omega$ 
on $V$ which is maximally nondegenerate in the sense that $\ker\omega=\{v\in TV:\omega(v,\cdot)=0\}$ is a one-dimensional distribution. 
Note that here we (and [CM]) differ slightly from [EKP]. The Hamiltonian structure is required to be stable in the sense that there exists a 
one-form $\lambda$ on $V$ such that $\ker\omega\subset\ker d\lambda$ and $\lambda\neq 0$ on $\ker\omega$. 
Every stable Hamiltonian structure $(\omega,\lambda)$ defines a symplectic hyperplane distribution $(\xi=\ker\lambda,\omega_{\xi})$, 
where $\omega_{\xi}$ is the restriction of $\omega$ to $\xi$, and 
a vector field $R$ on $V$ by requiring $R\in\ker\omega$ and $\lambda(R)=1$ which is called the Reeb vector field of the 
stable Hamiltonian structure. For the rest of this paper we assume that $(\omega,\lambda)$ is chosen such that all closed orbits $\gamma$ 
of the Reeb vector field are nondegenerate in the sense of the linearized Poincare return map along $\gamma$ has no eigenvalue equal to one. 
Furthermore the boundary of a symplectic manifold is called stable if the restriction of 
the symplectic form to the boundary defines a stable Hamiltonian structure on it. Note that this notion of stability agrees 
with the one in [HZ] as shown in lemma 2.3 in [CM]. Following [BEHWZ], examples of closed manifolds $V$ with a stable Hamiltonian structure 
are contact manifolds, principal circle bundles and symplectic mapping tori. For this note that when $\lambda$ is a contact form on $V$ 
then it is easy to check that 
$(\omega:=d\lambda,\lambda)$ is a stable Hamiltonian structure and that the symplectic hyperplane distribution agrees with the 
contact structure. For the other two cases, let $(M,\omega_M)$ be a symplectic manifold. Then every principal circle bundle 
$S^1\to V\to M$ and every symplectic mapping torus $M\to V\to S^1$, i.e., $V=\Mph=\IR\times M/\{(t,p)\sim(t+1,\phi(p))\}$ for 
$\phi\in\Symp(M,\omega_M)$ also carries a stable Hamiltonian structure. For the principal circle bundle the Hamiltonian 
structure $\omega$ is given by the pullback $\pi^*\omega_M$ under the bundle projection and the one-form $\lambda$ is given by the 
choice of a $S^1$-connection form. On the other hand, the stable Hamiltonian structure on the mapping torus $V=\Mph$ is given by 
lifting $\omega_M$ to $\omega\in\Omega^2(\Mph)$ via the natural flat connection $TV=TS^1\oplus TM$ and setting 
$\lambda=dt$ for the natural $S^1$-coordinate $t$ on $\Mph$. \\

Symplectic field theory assigns algebraic invariants to closed manifolds $V$ with a stable Hamiltonian structure, while 
to symplectic manifolds with stable boundary it assigns homomorphisms between the algebraic invariants of the closed 
manifolds forming the boundary. In this sense symplectic field theory is an abstract quantum field theory, since it can be viewed 
as a functor from a geometric category to an algebraic category. The algebraic invariants are defined by counting $J$-holomorphic 
curves in $\IR\times V$ with finite energy, where the underlying closed Riemann surfaces are explicitly allowed to have punctures. 
The almost complex structure $J$ on the cylindrical manifold $\IR\times V$ is required to be cylindrical and $(\omega,\lambda)$-compatible 
in the sense that it is $\IR$-independent, links the Reeb vector field $R$ with the $\IR$-direction $\del_s$ by $J\del_s=R$, and turns the 
symplectic hyperplane distribution on $V$ into a complex subbundle of $TV$, 
$\xi=TV\cap JTV$. It follows that a cylindrical almost complex structure $J$ on $\IR\times V$ is determined by its restriction 
$J_{\xi}$ to $\xi$ which is furthermore required to be compatible with the symplectic form $\omega_{\xi}$ on $\xi$. 
Following [BEHWZ] the energy $E(u)$ of a punctured $J$-holomorphic curve $u=(a,f): \Si\to \IR\times V$ is given by the sum of the $\lambda$- and the 
$\omega$-energy of $u$, 
\begin{eqnarray*}
 E_{\lambda}(u) = \sup_A \int_{\Si}\alpha(a)\; da\wedge f^*\lambda,\;\; E_{\omega}(u)=\int_{\Si} f^*\omega,
\end{eqnarray*} 
where $A$ denotes the set of all smooth functions $\alpha:\IR\to\IR^+_0$ with compact support and $L^1$-norm equal to one. 
It follows that $E_{\lambda}(u),E_{\omega}(u)$ are nonnegative and, following proposition 5.8 in [BEHWZ], that all 
punctured $J$-holomorphic curves with $E(u)<\infty$ are asymptotically cylindrical over a periodic orbit of the Reeb vector field $R$ 
in the neighborhood of each puncture as long as all periodic orbits are nondegenerate in the sense of [BEHWZ], i.e., one is not an eigenvalue of 
the linearized return map restricted to the symplectic hyperplane distribution. {\it For the following expositions we assume that the 
stable Hamiltonian structure is generic in the sense that all periodic orbits are nondegenerate.} \\ 
 
\subsection{Orbit curves in symplectic field theory} 
Beside the constant curves with no punctures, which do not contribute 
to the differential by algebraic reasons, note that for each closed orbit $\gamma$ of the vector field $R$ we have 
the trivial cylinder $\IR\times\gamma$ as basic example of a $J$-holomorphic curve in $\IR\times V$, where the $J$-holomorphicity follows from 
$J\partial_s = R = \dot{\gamma}$. While these trivial cylinders correspond to the trivial connecting orbits in Floer homology and by the same arguments 
turn out to be irrelevant for the algebraic invariants, it is important to observe that in contact homology and (rational) symplectic field theory we get 
from a single trivial cylinder infinitely many other basic examples of punctured $J$-holomorphic curves with two or more punctures by considering 
branched and unbranched covers of the given trivial cylinder. While the unbranched covers are again trivial cylinders over a multiple 
of the underlying Reeb orbit, it follows (see proposition 1.1) that the branched covers are in one-to-one correspondence with  
meromorphic functions on the underlying closed Riemann surface by removing zeroes and poles and identifying 
$\CP-\{0,\infty\}\cong\RS \cong(\IR\times\gamma,J)$. While these curves are clearly basic in the above sense, it is important to observe that they are also basic from another viewpoint: \\

Like the constant curves and cylinders staying over one orbit are the only holomorphic curves in Gromov-Witten theory and symplectic Floer homology 
with trivial energy, the branched and unbranched covers of trivial cylinders are the only punctured holomorphic curves with vanishing $\omega$-energy. 
Indeed, if $u=(a,f):\Si\to\IR\times V$ has $E_{\omega}(u)=0$ it follows, see lemma 5.4 in [BEHWZ], that $df\in\ker\omega=\IR\times R$, 
so that the image of the $V$-component $f$ 
is a closed Reeb orbit. On the other hand, assuming as in [EGH] that the first homology group of $V$ is torsion-free, observe that after choosing 
a basis for $H_1(V)$ and choosing for each simple orbit $\gamma$ a spanning surface $f_{\gamma}$ in $V$ realizing a cobordism between $\gamma$ 
and a suitable linear combination of these basis elements as in [EGH], we can define an action
\begin{equation*} S(\gamma) = \int f_{\gamma}^*\omega, \end{equation*} 
for every simple closed Reeb orbit $\gamma$. On the other hand, note that for a multiply covered orbit $\gamma^m$ we can use the formal multiple 
$f_{\gamma}^m$ of the spanning surface $f_{\gamma}$ to realize a cobordism between $\gamma^m$ and a linear combination of basis elements, so that 
$S(\gamma^m)=m\cdot S(\gamma)$. Then $E_{\omega}(u)$ can be expressed as the difference of the actions of the closed orbits 
$\gamma_1^{\pm},...,\gamma_{n^{\pm}}^{\pm}$ corresponding to positive, respectively negative punctures of $u$ 
and the $\omega$-area of the homology class $A\in H_2(V)$ which we can assign to $u$ using the spanning surfaces for the simple orbits underlying 
$\gamma_1^{\pm},...,\gamma_{n^{\pm}}^{\pm}$,  
\begin{equation*} 
 E_{\omega}(u) = \sum_{k=1}^{n^+} S(\gamma_k^+) - \sum_{\ell=1}^{n^-} S(\gamma_{\ell}^-) + \omega(A).
\end{equation*}
In particular, it follows that the moduli spaces $\IM_{g,0}(\gamma^{m^+_1},...,\gamma^{m^+_{n^+}};\gamma^{m^-_1},...,\gamma^{m^-_{n^-}})$ 
of $J$-holomorphic curves of genus $g$ in $\IR\times V$ which are asymptotically cylindrical over the multiple covers 
$\gamma^{m^+_1},...,\gamma^{m^+_{n^+}}$ of $\gamma$ at the positive, over $\gamma^{m^-_1},...,\gamma^{m^-_{n^-}}$ at the negative punctures 
and represent the homology class $A=0\in H_2(V)$, entirely consist of multiple covers of the trivial cylinder over $\gamma$. For this observe that 
$m^+_1 + ... + m^+_{n^+} = m^-_1 + ... + m^-_{n^-}$ since else the moduli space is empty by homological reasons, so that  
\begin{equation*}
 \sum_{k=1}^{n^+} S(\gamma^{m^+_k}) - \sum_{\ell=1}^{n^-} S(\gamma^{m^-_{\ell}}) = (\sum_{k=1}^{n^+} m^+_k - \sum_{\ell=1}^{n^-} m^-_{\ell}) \cdot S(\gamma) = 0. 
\end{equation*}

For the rest of this paper we restrict ourselves to the case of rational curves, i.e., with genus $g=0$. Note that the moduli space 
$\IM_{0,0}(\gamma^{m^+_1},...,\gamma^{m^+_{n^+}};\gamma^{m^-_1},...,\gamma^{m^-_{n^-}})$ contributes to the differential  
in rational symplectic field theory only when its virtual dimension given by the Fredholm index of the linearization of the Cauchy-Riemann 
operator $\CR_J$, 
\begin{equation*} 
 \ind\CR_J=\sum_{k=1}^{n^+} \mu_{CZ}(\gamma^{m^+_k}) - \sum_{\ell=1}^{n^-}\mu_{CZ}(\gamma^{m^-_{\ell}}) + (m-3) \cdot (2-n),
\end{equation*}
is equal to one, where $n=n^++n^-$ is the number of punctures and $\dim V=2m-1$. For this observe that under the assumption that  
the Cauchy-Riemann operator meets the zero section transversally in a suitable Banach space bundle over a Banach manifold of maps this 
implies that the moduli space is one-dimensional, i.e., discrete after quotiening out the natural $\IR$-action.
While for trivial cylinders the Fredholm index is always zero, there indeed exist examples of branched covers with 
Fredholm index one. For example it is easy to check that the moduli spaces $\IM_{0,0}(\gamma^2;\gamma,\gamma)$ and $\IM_{0,0}(\gamma,\gamma;\gamma^2)$ 
of pairs of pants mapping to the trivial cylinder over an arbitrary hyperbolic orbit $\gamma$ in a three-manifold have virtual dimension equal to one and 
therefore, in contrast to the underlying trivial cylinder, possibly contribute to the algebraic invariants of rational symplectic field 
theory. On the other hand we prove in proposition 1.1 that when the number of punctures $n=n^++n^-$ is greater or equal to three the moduli space is given by  
\begin{equation*}
 \IM_{0,0}(\gamma^{m^+_1},...,\gamma^{m^+_{n^+}};\gamma^{m^-_1},...,\gamma^{m^-_{n^-}}) 
  = \IR\times S^1 \times \IM_{0,n^++n^-} \times \IZ_{m^+} \times \IZ_{m^-},
\end{equation*}
where $\IM_{0,n^++n^-}$ is the moduli space of stable $n$-punctured spheres, which is a complex manifold of complex dimension $n-3$. 
In particular, the moduli space is a complex manifold 
of {\it complex} dimension greater or equal to one so that, when the Fredholm index is assumed to be one, the actual dimension of the moduli space 
must be strictly larger than its virtual dimension expected by the Fredholm index. Note that this in turn implies that the moduli cannot be 
transversally cut out by the Cauchy-Riemann operator, in other words: {\it Even for generic choices of $J$, each moduli space of orbit curves 
with Fredholm index one must be nonregular in the sense that the the Cauchy-Riemann operator does not meet the zero section transversally.} \\
 
In order to see why the Fredholm index can be smaller than the actual dimension, observe that the index is sensitive to the underlying periodic orbit 
$\gamma$ and the dimension of $V$, while the actual dimension is not. On the other hand the nontrivial behaviour of the Conley-Zehnder index under replacing an 
orbit by some multiple cover makes it hard to exclude orbit curves with Fredholm index one. Restricting to contact homology for simplicity, note that 
the best way to get a hand on the possible range of the Fredholm index of orbit curves for the general case, i.e., without further assumptions on the 
underlying Reeb orbit $\gamma$, is to combine the formula for the virtual dimension of the moduli space $\IM_{0,0}(\gamma^{n-1};\gamma,...,\gamma)$, 
\begin{equation*} 
 \ind\CR_J = \mu_{CZ}(\gamma^{n-1}) - (n-1)\cdot \mu_{CZ}(\gamma) + (m-3) \cdot (2-n)
\end{equation*}
with the estimate for the Conley-Zehnder index of multiply covered orbits in [L], 
\begin{eqnarray*} 
 && (n-1)(\mu_{CZ}(\gamma)-(m-1))+(m-1)\;\leq\; \mu_{CZ}(\gamma^{n-1})\\&&\leq\; (n-1)(\mu_{CZ}(\gamma)+(m-1))-(m-1) 
\end{eqnarray*}
to obtain 
\begin{equation*} 
  (2-n)(2m-4)\,\leq\,\ind\CR_J\,\leq\,2n-4.
\end{equation*} 
While the right hand side agrees with the actual dimension of the moduli space $\IM_{0,0}(\gamma^{n-1};\gamma,...,\gamma)$ and 
is strictly greater than one, the left hand side is nonpositive for $m\geq 2$, i.e., $\dim V\geq 3$. 
Hence we cannot exclude branched covers of trivial cylinders with Fredholm index one for any number of punctures greater or equal to three as well as 
any dimension of $V$ greater or equal to three (without imposing further assumptions on the underlying Reeb orbit). \\    
\\
\subsection{Coherent compact perturbations} 
Since the actual dimension of the moduli spaces does not agree with the virtual dimension expected by the Fredholm index, we already deduced that  
the Cauchy-Riemann operator $\CR_J$ cannot be 
transversal to the zero section in a suitable Banach space bundle over a Banach manifold of maps. The general way to remedy this is to add 
compact perturbations to the Cauchy-Riemann operator $\CR_J$ so that it becomes transversal. Since the linearization of the perturbed Cauchy-Riemann 
operator then differs from the linearization of the original one only by a compact operator, it is still a Fredholm operator with the same 
index, which now by the implicit function theorem agrees with the local dimension of the zero set of the underlying nonlinear perturbed operator. 
In order to obtain a compactness 
result for this new zero set one also has to add compact perturbations to the Cauchy-Riemann operator 
over the moduli spaces forming the boundary. In particular, the compact perturbations chosen for any moduli space must be compatible 
with the compact perturbations chosen for the moduli spaces forming its boundary. 
The algebraic invariants are then defined by replacing the original compactified moduli space by the compactified zero set of the 
perturbed Cauchy-Riemann operator. Note that this can be achieved by either thinking about the specialities of the problem and then using special 
perturbations as in [F] or by building a general framework allowing for arbitrary compact perturbations.  
The observation that one is only interested in the zero set of the perturbed Cauchy-Riemann operator led 
to the (relative) virtual moduli cycle techniques in symplectic Floer homology and Gromov-Witten theory for general symplectic manifolds, see 
[LiuT], [LT], [FO], [MD], where the construction of the relative virtual moduli cycles in symplectic field theory is sketched in [B]. 
On the other hand, the wish to obtain the (relative) virtual moduli 
cycle directly as the zero set of the perturbed Cauchy-Riemann operator, viewed as a section in some kind of infinite-dimensional bundle over an 
infinite-dimensional space of maps, led to the invention of polyfolds by Hofer, Wysocki and Zehnder, see [HWZ] and the references therein. \\   

While the virtual moduli cycles techniques as well as the polyfold theory provide us with the correct setup to handle the problem of 
transversality in symplectic field theory, it seems that one has to give up any hope to finally compute the desired algebraic invariants. 
However it is is a folk's theorem in Gromov-Witten theory, see e.g. [MD], [MDSa], that in some good cases the situation can be drastically simplified: \\

Although the Cauchy-Riemann operator $\CR_J$ is not transversal to the zero section, it might happen that its zero set is still a manifold and that the 
virtual moduli cycle can be represented by the zero set of a generic section in a finite-dimensional obstruction bundle over the compactification 
of the nonregular moduli space. In particular, the zero set agrees with (the compactification of) the regular moduli space obtained by adding to the 
Cauchy-Riemann operator a suitably extension of the given obstruction bundle section. The standard example of such an obstruction bundle is the cokernel bundle, 
where one has to show that 
the cokernels of the linearization of $\CR_J$ at every zero always have the right dimension so that, in particular, they fit together to give a  
finite-dimensional vector bundle. Note however that the dimension of the cokernel usually jumps, so that the cokernels in general only fit together to 
local obstruction bundles, which leads to the definition of Kuranishi structures in [FO]. \\

Using the characterization of orbit curves as curves with trivial $\omega$-energy we can prove that we indeed have a global obstruction bundle 
over the compactification of every moduli space of orbit curves. 
While in Gromov-Witten theory the count of elements in the moduli space, more general, the cobordism class of the moduli 
space, is independent of the chosen abstract perturbation of the Cauchy-Riemann operator, this no longer holds for the moduli spaces in symplectic 
field theory. This follows from the fact that the moduli spaces in symplectic field theory typically have codimension one  
boundary strata, while in Gromov-Witten theory the regular moduli spaces form pseudo-cycles in the sense that the boundary strata have 
codimension at least two, i.e., from the homological point of view have no boundary. So while in Gromov-Witten 
theory the moduli spaces can be studied separately, the interplay between the different moduli spaces is the reason why  
the algebraic invariants of symplectic field theory are defined as differential algebras, which can be shown to be independent of extra choices 
like the cylindrical almost complex structure and the compact perturbation. In our case this problem is expressed by the fact that we have to study sections 
in vector bundles over moduli spaces with codimension one boundary, so that the count of zeroes in general depends on the choice of 
sections in the boundary, i.e., the chosen perturbations of the Cauchy-Riemann operator used to define the regular moduli spaces in the boundary. 
However we outline below that in our case we indeed have a well-defined count of zeroes so that, as in the Gromov-Witten case, we can (iteratively) define Euler 
numbers for our Fredholm problems. 

\section{Moduli space of orbit curves}

\subsection{Branched covers of trivial cylinders}

Choosing closed orbits $\gamma_{1,\pm}^{m^{\pm}_1},...,\gamma_{n^{\pm},\pm}^{m^{\pm}_{n^{\pm}}}$ 
of the vector field $R$ on $V$, where $\gamma^m$ denotes the $m$.th iterate of the simple orbit $\gamma$, and a homology class $A\in H_2(V)$, the moduli space 
$\IM_{A,0}(\gamma_{1,+}^{m^+_1},...,\gamma_{n^+,+}^{m^+_{n^+}};\gamma_{1,-}^{m^-_1},...,\gamma_{n^-,-}^{m^-_{n^-}})$ of punctured $J$-holomorphic curves in 
$\IR\times V$ of genus zero is defined as follows (see [EGH]): \\

Fix positive and negative punctures $z^{\pm}_1,...,z^{\pm}_{n^{\pm}} \in S^2$ and pairwise disjoint embeddings of half-cylinders 
$\psi^{\pm}_k:\IR^{\pm}\times S^1\hookrightarrow \Si$ with $\lim_{r\to\pm\infty}\psi^{\pm}_k(r,\cdot)=z^{\pm}_k$, where $\Si=S^2-\{z^{\pm}_1,...,z^{\pm}_{n^{\pm}}\}$.
Then the moduli space $\IM^0_{A,0}(\gamma_{1,+}^{m^+_1},...,\gamma_{n^+,+}^{m^+_{n^+}};\gamma_{1,-}^{m^-_1},...,\gamma_{n^-,-}^{m^-_{n^-}})$ of parametrized curves
consists of tuples $u=(u,j,\mu^{\pm})$, where $j$ denotes a complex structure on the punctured sphere $\Si$ which agrees with the standard complex structure 
on the cylindrical coordinate neighborhoods of the punctures, $\mu^{\pm}=(\mu^{\pm}_1,...,\mu^{\pm}_{n^{\pm}})$, 
$\mu^{\pm}_k \in (T_{z^{\pm}_k} S^2-\{0\})/\IR_+ \cong S^1$ 
is a collection of directions at the punctures $z^{\pm}_1,...,z^{\pm}_{n^{\pm}}$, called asymptotic markers, and 
$u:(\Si,j)\to(\IR\times V,J)$ is a $(j,J)$-holomorphic map which is asymptotically cylindrical over the closed orbit $\gamma_{k,\pm}^{m^{\pm}_k}$ at the 
puncture $z^{\pm}_k$,  
\begin{equation*} (u\circ\psi^{\pm}_k)(s,t+\mu^{\pm}_k) \to \gamma(m^{\pm}_k T^{\gamma_{\pm,k}} t), \;\;k=1,...,n^{\pm}. \end{equation*} 
Here $T^{\gamma}$ denotes the period of the simple orbit $\gamma$ and it follows from the chosen $S^1$-shift in the 
asymptotic condition that the asymptotic marker $\mu^{\pm}_k\in S^1$ is mapped to the point $z^{\gamma_{\pm,k}}=\gamma_{\pm,k}(0)$ on 
the underlying simple orbit. Note that when the asymptotic condition is fulfilled with the asymptotic marker $\mu^{\pm}_k$, then 
it also holds for the asymptotic markers $\mu^{\pm}_k + \ell/m^{\pm}_k$, $\ell=1,...,m^{\pm}_k-1$. Representing a basis of $H_1(V)$, which is assumed to be 
torsion-free as in [EGH], by  
circles in $V$ and choosing for each simple orbit $\gamma$ a spanning surface in $V$ between $\gamma$ and a suitable linear combination of these 
circles as in 0.3, one can assign an absolute homology class in $H_2(V)$ to each map $u$. With this we require that the map $u$ represents 
the given homology class $A\in H_2(V)$. \\

Note that when $n^++n^-\leq 3$ we have a unique complex structure $i$ on $\Si$ and we obtain the moduli 
space $\IM_{A,0}(\gamma_{1,+}^{m^+_1},...,\gamma_{n^+,+}^{m^+_{n^+}};\gamma_{1,-}^{m^-_1},...,\gamma_{n^-,-}^{m^-_{n^-}})$ as quotient of 
$\IM^0_{A,0}(\gamma_{1,+}^{m^+_1},...,\gamma_{n^+,+}^{m^+_{n^+}};\gamma_{1,-}^{m^-_1},...,\gamma_{n^-,-}^{m^-_{n^-}})$ under the obvious action of 
the automorphism group $\Aut(\Si,i)$. On the other hand, when $n^++n^-\geq 3$ the automorphism group of $(\Si,j)$ is trivial, so that the 
desired moduli space $\IM_{A,0}(\gamma_{1,+}^{m^+_1},...,\gamma_{n^+,+}^{m^+_{n^+}};\gamma_{1,-}^{m^-_1},...,\gamma_{n^-,-}^{m^-_{n^-}})$ agrees 
with the moduli space $\IM^0_{A,0}(\gamma_{1,+}^{m^+_1},...,\gamma_{n^+,+}^{m^+_{n^+}};\gamma_{1,-}^{m^-_1},...,\gamma_{n^-,-}^{m^-_{n^-}})$ of 
parametrized curves from before. \\

When all chosen simple orbits agree, $\gamma_{\pm,k}=\gamma$, $k=1,...,n^{\pm}$, and $A=0\in H_2(V)$, we already outlined in 0.2 that all curves have 
trivial $\omega$-energy $E_{\omega}(u)=0$, and therefore have $V$-image contained in a trajectory of the Reeb vector field. 
When there is at least one puncture it follows that the 
moduli space $\IM_{0,0}(\gamma^{m^+_1},...,\gamma^{m^+_{n^+}};\gamma^{m^-_1},...,\gamma^{m^-_{n^-}})$ entirely consists of branched covers of the trivial 
cylinder over a single closed orbit $\gamma$. For every (simple) closed orbit $\gamma$ of the vector field $R$, the trivial cylinder $\IR\times\gamma$ 
represents a curve in the above sense with $u_0: (\RS,i)\to(\IR\times V,J)$, 
$(s,t)\mapsto (T^{\gamma}s,\gamma(T^{\gamma}t))$, which is holomorphic by $J\del_s = \dot{\gamma}=R$. It follows that every curve $u$ in 
$\IM_{0,0}(\gamma^{m^+_1},...,\gamma^{m^+_{n^+}};\gamma^{m^-_1},...,\gamma^{m^-_{n^-}})$ is of the form $u=h\circ u_0$ 
with the branched covering map 
\begin{equation*}
 h: (\Si,j) \to \IR\times S^1
\end{equation*}
between the punctured Riemann spheres $(\Si,j)$ and $\RS\cong \IC^*=\CP-\{0,\infty\}$. \\

It directly follows from the asymptotic conditions for the curve 
$u$ in $\IM_{0,0}(\gamma^{m^+_1},...,\gamma^{m^+_{n^+}};\gamma^{m^-_1},...,\gamma^{m^-_{n^-}})$ that $h$ extends to a holomorphic map 
from $(S^2,j)\cong \CP$ to $\CP$. More precisely, it represents a meromorphic function $h$ on $(S^2,j)$, where the positive punctures $z^+_1,...,z^+_{n^+}$ 
are poles of order $m^+_1,...,m^+_{n^+}$, the negative punctures $z^-_1,...,z^-_{n^-}$ are zeroes of order $m^-_1,...,m^-_{n^-}$. For the rest of the paper we make the 
convention to identify $u$ directly with the branched covering $h$. Furthermore we make the convention 
that, unless otherwise mentioned, all considered branched covers are connected and have no nodes. 
Choosing the standard complex structure $i$ on $S^2$, $(S^2,i)=\CP$, and letting the positions of 
$z^{\pm}_1,...,z^{\pm}_{n^{\pm}}\in\CP$ vary, it follows that the moduli space $\IM_{0,0}(\gamma^{m^+_1},...,\gamma^{m^+_{n^+}};
\gamma^{m^-_1},...,\gamma^{m^-_{n^-}})$ agrees with the moduli space of meromorphic functions on $\CP$ with the given number of 
poles and zeroes with multiplicities $m^{\pm}_1,...,m^{\pm}_{n^{\pm}}$, where we just must take care of the possible 
different choices for the asymptotic markers.\\

For the following expositions we assume that $m^+_1+...+m^+_{n^+} = m^-_1+...+m^-_{n^-}$ since else the 
moduli space is obviously empty by homological reasons. In particular, there are no holomorphic planes ($n=n^++n^-=1$). For $n=2$ the moduli 
space $\IM_{0,0}(\gamma^m;\gamma^m)/\IR$ consists precisely of $m^2$ elements, namely the unique trivial cylinder over the iterated orbit 
$\gamma^m$ together with the $m^2$ possible choices for the asymptotic marker above and below. 
Note that here the actual and the virtual dimension given by the Fredholm index agree to be zero, so that they are not 
interesting from the viewpoint of symplectic field theory. Hence it suffices to restrict our considerations to the stable case $n\geq 3$. \\ 
\\
{\bf Proposition 1.1:} {\it For $n=n^++n^- \geq 3$ the moduli space of orbit curves (connected, without nodes) with 
fixed multiplicities $m^{\pm}_1,...,m^{\pm}_{n^{\pm}}$ is given by}
\begin{equation*}
 \IM_{0,0}(\gamma^{m^+_1},...,\gamma^{m^+_{n^+}};\gamma^{m^-_1},...,\gamma^{m^-_{n^-}})/\IR 
  \cong S^1 \times \IM_{0,n^++n^-} \times \IZ_{m^+} \times \IZ_{m^-},
\end{equation*}
{\it where $\IM_{0,n}$ denotes the moduli space of $n$-punctured spheres and $m^{\pm}=m^{\pm}_1\cdot ...\cdot m^{\pm}_{n^{\pm}}$.} \\  
\\
{\it Proof:} 
For the proof we fix the natural complex structure $j=i$ on $S^2$, $(S^2,i)=\CP$, and let instead the positions of the punctures $z_1^{\pm},...,z_{n^{\pm}}^{\pm}$ 
vary. Since the zeroth Picard group $\Pic^0(\CP)$ is trivial, 
i.e., all degree zero divisors on $\CP$ are in fact principal divisors, it follows that a meromorphic function exists 
for any choice of zeroes and poles with multiplicities, as long as the number of poles with multiplicities agrees with the number 
of zeroes with multiplicities. More explicitly, an example of $h$ is 
\begin{equation*} 
        h^0(z) = \frac{\prod_{k=1}^{n^-} (z-z^-_k)^{m^-_k}}{\prod_{k=1}^{n^+} (z-z^+_k)^{m^+_k}}
\end{equation*} 
and it follows from Liouville's theorem that such a map is uniquely determined up to a nonzero 
multiplikative factor, i.e., $h = a \cdot h_0$ with $a \in \IC^*$. Since for $n\geq 3$ the automorphism group $\Aut(\CP)$ already acts freely 
on the ordered set of punctures $(z_1^{\pm},...,z_{n^{\pm}}^{\pm})$, it follows that the moduli space   
agrees with the product $\IC^*\times\IM_{0,n^++n^-}$ with $\IC^*\cong\RS$. 
On the other hand there are $m^{\pm}_k$ possible directions for the asymptotic marker $\mu^{\pm}_k$ 
at each puncture $z^{\pm}_k$, $k=1,...,n^{\pm}$, for each $(h,j)$ as outlined in the definition of the moduli spaces, so that $\mu^{\pm}_k\in\IZ_{m^{\pm}_k}$, 
i.e., $\mu^{\pm}=(\mu^{\pm}_1,...,\mu^{\pm}_{n^{\pm}})\in \IZ_{m_1^{\pm}}\times ... \times \IZ_{m_{n^{\pm}}^{\pm}}\cong\IZ_{m^{\pm}}$. $\qed$ \\

In what follows we fix the multiplicities $m_1^{\pm},...,m_{n^{\pm}}^{\pm}$ and abbreviate the corresponding moduli space of orbit curves by 
\begin{equation*}
\IM= \IM_{0,0}(\gamma^{m^+_1},...,\gamma^{m^+_{n^+}};\gamma^{m^-_1},...,\gamma^{m^-_{n^-}})/\IR. 
\end{equation*}
Note that here we view the target $\RS$ as a cylindrical cobordism in the sense of [BEHWZ], so that we quotient out the corresponding $\IR$-symmetry on 
the moduli space.  
Later, for the proof of the main theorem, we further have to consider the corresponding moduli space of holomorphic curves in $\RS$ 
without quotiening out the $\IR$-translations,
\begin{equation*}
\IM^0 = \IM_{0,0}(\gamma^{m^+_1},...,\gamma^{m^+_{n^+}};\gamma^{m^-_1},...,\gamma^{m^-_{n^-}}), 
\end{equation*}
i.e., we view the holomorphic curves as sitting in a noncylindrical cobordism by just ignoring the natural $\IR$-action.
 
\subsection{Compactification}

While introducing abstract perturbations we must asure that these are compatible with the curve splitting    
phenomena described in the compactness theorem of symplectic field theory. Hence we must also include the compactification of the moduli space of orbit curves  
into our considerations which is, of course, not too bad. Recall that by [BEHWZ] the compactification of a moduli space of curves 
in a cylindrical or noncylindrical cobordism consists of holomorphic curves in cobordisms together with a level structure. Calling 
a level (non-)cylindrical whenever the corresponding cobordism is (non-)cylindrical, observe that when we start with curves in a cylindrical cobordism 
the resulting levels are all cylindrical. On the other hand, when we start with curves in a noncylindrical cobordism, there is precisely one noncylindrical level, while 
all other levels are cylindrical. Furthermore we call a connected component of a holomorphic curve (non-)cylindrical when it is (not) a cylinder.  
This leads to the following compactness statement: \\
\\
{\bf Proposition 1.2:} {\it The boundary of the compactified moduli space $\CM$ consists of level holomorphic curves in the 
sense of [BEHWZ], which are connected or disconnected nodal branched covers of the same orbit cylinder, such that the punctured spheres underlying all 
noncylindrical components are stable and on each level there is at least one noncylindrical component. The same holds true for the 
compactification $\overline{\IM^0}$, except that the last part of the statement need not be satisfied for the noncylindrical level. 
For $\CM$ it follows that all connected components carry strictly less than $n$ punctures, whereas for $\overline{\IM^0}$ 
this is true only up to the case of a two level curve where all curves on the noncylindrical level are cylinders. } \\
\\  
{\it Proof:} Choosing a sequence of holomorphic curves in $\IM$, it follows from the compactness theorem in [BEHWZ] that 
a suitable subsequence converges 
to a level holomorphic map in the sense of [BEHWZ]. It follows from lemma 5.4 in [BEHWZ] together with the preservation of the 
$\omega$-energy that the connected components in each level of the limiting level curve  
are again, after resolving the nodes, multiple covers of the corresponding orbit cylinder. Since there are no multiple covers 
with one puncture and every curve with no punctures is constant it follows that every component of the limit level holomorphic map 
has at least two punctures, i.e., that every noncylindrical component has positive Euler characteristic. Furthermore there always must be a 
noncylindrical component on each cylindrical level, since otherwise the $\IR$-action is trivial. The remaining statements on the 
number of punctures follow from the additivity of the Euler characteristic. $\qed$ \\ 
\\
{\bf Definition 1.3:} {\it A $(n^+,n^-)$-labelled tree with level structure is a tuple $(T,\LL)=(T,E,\Lambda^+,\Lambda^-,\LL)$, where $(T,E)$ is a tree 
with the set of vertices $T$ and the edge relation $E\subset T\times T$, the sets $\Lambda^{\pm} = 
(\Lambda^{\pm}_{\alpha})_{\alpha\in T}$ are decompositions of $\{1,...,n^{\pm}\}$, i.e.,}  
\begin{equation*} 
 \bigcup_{\alpha\in T} \Lambda^{\pm}_{\alpha} = \{1,...,n^{\pm}\},\,\,
 \Lambda^{\pm}_{\alpha}\cap\Lambda^{\pm}_{\beta}=\emptyset \,\textrm{when}\,\alpha\neq\beta,  
\end{equation*} 
{\it and $\LL: T \to \{1,...,L\}$ is surjective map, which is called a level structure. Furthermore, a tuple $(T,\LL,\ell_0)=(T,E,\Lambda^+,\Lambda^-,\LL,\ell_0)$ 
with $\ell_0\in\{1,...,L\}$ is called a $(n^+,n^-)$-labelled tree with based level structure. } \\

Observe that every level branched cover in $\CM$ represents a $(n^+,n^-)$-labelled tree with level structure, where the tree structure 
$(T,E)$ represents the underlying nodal curve, i.e., bubble tree, and the elements $k\in\{1,...,n^{\pm}\}$ represent 
positive or negative punctures. On the other hand, a level branched cover in the boundary of $\IM^0$ represents a tree with based level structure 
$(T,\LL,\ell_0)$ with $\ell_0$ denoting the noncylindrical level. It follows that $\CM$ and $\overline{\IM^0}$ carry natural stratifications
\begin{equation*}
\CM=\bigcup_{T,\LL} \IM_{T,\LL},\;\; \overline{\IM^0}=\bigcup_{T,\LL,\ell_0} \IM^0_{T,\LL,\ell_0}
\end{equation*}
where $\IM_{T,\LL}$ and $\IM^0_{T,\LL,\ell_0}$ can be described as follows: \\

First we can assign to every labelled tree with level structure $(T,\LL)=(T,E,\Lambda^{\pm},\LL)$ a nodal surface with positive and negative punctures 
by assigning to each $\alpha\in T$ a sphere $S_{\alpha}=S^2$, to any edge $(\alpha,\beta)\in E$ a marked point $z_{\alpha\beta}\in S_{\alpha}$ 
and to any $k\in\Lambda^{\pm}_{\alpha}$, $\alpha\in T$ a positive, respectively negative puncture $z^{\pm}_k \in S_{\alpha}$. 
Since to each positive, respectively negative puncture we assign a fixed multiple $\gamma^{m^{\pm}_k}$ of the underlying simple orbit $\gamma$, we can naturally 
assign a multiplicity with sign $m_{\alpha\beta}\in \IZ$ to each edge in $E$ by requiring for each 
$\alpha\in T$ that
\begin{equation*} 
 \sum_{\beta: \alpha E\beta} m_{\alpha\beta} + \sum_{k\in\Lambda^+_{\alpha}} m^+_k - \sum_{k\in\Lambda^-_{\alpha}} m^-_k = 0. 
\end{equation*}
Note that each edge $(\alpha,\beta)$ with $m_{\alpha\beta}\neq 0$ corresponds to a positive or negative puncture 
for the components $\alpha$ and $\beta$ and $m_{\alpha\beta}=-m_{\beta\alpha}$ denotes the period with sign. In particular, when $m_{\alpha\beta}>0$ 
then $\LL(\alpha)>\LL(\beta)$, whereas by similar arguments the edges with trivial multiplicity $m_{\alpha\beta}=0$ corresponds to  
nodes between components $\alpha$ and $\beta$ in the same level, $\LL(\alpha)=\LL(\beta)$. With 
this we define sets of positive, respectively negative punctures on $S_{\alpha}$ by 
\begin{eqnarray*}
  Z^+_{\alpha} &=& \{z^+_k: k\in\Lambda^+_{\alpha}\}\cup\{z_{\alpha\beta}: \LL(\beta)>\LL(\alpha)\} \\
                   &=& \{z^+_{\alpha,k}: k=1,...,n^+_{\alpha}\}, \\
  Z^-_{\alpha} &=& \{z^-_k: k\in\Lambda^-_{\alpha}\}\cup\{z_{\alpha\beta}: \LL(\beta)<\LL(\alpha)\} \\
                   &=& \{z^-_{\alpha,k}: k=1,...,n^-_{\alpha}\}
\end{eqnarray*}
and denote the corresponding punctured sphere by $\Si_{\alpha}=S_{\alpha}-\{z^{\pm}_{\alpha,1},...,z^{\pm}_{\alpha,n^{\pm}_{\alpha}}\}$, 
while $Z^0_{\alpha}=\{z_{\alpha\beta}: \LL(\alpha)=\LL(\beta)\}$ is the set of nodes connecting $S_{\alpha}$ with $S_{\beta}$ of the 
same level. Note that by the above definitions we assign a positive multiplicity $m^{\pm}_{\alpha,k}$ to any point $z^{\pm}_{\alpha,k}$ in $Z^{\pm}_{\alpha}$. 
Finally note that we did not fix the complex structure on any of the punctured spheres $S_{\alpha}$. \\

We want to describe the moduli space $\IM_{T,\LL}$ using the corresponding moduli spaces of nodal curves on the different levels. 
For this observe that to any labelled tree with level structure $(T,E,\Lambda^{\pm},\LL)$ we can assign a 
tuple of labelled trees $T_{\ell}=(T_{\ell},E_{\ell},\Lambda^{\pm}_{\ell})$, $\ell\in\{1,...,L\}$, where 
$T_{\ell}=\{\alpha\in T: \LL(\alpha)=\ell\}$, $E_{\ell}=E\cap(T_{\ell}\times T_{\ell})$ and 
$\Lambda^{\pm}_{\ell}=(\Lambda^{\pm}_{\ell,\alpha})_{\alpha\in T_{\ell}}$ with 
$\Lambda^{\pm}_{\ell,\alpha}=\Lambda^{\pm}_{\alpha} \cup\{\beta \in T_{\ell\pm 1}: \alpha E \beta\}$. \\
 
For every $T_{\ell}=(T_{\ell},E_{\ell},\Lambda^{\pm}_{\ell})$, $\ell\in\{1,...,L\}$ we now introduce the moduli space $\IM_{T_{\ell}}$ as follows: Every element 
in $\IM_{T_{\ell}}$ is a tuple $(h_{\ell},j_{\ell},\mu_{\ell}^{\pm})=(h_{\alpha},j_{\alpha},\mu_{\alpha}^{\pm})_{\alpha\in T_{\ell}}$, where 
$j_{\alpha}$ is a complex structure on $\Si_{\alpha}$ and $h_{\alpha}: (\Si_{\alpha},j_{\alpha})\to\RS$ extends to a meromorphic function on 
$(S_{\alpha}=S^2,j_{\alpha})$ with poles, respectively zeroes $z^{\pm}_{\alpha,1},...,
z^{\pm}_{\alpha,n^{\pm}_{\alpha}}$ of multiplicities $m^{\pm}_{\alpha,1},...,m^{\pm}_{\alpha,n^{\pm}_{\alpha}}$, such that 
$h_{\alpha}(z_{\alpha\beta})=h_{\beta}(z_{\beta\alpha})$ if $z_{\alpha\beta}\in Z^0_{\alpha}$, i.e., $z_{\beta\alpha}\in Z^0_{\beta}$.
Further $\mu_{\alpha}^{\pm}=(\mu^{\pm}_{\alpha,1},...,\mu^{\pm}_{\alpha,n^{\pm}_{\alpha}})$ denotes the collection of asymptotic markers 
$\mu^{\pm}_{\alpha,k}\in\IZ_{m^{\pm}_{\alpha,k}}$. \\

Note that in general the trees $T_{\ell}$ are not connected. Denoting the connected components by $T_{\ell,1},...,T_{\ell,N_{\ell}}$, 
the moduli space $\IM_{T_{\ell}}$ can be written as direct product
\begin{equation*} 
 \IM_{T_{\ell}} = \IM_{T_{\ell,1}}\times ... \times \IM_{T_{\ell,N_{\ell}}} \times \IR^{N_{\ell}-1}
\end{equation*} 
of moduli spaces $\IM_{T_{\ell,k}}$, $k=1,...,N_{\ell}$ of connected nodal branched covers, where the $\IR$-factors 
encode the relative $\IR$-position of the connected components of the curves in $\IM_{T_{\ell}}$. \\
 
With the moduli spaces $\IM_{T_1},...,\IM_{T_L}$ we can finally describe the moduli spaces $\IM_{T,\LL}$ and $\IM^0_{T,\LL,\ell_0}$: \\

While the definitions of complex structures and holomorphic maps is straightforward, we explicitly want that two tuples 
$(h_{\ell},j_{\ell},\mu_{\ell})_{\ell=1,...,L}$ represent the same element in $\IM_{T,\LL}$ if the asymptotic markers at pairs 
of positive and negative punctures, which correspond to edges between components in neighboring levels, describe the same decorations. Note that 
this convention is implicit in the proof of the master equation 
of (rational) symplectic field theory, which is derived by studying the codimension boundary strata of moduli spaces. Indeed we will show 
below that this convention guarantees that the compactified moduli space $\CM$ (and $\overline{\IM^0}$) carries the 
structure of a manifold with boundary. Going back to the goal of describing $\IM_{T,\LL}$ explicitly, we assign to any tuple 
$(h_{\ell},j_{\ell},\mu_{\ell}^{\pm})_{\ell=1,...,L}\in\IM_{T_1}\times ...\times\IM_{T_L}$ a tuple $(h,j,\mu^{\pm},\theta)\in\IM_{T,\LL}$, where 
$(h,j)=(h_{\ell},j_{\ell})_{\ell=1,...,L}=(h_{\alpha},j_{\alpha})_{\alpha\in T}$. For the asymptotic markers $\mu^{\pm}$ and decorations $\theta$ we 
recall that 
\begin{eqnarray*} 
 \mu_{\ell}^{\pm}=(\mu_{\alpha}^{\pm})_{\alpha\in T_{\ell}}, 
 &&\mu^+_{\alpha}=((\mu^+_k)_{k\in\Lambda^+_{\alpha}},(\mu_{\alpha\beta})_{\LL(\beta)>\LL(\alpha)}),\\
 &&\mu^-_{\alpha}=((\mu^-_k)_{k\in\Lambda^-_{\alpha}},(\mu_{\alpha\beta})_{\LL(\beta)<\LL(\alpha)}).
\end{eqnarray*}
>From this we get asymptotic markers $\mu^{\pm}=(\mu^{\pm}_k)_{k=1,...,n^{\pm}}$ and decorations $\theta=(\theta_{\alpha\beta})_{\LL(\alpha)>\LL(\beta)}$ by setting 
\begin{equation*}
   \theta_{\alpha\beta}=[(\mu_{\alpha\beta},\mu_{\beta\alpha})]\in\frac{\IZ_{|m_{\alpha\beta}|}\times\IZ_{|m_{\alpha\beta}|}}{\Delta_{\alpha\beta}},
\end{equation*}
where $\Delta_{\alpha\beta}=\Delta_{\beta\alpha}$ denotes the diagonal in $\IZ_{|m_{\alpha\beta}|}\times\IZ_{|m_{\beta\alpha}|}$.
For this recall that $m_{\alpha\beta}=-m_{\beta\alpha}$ and observe that two pairs of asymptotic markers $(\mu_{\alpha\beta},\mu_{\beta\alpha})$ 
and $(\mu'_{\alpha\beta},\mu'_{\beta\alpha})$ represent the same decoration if there exists some $\mu_0 \in \IZ_{|m_{\alpha\beta}|}$ with 
$(\mu'_{\alpha,\beta},\mu'_{\beta,\alpha})=(\mu_{\alpha\beta}+\mu_0,\mu_{\beta\alpha}+\mu_0)$. 
With this it follows that the moduli space $\IM_{T,\LL}$ is given by 
\begin{equation*} \IM_{T,\LL}=\frac{\IM_{T_1}\times ... \times \IM_{T_L}}{\Delta}.
\end{equation*}
with $\Delta=\prod_{\LL(\alpha)>\LL(\beta)} \Delta_{\alpha\beta}$. On the other hand, it follows from the same arguments that 
$\IM^0_{T,\LL,\ell_0}$ is given by 
\begin{equation*} \IM^0_{T,\LL,\ell_0}=\frac{\IM_{T_1}\times ... \times \IM^0_{T_{\ell_0}}\times ... \times \IM_{T_L}}{\Delta},
\end{equation*}
Here $\IM^0_{T_{\ell_0}}$ is the moduli space of orbit curves on the noncylindrical level, so that $\IM^0_{T_{\ell_0}}=\IR\times\IM_{T_{\ell_0}}$ 
whenever $T_{\ell_0}$ represents a curve with at least one noncylindrical component, and just consists of a point if all components are trivial cylinders. \\

Observe that each $\IM_{T,\LL}$ is a smooth manifold of codimension 
\begin{equation*} \dim \IM - \dim \IM_{T,\LL} = L - 1 + 2 N , \end{equation*}
where $L$ is the number of levels and $N=\frac{1}{2}\sharp\{\alpha E\beta: \LL(\alpha)=\LL(\beta)\}$ denotes the number of nodes between 
components in the same level. For this observe that creating a new level we indeed only loose one dimension corresponding to the 
$\IR$-coordinate on the new level which is quotiented out. It follows that the compactified moduli space $\CM$ is a stratified space with 
natural stratification 
\begin{equation*} 
    \IM=\CM^0 \subset \CM^1 \subset \CM^2 \subset ... \subset \CM^k \subset ... \subset \CM^{\infty}=\CM,
\end{equation*} 
where 
\begin{equation*} \CM^k=\bigcup_{(T,\LL):L-1+2N\leq k} \IM_{T,\LL}. \end{equation*} 
contains the components of the compactified moduli space of codimension at most $k$.
In the same way we have 
\begin{equation*} 
    \IM^0=\overline{\IM^0}^0 \subset \overline{\IM^0}^1 \subset \overline{\IM^0}^2 \subset ... \subset \overline{\IM^0}^k \subset ... 
    \subset \overline{\IM^0}^{\infty}=\overline{\IM^0},
\end{equation*} 
where 
\begin{equation*} \overline{\IM^0}^k=\bigcup_{(T,\LL,\ell_0):L-1+2N\leq k} \IM^0_{T,\LL,\ell_0}. \end{equation*} 

Observe that $\CM^1$, defined as disjoint union of the 
moduli space with the codimension one boundary components, consists of curves with two level and no nodes. 
More precisely, the connected components of the codimension one boundary are given by fibre products 
\begin{equation*}
\IM_1 \times_{\IZ_{m_{1,2}}} \IM_2 = \frac{\IM_1\times\IM_2}{\Delta}, 
\end{equation*}
where $\IM_1=\IM_{T_1}$, $\IM_2=\IM_{T_2}$ denote moduli spaces of possibly disconnected branched covers without nodes. Note that here   
$T_1,T_2$ are trees with trivial edge relation and $\IZ_{m_{1,2}}=\prod_{\LL(\alpha)=2,\LL(\beta)=1} \IZ_{|m_{\alpha\beta}|}$ acts on 
$\IM_1$ and $\IM_2$ in the obvious way. On the other hand, observe that the connected components of the codimension one boundary of $\IM^0$ are 
given either given by products of the form 
\begin{equation*}
\IM^0_1 \times_{\IZ_{m_{1,2}}} \IM_2,\;\; \IM_1 \times_{\IZ_{m_{1,2}}} \IM^0_2 
\end{equation*}
with $\IM^0_1=\IR\times\IM_1$ and $\IM^0_2=\IR\times\IM_2$ or 
\begin{equation*}
\{\point\}\times \IM,\;\; \IM\times\{\point\}
\end{equation*} 
corresponding to $\IM^0_1=\{\point\}$, $\IM^0_2=\{\point\}$, respectively, i.e., where on the noncylindrical level we just find trivial cylinders. \\

We close this section with an important technical lemma about the compactified moduli spaces $\CM$ and $\overline{\IM^0}$. \\
\\
{\bf Proposition 1.4:} {\it The compactified moduli spaces $\CM$ and $\overline{\IM^0}$ naturally carry the structure of a manifold with corners.} \\
\\
{\it Proof:} We prove the statement only for the compactification of $\IM$, since the statement about the compactification of $\IM^0$ follows the same 
arguments. Essentially it follows from an explicit description of the moduli space $\IM$ and its compactification in terms of Fenchel-Nielsen 
coordinates: \\

Recall from the definition of the moduli spaces that we fixed $n^+$ positive and $n^-$ negative punctures $z^{\pm}_1,...,z^{\pm}_{n^{\pm}}\in S^2$ and 
fixed cylindrical coordinates 
\begin{equation*} \psi^{\pm}_k: \IR^{\pm}_0\times S^1\hookrightarrow\Si \end{equation*} 
around each puncture $z^{\pm}_k$, $k\in\{1,...,n^{\pm}\}$ on the punctured sphere $\Si=S^2-\{z^{\pm}_1,...,z^{\pm}_{n^{\pm}}\}$. Beside the 
mentioned embeddings of half-cylinders we now embed $n-3$ finite cylinders 
$\psi_k: [-1,+1] \times S^1 \hookrightarrow \Si$, $k\in\{1,...,n-3\}$ such that their images are pairwise disjoint, 
disjoint from the cylindrical coordinate neighborhoods of the punctures and such that the circles $\psi_k(\{0\}\times S^1) \subset \Si$, $k=\in\{1,...,n-3\}$ 
define a pair of pants decomposition of $\Si$. 
Observe that this naturally defines a $(n^+,n^-)$-labelled tree $(T^0,E^0,\Lambda^0)$, where $T^0$ is the set of pair-of-pants components, 
\begin{equation*} \Si=\bigcup_{\alpha\in T^0} Y_{\alpha} \end{equation*} 
with the obvious edge relation 
\begin{equation*} (\alpha,\beta)\in E^0\;\;\Leftrightarrow\;\; Y_{\alpha}\cap Y_{\beta}\neq\emptyset, \end{equation*}
and the decompositions $\Lambda^{0,\pm}=(\Lambda^{0,\pm}_{\alpha})_{\alpha\in T^0}$ of the sets $\{1,...,n^{\pm}\}$ given by  
\begin{equation*} k\in\Lambda^{0,\pm}_{\alpha}\subset\{1,...,n^{\pm}\}\;\;\Leftrightarrow\;\;z^{\pm}_k\in Y_{\alpha}. \end{equation*}
We fix a complex structure $j_0$ on $\Si$ such that it agrees with the natural complex structures on the embedded cylinders. 
Let $\bar{E}^0=E^0/\{(\alpha,\beta)\sim(\beta,\alpha)\}$ be the set of undirected edges and for 
every $\tau\in \bar{E}^0$ let $\psi_{\tau}: [-1,+1]\times S^1 \hookrightarrow\Si$ denote the corresponding embedding of the finite cylinder. 
For every $(r_{\tau},t_{\tau})\in (\IR^+_0\times S^1)^{\bar{E}^0}$ let $\Si_{(r_{\tau},t_{\tau})}$ denote the punctured Riemann sphere obtained from 
$\Si$ by replacing for each $\tau\in \bar{E}^0$ the embedded cylinders $\psi_{\tau}([-1,0]\times S^1)$ by 
$[-r_{\tau},0]\times S^1$, $\psi_{\tau}([0,+1]\times S^1)$ by $[0,+r_{\tau}]\times S^1$, and gluing 
$[-r_{\tau},0]\times S^1$ and $[0,+r_{\tau}]\times S^1$ with a twist $t_{\tau}\in S^1$. Note that for any $(r_{\tau},t_{\tau})\in (\IR^+_0\times S^1)^{\bar{E}^0}$ 
the punctured Riemann sphere $\Si_{(r_{\tau},t_{\tau})}$ represents an element in $\IM_{0,n}$ and we assume without loss of generality 
that the complex structure $j_0$ on the noncylindrical part of $\Si$ is chosen such that the map from $(\IR^+_0 \times S^1)^{\bar{E}^0}$ to $\IM_{0,n}$ 
is indeed a coordinate chart for $\IM_{0,n}$. \\

Assuming that we have covered $\IM_{0,n}$ by coordinate charts of the above form, we are now ready to describe the compactification 
$\CM$ of $\IM$ by compactifying each coordinate neighborhood in the following nonstandard way. First observe (compare [BEHWZ]) that when we compactify each coordinate 
neighborhood by viewing it as a subset of $(\RS)^{\bar{E}^0}$ with  
compactification $(\overline{\IR}\times S^1)^{\bar{E}^0}$, $\overline{\IR}=\IR\cup\{\pm\infty\}$, then we obtain the Deligne-Mumford compactification 
$\overline{\IM}^\$_{0,n}$ with decorations at each node. On the other hand, note that when we use the compactification $(\CP)^{\bar{E}^0}$ of  
$(\RS)^{\bar{E}^0}$ by identifying $\RS\cong\IC^*$, then we obtain the 
usual Deligne-Mumford compactification $\CM_{0,n}$ without decorations. In order to obtain $\CM=S^1\times\TiM_{0,n}\times\IZ_{m^+}\times\IZ_{m^-}$ we 
need yet another compactification $\TiM_{0,n}$ of $\IM_{0,n}$. Besides that we want decorations only at those nodes which correspond to a pair of a positive and 
a negative puncture, we must keep track of the relative $\IR$-shift of the different components when they are mapped to the trivial cylinder. \\
 
To this end, recall that each $k\in\Lambda^{0,\pm}_{\alpha}$ represents a positive, respectively negative puncture to which 
we assign a fixed multiple $\gamma^{m^{\pm}_k}$ of the underlying simple orbit $\gamma$. Hence we can again naturally 
assign a multiplicity with sign $m_{\alpha\beta}\in \IZ$ to each directed edge in $E^0$ by requiring for each 
$\alpha\in T^0$ that
\begin{equation*} 
 \sum_{\beta: \alpha E^0\beta} m_{\alpha\beta} + \sum_{k\in\Lambda^{0,+}_{\alpha}} m^+_k - \sum_{k\in\Lambda^{-,0}_{\alpha}} m^-_k = 0. 
\end{equation*}
Note that $m_{\beta\alpha}=-m_{\alpha\beta}$. Now we identify the coordinate subset of $\IM_{0,n}$ not with 
$(\IR^+_0 \times S^1)^{\bar{E}^0}$, but view it as a linear subspace of 
$(\IR^+_0 \times S^1)^{\bar{E}^0}\times \IR^{T^0\times T^0}$ by setting for $(\alpha,\beta)\in T^0\times T^0$  
\begin{equation*} s_{\alpha\beta} = \sum_{i=1}^k m_{(\gamma_{i-1},\gamma_i)} r_{[\gamma_{i-1},\gamma_i]}, \end{equation*} 
where $\alpha=\gamma_0,...,\gamma_k=\beta$ is the enumeration of vertices on the unique directed path in $(T^0,E^0)$ from $\alpha$ to $\beta$. \\

Distinguishing further the undirected edges in $\bar{E}^0$ by whether their multiplicity is zero or not, 
$\bar{E}^0=\bar{E}^0_0\cup \bar{E}^0_{\pm}$, we now obtain $\TiM_{0,n}$ by viewing it as a subset of   
$(\IR\times S^1)^{\bar{E}^0_0}\times(\IR\times S^1)^{\bar{E}^0_{\pm}} \times \IR^{T^0\times T^0}$ with compactification given by 
$(\CP)^{\bar{E}^0_0}\times (\overline{\IR}\times S^1)^{\bar{E}^0_{\pm}} \times\overline{\IR}^{T^0\times T^0}$. 
It directly follows from the construction of $\TiM_{0,n}$ that $\TiM_{0,n}$ carries the structure of a manifold 
with corners. Further the boundary of $\IM_{0,n}$ in $\TiM_{0,n}$ consists of tuples 
$((r_{\tau},t_{\tau}),(s_{\alpha\beta}))$ with $r_{\tau}=\infty$ for some edge $\tau\in \bar{E}^0$.  While the coordinates $(r_{\tau},t_{\tau})$ describe a 
nodal curve with decorations at nodes corresponding to edges in $\bar{E}^0_{\pm}$, we show that the coordinates $(s_{\alpha\beta})$ describes a 
level structure with relative $\IR$-shifts. More precisely, recalling that $\IM\cong S^1 \times \IM_{0,n}\times \IZ_{m^+}\times \IZ_{m^-}$, 
we show in the following that there is a natural identification of $S^1 \times \TiM_{0,n}\times\IZ_{m^+}\times \IZ_{m^-}$ with the compactified 
moduli space $\CM$ of orbit curves. To this end we assign to any tuple $(t_0,((r_{\tau},t_{\tau}),(s_{\alpha\beta})),\mu^{\pm})$ a level branched covering 
$(h,j,\mu^{\pm},\theta)$ as follows: \\

First observe that the underlying nodal curve is described by the coordinates $(r_{\tau},t_{\tau})\in (\CP)^{\bar{E}^0_0}\times 
(\overline{\IR}\times S^1)^{\bar{E}^0_{\pm}}$, where $\alpha,\beta\in T^0$ belong to the same connected component when 
$r_{\tau}<\infty$ for each edge on the unique path from $\alpha$ to $\beta$. Note that the latter defines an equivalence 
relation $\approx$ on $T^0$, such that the quotient $T=T^0/\approx$ with induced edge relation $E\subset T\times T$ is the tree representing the nodal curve. 
Distinguishing the undirected edges in $\bar{E}$ by whether they have a nonzero multiplicity or not, $\bar{E}=\bar{E}_0\cup\bar{E}_{\pm}$, note 
that the undirected edges in $\bar{E}_0$ now correspond to nodes connecting components in the same level, while the edges in $\bar{E}_{\pm}$ 
correspond to pairs of components living on neighboring levels connected by a positive, respectively negative puncture. 
Since each branched cover of the trivial cylinder is determined up to $\RS$-shift by the underlying punctured sphere in 
$\IM_{0,n}$, it follows that the level branched cover in $\CM$ is already known up to the $S^1$-shifts, decorations in 
$\IZ_{|m_{\alpha\beta}|}\times \IZ_{|m_{\alpha\beta}|}/\Delta\cong \IZ_{|m_{\alpha\beta}|}$ at the punctures 
between levels and the level structure with the relative $\IR$-shifts. \\ 

First, in order to see how the coordinates $s_{\alpha\beta}\in\bar{\IR}$, $\alpha,\beta\in T^0$ fix the level structure and 
the relative $\IR$-shifts, let $((r_{\tau}^n,t_{\tau}^n),(s_{\alpha\beta}^n))\in (\IR^+_0\times S^1)^{E^0}\times 
\IR^{T^0\times T^0}$ be a sequence converging to $((r_{\tau},t_{\tau}),(s_{\alpha\beta}))$, where without loss 
of generality $t_{\tau}^n=t_{\tau}$. 
Let $\Si_n = \Si_{(r^n_{\tau},t^n_{\tau})}$ be the corresponding 
sequence of punctured spheres converging to the punctured nodal surface $\Si$ with connected components 
$\Si_{[\alpha]}$, $[\alpha]\in T=T^0/\approx$, and let $h^n = (h^n_1,h^n_2): \Si_n \to \RS$ be a corresponding 
sequence of branched covering maps converging to the level branched cover $h = (h^{[\alpha]})_{[\alpha]\in T}: 
\Si\to\RS$. In order 
to see the relation between $(s_{\alpha\beta}^n)_{\alpha,\beta}$ and the level structure and relative $\IR$-shifts of 
the limit level curve $h$, fix points $z_{\alpha},z_{\beta}$ on the pair of pants components corresponding to two 
chosen $\alpha,\beta \in T^0$. For each $(\gamma,\delta)\in E^0$ on the unique path from $\alpha$ to $\beta$, 
set $h^n_{\gamma\delta} = \int_0^1 h^n_1\circ\psi^n_{\gamma\delta}(r^n_{\gamma\delta},t) dt$, with the embedding 
$\psi^n_{\gamma\delta}: [-r^n_{\gamma\delta},+r^n_{\gamma\delta}]\times S^1\to \Si_n$ of the 
finite cylinder at the edge $(\gamma,\delta)\in E^0$, where $r^n_{\gamma\delta}=r^n{[\gamma,\delta]}$ and 
$\psi^n_{\delta\gamma}:[-r^n_{\delta\gamma},r^n_{\delta\gamma}]\times S^1\to \Si_n$, 
$\psi^n_{\delta\gamma}(s,t)=\psi^n_{\gamma\delta}(-s,-t)$. Observe that we have 
\begin{eqnarray*} 
h^n_{\gamma\delta} - h^n_{\delta\gamma} 
&=& \int_0^1 \int_{-r^n_{\gamma\delta}}^{+r^n_{\gamma\delta}} \partial_s (h^n_1\circ\psi^n_{\gamma\delta})(s,t)\,ds\,dt \\
&=& \int_{-r^n_{\gamma\delta}}^{+r^n_{\gamma\delta}} \int_0^1 \partial_t (h^n_2\circ\psi^n_{\gamma\delta})(s,t)\,dt\,ds \\
&=& \int_{-r^n_{\gamma\delta}}^{+r^n_{\gamma\delta}} 
((h^n_2\circ\psi^n_{\gamma\delta})(s,1)-(h^n_2\circ\psi^n_{\gamma\delta})(s,0))
\,ds \\
&=& 2\,\cdot\, m_{\gamma\delta} \,\cdot\, r^n_{\gamma\delta}. 
\end{eqnarray*}
Now let $\alpha=\gamma_0,\gamma_1,...,\gamma_k=\beta$ be the enumeration of vertices in $T^0$ on the 
unique path from $\alpha$ to $\beta$ and set $h^n_{i,j} = h^n_{\gamma\delta}$ for $\gamma=\gamma_i$, $\delta=\gamma_j$.
Then we have 
\begin{eqnarray*} 
 h^n_1(z_{\alpha}) - h^n_1(z_{\beta}) &=& h^n_1(z_{\alpha}) - h^n_{0,1} 
 \,+\, \sum_{i=0}^{k-1} \bigl(h^n_{i,i+1} - h^n_{i+1,i}\bigr) \\&&+\, \sum_{i=1}^{k-1} \bigl(h^n_{i,i-1} - h^n_{i,i+1}\bigr) 
 \,+\, h^n_{k,k-1} - h^n_1(z_{\beta}). 
\end{eqnarray*}
With $m_{i,j} = m_{\gamma\delta}$, $r^n_{i,j} = r^n_{\gamma\delta}$ for $\gamma=\gamma_i$, $\delta=\gamma_j$ we have 
\begin{equation*}
\sum_{i=0}^{k-1} \bigl(h^n_{i,i+1} - h^n_{i+1,i}\bigr) = \sum_{i=0}^{k-1} 2 m_{i,i+1} r^n_{i,i+1} = 2 s^n_{\alpha\beta},
\end{equation*} 
so that 
\begin{eqnarray*} 
 && (h^n_1(z_{\alpha}) - h^n_1(z_{\beta})) - 2 s^n_{\alpha\beta} \\
 && = \bigl(h^n_1(z_{\alpha}) - h^n_{0,1}\bigr) \,+\, \sum_{i=1}^{k-1} \bigl(h^n_{i,i-1} - h^n_{i,i+1}\bigr) \,+\, \bigl(h^n_{k,k-1} - h^n_1(z_{\beta})\bigr) \\
 && \stackrel{n\to\infty}{\longrightarrow}
 \bigl(h_{[\alpha],1}(z_{\alpha}) - h_{[\alpha],0,1}\bigr) \,+\, 
 \sum_{i=1}^{k-1} \bigl(h_{[\gamma_i],i,i-1} - h_{[\gamma_i],i,i+1}\bigr) \\
 && + \; \bigl(h_{[\beta],k,k-1} - h_{[\beta],1}(z_{\beta})\bigr).
\end{eqnarray*}
Note that the last expression depends only on the underlying nodal curve and is independent of the $\RS$-shifts. But this shows how 
the coordinates $s_{\alpha\beta}\in\IR$ describe the level structure and the relative $\IR$-shifts, in particular, 
two connected components belong to the same level precisely when $-\infty< s_{\alpha\beta} < +\infty$ for each 
$\alpha,\beta\in T^0$ representing the connected components in $T=T^0/\approx$. \\   

In order to fix the $S^1$-shifts and decorations in $\IZ_{|m_{\alpha\beta}|}$ at punctures between levels, observe that the 
coordinates $t_{\tau}\in S^1$ with $\tau\in \bar{E}^0_{\pm}$ determine decorations $t_{\tau}$ at the nodes $\tau\in\bar{E}_{\pm}$ 
corresponding to pairs of punctures connecting components on neighboring levels. Together with the 
$S^1$-coordinate $t_0$ they fix the $S^1$-shifts on each connected component of the level branched covering map as follows: \\

First for $\alpha\in T$ with $1\in\Lambda^+_{\alpha}$ we fix $h_{\alpha}$ 
by requiring that $h_{\alpha}$ maps the asymptotic marker at $z^+_1$ to $t_0\in S^1$. On the other hand, 
if $h_{\alpha}$ is fixed for some $\alpha\in T$, we can fix the $S^1$-shift for maps $h_{\beta}$ with $\alpha E\beta$ 
as follows: On the one hand, when $m_{\alpha\beta}=0$, i.e., when $\alpha$ and $\beta$ represent curves in the same level connected 
by a node $z_{\alpha\beta}\sim z_{\beta\alpha}$, the condition $h_{\alpha}(z_{\alpha\beta})=h_{\beta}(z_{\beta\alpha})$ 
immediately fixes the $S^1$-shift for $h_{\beta}$. Now consider the case when $m_{\alpha\beta}\neq 0$, i.e., $z_{\alpha\beta}$ 
and $z_{\beta\alpha}$ are positive or negative punctures. 
After choosing an asymptotic marker at $z_{\alpha\beta}$, which is mapped to $0\in S^1$ under $h_{\alpha}$, we 
can use the decoration $t_{[\alpha,\beta]}\in S^1$, $[\alpha,\beta]\in\bar{E}_{\pm}$ to get an asymptotic marker at $z_{\beta\alpha}$, and choose 
$h_{\beta}: (\Si_{\beta},j_{\beta})\to\RS$ so that it maps the asymptotic marker at $z_{\beta\alpha}$ to 
$0\in S^1$. Since $h_{\alpha}: (\Si_{\alpha},j_{\alpha})\to \RS\cong\IR\times\gamma$ is asymptotically cylindrical 
over the multiple $\gamma^{|m_{\alpha\beta}|}$, it follows that there are $|m_{\alpha\beta}|$ different possible 
choices for the asymptotic marker at $z_{\alpha\beta}$. Using the decoration $t_{\alpha\beta}$ this leads to 
$|m_{\beta\alpha}| = |m_{\alpha\beta}|$ different possible choices for the asymptotic marker at $z_{\beta\alpha}$, 
which however all lead to the same map $h_{\beta}: (\Si_{\beta},j_{\beta})\to\RS$.
Note that in this way we do not only get the holomorphic maps $h_{\alpha}:(\Si_{\alpha},j_{\alpha})\to \RS$ (up to the 
common $\IR$-shift in each level), but also the decorations $\theta_{\alpha\beta}\in \IZ_{|m_{\alpha\beta}|}$, i.e., we see that each element 
$(t_0,((r_{\tau},t_{\tau}),(s_{\alpha\beta})),\mu^{\pm}) \in S^1\times\TiM_{0,n}\times \IZ_{m^+}\times\IZ_{m^-}$ uniquely defines an 
element $(h,j,\mu,\theta)\in\CM$. \\

For the reverse, assume we are given an element $(h,j,\mu,\theta)\in\CM$, i.e., we are given maps $h_{\alpha}$ and $h_{\beta}$ for two components  
$\alpha,\beta$ connected by an edge in $(T,\LL)$, where we must only consider the case where $\alpha$ and $\beta$ live on different levels. Here we 
simultaneously have $|m_{\alpha\beta}|$ different possible choices for the asymptotic marker at $z_{\alpha\beta}$ and 
$|m_{\alpha\beta}|$ different possible choices for the asymptotic marker at $z_{\alpha\beta}$, which lead to 
$|m_{\alpha\beta}|$ different possible choices for the decoration $t_{[\alpha,\beta]}\in S^1$, which is then fixed using $\theta_{\alpha\beta}\in 
\IZ_{|m_{\alpha\beta}|}$. $\qed$ \\
   
\section{Obstruction bundle and Fredholm theory}

For determining the contribution of the moduli spaces of branched covers of trivial cylinders to the differential in rational 
symplectic field theory and contact homology, we show in section 3.1 that it suffices to study sections in a natural candidate for 
an obstruction bundle over the compactified moduli space of branched covers, the so-called cokernel bundle $\overline{\Coker}\CR_J$ 
of the Cauchy-Riemann operator $\CR_J$. Hence we follow the standard approach in Gromov-Witten theory of using obstruction bundles 
in order to deal with moduli spaces which are not regular in the sense that they are not transversally cut out by the Cauchy-Riemann operator. \\

\subsection{Cokernel bundle}

Denoting by $D_{h,j}: T_{h,j}\BB^{p,d}\to\EE^{p,d}_{h,j}$ the linearization of the nonlinear Cauchy-Riemann operator $\CR_J$ in the Banach space bundle $\EE^{p,d}\to\BB^{p,d}$ at $(h,j,\mu^{\pm})\in\IM\subset\BB$, which we discuss in detail 
in the upcoming subsection, the fibre at $(h,j,\mu^{\pm})$ of the bundle $\Coker\CR_J$ over $\IM$ as well as the bundle $\Coker_0\CR_J$ over $\IM^0$ 
is given by the cokernel of $D_{h,j}$
\begin{equation*}
  (\Coker\CR_J)_{(h,j,\mu^{\pm})} = (\Coker_0\CR_J)_{(h,j,\mu^{\pm})}= \coker D_{h,j}.
\end{equation*}   

For the extensions $\overline{\Coker}\CR_J$, $\overline{\Coker_0}\CR_J$ over the compactifications $\CM$ and $\overline{\IM^0}$ we make use of the fact that for every stratum given by the tree with level structure $(T,\LL)$ there also exists a natural Banach space bundle $\EE^{p,d}_{T,\LL}\to\BB^{p,d}_{T,\LL}$, so that we can require that the fibre over $(h,j,\mu,\theta)$ in the stratum $\IM_{T,\LL}$ or $\IM^0_{T,\LL,\ell_0}$ is given by the cokernel of the linearized operator $D_{h,j}: T_{h,j}\BB^{p,d}_{T,\LL}\to(\EE^{p,d}_{T,\LL})_{h,j}$. \\

Since the cokernel does not depend on the position of the asymptotic markers $\mu^{\pm}\in\IZ_{m^{\pm}}$, it 
follows that $\overline{\Coker}\CR_J$ ($\overline{\Coker_0}\CR_J$) naturally lives over the quotient 
$\CM/(\IZ_{m^+}\times\IZ_{m^-})$ ($\overline{\IM^0}/(\IZ_{m^+}\times\IZ_{m^-})$) rather than $\CM$ ($\overline{\IM^0}$) and we will view it this way. However it will later become important to consider it as a bundle over $\CM$ ($\overline{\IM^0}$) when we talk about orientations. \\

Denoting by $\Coker^{T,\LL}\CR_J$, $\Coker^{T,\LL}_0\CR_J$ the restrictions of $\overline{\Coker}\CR_J$, $\overline{\Coker_0}\CR_J$ to $\IM_{T,\LL}$, 
$\IM^0_{T,\LL}$, it follows from  
\begin{eqnarray*}
\BB^{p,d}_{T,\LL}/(\IZ_{m^+}\times\IZ_{m^-})&=& \frac{\BB^{p,d}_{T_1}\times ... \times\BB^{p,d}_{T_L}}{\Delta \times \IZ_{m^+}\times\IZ_{m^-}},\\
\EE^{p,d}_{T,\LL}&=& \pi_{1,1}^* \EE^{p,d}_{T_{1,1}}\oplus ...\oplus\pi_{L,N_L}^*\EE^{p,d}_{T_{L,N_L}}
\end{eqnarray*}
with the projections 
\begin{equation*} 
   \pi_{\ell,k}: \BB^{p,d}_{T,\LL}/(\IZ_{m^+}\times\IZ_{m^-})= \frac{\BB^{p,d}_{T_1}\times ...
   \times\BB^{p,d}_{T_L}}{\Delta \times \IZ_{m^+}\times\IZ_{m^-}} 
   \to \BB^{p,d}_{T_{\ell,k}}/(\IZ_{m^+_{\ell,k}}\times\IZ_{m^-_{\ell,k}}), 
\end{equation*}
that they are given as direct sums
\begin{eqnarray*}
 \Coker^{T,\LL}\CR_J &=& \pi_{1,1}^* \Coker^{T_{1,1}}\CR_J\oplus ...\oplus\pi_{L,N_L}^*\Coker^{T_{L,N_L}}\CR_J, \\
 \Coker^{T,\LL}_0\CR_J &=& \pi_{1,1}^* \Coker^{T_{1,1}}\CR_J\oplus ...\oplus\pi_{\ell_0,N_{\ell_0}}^*\Coker^{T_{\ell_0,N_{\ell_0}}}_0\CR_J
      \oplus \\ && ...\oplus\pi_{L,N_L}^*\Coker^{T_{L,N_L}}\CR_J,
\end{eqnarray*}
where $m^{\pm}_{\ell,k}=\prod_{\alpha\in T_{\ell,k}} m^{\pm}_{\alpha}$, 
$m^{\pm}_{\alpha}=m^{\pm}_{\alpha,1}\cdot ...\cdot 
m^{\pm}_{\alpha,n^{\pm}_{\alpha}}$ and $\Coker^{T_{\ell,k}}$ ($\Coker^{T_{\ell,k}}_0$) denotes the cokernel bundle over 
$\IM_{T_{\ell,k}}/(\IZ_{m^+_{\ell,k}}\times\IZ_{m^-_{\ell,k}})$ ($\IM^0_{T_{\ell,k}}/(\IZ_{m^+_{\ell,k}}\times\IZ_{m^-_{\ell,k}})$) 
for $\ell=1,...,L$, $k=1,...,N_{\ell}$. Note that there exists no natural map from $\IM_{T,\LL}$ ($\IM^0_{T,\LL}$) to 
$\IM_{T_1},...,\IM_{T_L}$ and hence to $\IM_{T_{1,1}},...,\IM_{T_{L,N_L}}$, since we quotient out the diagonal $\Delta$, i.e., identify pairs of asymptotic markers if they represent the same decoration. \\

Recall from subsection 1.2 that when $\IM_{T,\LL}$ belongs to the codimension one boundary of $\IM$ it is of the form 
$\IM_{T,\LL} = \IM_1\times_{\IZ_{m_{1,2}}}\IM_2$, where $\IM_1$ and $\IM_2$ are moduli spaces of possibly disconnected orbit curves  
without nodes. Note that the compactification of the fibre product $\overline{\IM_1\times_{\IZ_{m_{1,2}}}\IM_2}\subset\del\CM$ can directly be identified with the 
fibre product of the compactifications,
\begin{equation*} \overline{\IM_1\times_{\IZ_{m_{1,2}}}\IM_2} = \CM_1\times_{\IZ_{m_{1,2}}}\CM_2. \end{equation*}
For this observe that the partitioning of the levels of a limiting curve in $\overline{\IM_1\times_{\IZ_{m_{1,2}}}\IM_2}$ into levels belonging 
to the compactification $\CM_1$ or $\CM_2$, respectively, follows from the conservation of the total Euler characteristic under degeneration of 
punctured Riemann surfaces. Denoting by $\overline{\Coker}^1\CR_J$ and $\overline{\Coker}^2\CR_J$ the extensions of the cokernel bundles 
over $\IM_1$, $\IM_2$ to the corresponding compactified moduli spaces, it directly follows from the form of 
$\overline{\Coker}\CR_J$ over the strata $\IM_{T,\LL}$ that  
\begin{equation*} 
 \overline{\Coker}\CR_J|_{\CM_1\times_{\IZ_{m_{1,2}}}\CM_2} \;=\; 
 \pi_1^*\overline{\Coker}^1\CR_J \oplus \pi_2^*\overline{\Coker}^2\CR_J,
\end{equation*}
with the projections $\pi_{1,2}: \CM/\IZ_{m^{\pm}}\to\CM_{1,2}/\IZ_{m^{\pm}_{1,2}}$. 
For the cokernel bundle $\overline{\Coker_0}\CR_J$ it follows in the same way that 
\begin{eqnarray*} 
 \overline{\Coker_0}\CR_J|_{\overline{\IM^0_1}\times_{\IZ_{m_{1,2}}}\CM_2} &=& \pi_1^*\overline{\Coker_0}^1\CR_J \oplus \pi_2^*\overline{\Coker}^2\CR_J, \\
 \overline{\Coker_0}\CR_J|_{\CM_1\times_{\IZ_{m_{1,2}}}\overline{\IM^0_2}} &=& \pi_1^*\overline{\Coker}^1\CR_J \oplus \pi_2^*\overline{\Coker_0}^2\CR_J
\end{eqnarray*}
and 
\begin{eqnarray*}
 \overline{\Coker_0}\CR_J|_{\{\point\}\times \CM} \;=\; \overline{\Coker_0}\CR_J|_{\CM\times\{\point\}} \;=\;\overline{\Coker}\CR_J.
\end{eqnarray*}
$ $\\

In order to show that $\overline{\Coker}\CR_J$ indeed serves as an obstruction bundle, we show in the upcoming subsection 2.2
that on every stratum $\IM_{T,\LL}\subset\CM$ we have  
\begin{equation*}
        \ker D_{h,j} = T_{h,j}\IM_{T,\LL}
\end{equation*} 
at every $(h,j,\mu^{\pm},\theta)\in\IM_{T,\LL}$, see subsection 3.1 below, which then automatically implies that $\Coker^{T,\LL}$ is indeed a smooth vector bundle over 
$\IM_{T,\LL}$. In order to show that these vector bundle naturally fit together to a smooth vector bundle $\overline{\Coker}\CR_J$ over the manifold with corners 
$\CM$, we prove in subsection 2.3 a linear gluing theorem relating the cokernel bundles over different strata of $\CM$. \\

\subsection{Linearized operator}

For all this we first need to understand the linearization $D_{h,j}$ of the Cauchy-Riemann operator $\CR_J$ 
at $(h,j) \in \IM/(\IZ_{m^+}\times\IZ_{m^-})$. In what follows we formulate our statements only for 
the cokernel bundle $\overline{\Coker}\CR_J$, since the statements for $\overline{\Coker_0}\CR_J$ then follow immediately. 
For the Banach manifold setup we follow [BM] and the author's expositions in [F].\\

Recall from the definition of the moduli spaces that we fixed $n^+$ positive and $n^-$ negative punctures $z^{\pm}_1,...,z^{\pm}_{n^{\pm}}\in S^2$ and 
fixed cylindrical coordinates 
\begin{equation*} \psi^{\pm}_k: \IR^{\pm}_0\times S^1\hookrightarrow\Si \end{equation*} 
around each puncture $z^{\pm}_k$, $k\in\{1,...,n^{\pm}\}$ on the punctured sphere $\Si=S^2-\{z^{\pm}_1,...,z^{\pm}_{n^{\pm}}\}$.
Let the space $H^{1,p,d}_{\cst}(\Si,\IC)$ consist of all 
maps from $\Si$ to $\IC$ differing asymptotically from a 
constant one by a function, which is still in $H^{1,p}$ after multiplication with an asymptotic weight. To be precise,  
any $v \in H^{1,p,d}_{\cst}(\Si,\IC)$ is in $H^{1,p}_{\loc}$ and for any puncture $z^{\pm}_k$ 
there exist $(s_0^{\pm,k},t_0^{\pm,k})\in \IR^2\cong\IC$, so that the function 
\begin{equation*}
 \IR^{\pm}\times S^1 \to \IC, \,\, 
 (s,t)\mapsto [(v\circ\psi^{\pm}_k)(s,t)-(s_0^{\pm,k},t_0^{\pm,k})] \cdot e^{\pm d\cdot s} 
\end{equation*}
is in $H^{1,p}$. Let further $L^{p,d}(T^*\Si\otimes_{j,i}\IC)$ denote the space of 
$(j,i)$-antiholomorphic one-forms on $\Si$ with values in $\IC$, which are still in $L^p$ after multiplication 
with the asymptotic weight $e^{\pm d\cdot s}$.  \\

With $h^*\xi$ denoting the pullback of the subbundle $\xi \subset TV$ under the branched covering map $h: (\Si,j) \to (\RS,i) \cong (\IR\times\gamma,J)$, 
we introduce the spaces $H^{1,p}(h^*\xi)$ of sections and $L^p(T^*\Si\otimes_{j,J_{\xi}}h^*\xi)$ of $(j,J_{\xi})$-antiholomorphic one-forms 
on $\Si$ with values in $h^*\xi$, where the $H^{1,p}$- and $L^p$-topologies are defined with respect to any trivialization of $\xi$ along the fixed
Reeb orbit $\gamma$. \\

Following [Sch] and [BM], [F] there exists a Banach space bundle $\EE^{p,d}$ over a Banach manifold of maps $\BB^{p,d}$ in which the Cauchy-Riemann operator 
$\CR_J$ extends to a smooth section. In our special case it follows that the fibre is given by  
\begin{equation*} 
 \EE^{p,d}_{h,j} = L^{p,d}(T^*\Si\otimes_{j,i}\IC) \oplus L^p(T^*\Si\otimes_{j,J_{\xi}}h^*\xi),
\end{equation*}
while the tangent space to the Banach manifold of maps $\BB^{p,d}= \BB^{p,d}(\gamma^{m^+_1},...,\gamma^{m^+_{n^+}};\gamma^{m^-_1},...,\gamma^{m^-_{n^-}})$ 
at $(h,j) \in \IM/(\IZ_{m^+}\times\IZ_{m^-})$ is given by   
\begin{equation*}
  T_{h,j}\BB^{p,d}(V;(\gamma^{m^{\pm}_i})) = H^{1,p,d}_{\cst}(\Si,\IC)\oplus H^{1,p}(h^*\xi)\oplus T_j\IM_{0,n}.
\end{equation*} 
Note that we use the complex splitting of the tangent bundle $T(\IR\times V) = \IC \oplus \xi$ in order to write tangent 
spaces and fibres as direct sums.\\

In order to give an explicit formula for the linearization $D_{h,j}$ of $\CR_J$ we choose a complex connection on 
$(\xi,J_{\xi})$ which we extend to a connection $\nabla$ on $T(\IR\times V) =\IC\oplus\xi$, 
$\IC=\IR\cdot\partial_s \oplus \IR\cdot R$ by requiring $\IR$-invariance and $\nabla\partial_s =\nabla R =0$, where $\del_s$ is the 
$\IR$-direction and $R$ the Reeb vector field of the stable Hamiltonian structure. 
For this connection it follows that the linearization $D_{h,j}$ of $\CR_J$ at branched covers of orbit cylinders $(h,j)$ is of a special form. \\
\\
{\bf Proposition 2.1:} {\it With respect to the complex connection $\nabla$ on $T(\IR\times V)$ from above, the linearization $D_{h,j}$ 
of $\CR_J$ at $(h,j)\in\IM/(\IZ_{m^+}\times\IZ_{m^-})$ is given by} 
\begin{eqnarray*}
 &D_{h,j}:& H^{1,p,d}_{\cst}(\Si,\IC)\oplus H^{1,p}(h^*\xi)\oplus T_j\IM_{0,n} \\
 &&\to L^{p,d}(T^*\Si\otimes_{j,i}\IC) \oplus L^p(T^*\Si\otimes_{j,J_{\xi}}h^*\xi),\\
 && D_{h,j} \cdot (v_1,v_2,y) = (\CR v_1+ D_j y, D_h^{\xi} v_2),
\end{eqnarray*}
{\it where $\CR: H^{1,p,d}_{\cst}(\Si,\IC) \to L^{p,d}(T^*\Si\otimes_{j,i}\IC)$ is the standard 
Cauchy-Riemann operator,} 
\begin{eqnarray*}
 &&D_h^{\xi}: H^{1,p}(h^*\xi) \to L^p(T^*\Si\otimes_{j,J_{\xi}}h^*\xi),\,\,\\ 
 &&D_h^{\xi} v_2 = \nabla v_2 + J_{\xi} \cdot \nabla v_2 \cdot j + \nabla_{dh}v_2 + J_{\xi}\nabla_{i\,dh}v_2 
\end{eqnarray*}
{\it describes the linearization of $\CR_J$ in the direction of $\xi \subset TV$ and} 
\begin{equation*} 
D_j: T_j\IM_{0,n} \to L^{p,d}(T^*\Si\otimes_{j,i}\IC),\,\,D_j y = i\cdot dh\cdot y.
\end{equation*}
{\it describes the variation of $\CR_J$ with $j\in\IM_{0,n}$.} \\
\\
{\it Proof:} Since $\nabla$ is a complex connection, it is well-known, see e.g. [Sch], that the linearization $D_h: H^{1,p,d}_{\cst}(\Si,\IC)\oplus H^{1,p}(h^*\xi)\to 
L^{p,d}(T^*\Si\otimes_{j,i}\IC) \oplus L^p(T^*\Si\otimes_{j,J_{\xi}}h^*\xi)$ of $\CR_J$ fixing the complex structure $j\in\IM_{0,n}$ is given by
\begin{equation*} 
 D_h\cdot v \,=\, \nabla v + J\cdot\nabla v\cdot j + \Tor(dh,v) + J\Tor(J\cdot dh,v),
\end{equation*}
where $\Tor(X,Y) = \nabla_X Y - \nabla_Y X - [X,Y]$. First it follows from the special form of $\nabla$ that 
\begin{equation*}
 \nabla v + J\cdot\nabla v\cdot j = (\CR v_1, \nabla v_2 + J_{\xi}\cdot\nabla v_2\cdot j).
\end{equation*}
for $v=(v_1,v_2) \in H^{1,p,d}_{\cst}(\Si,\IC)\oplus H^{1,p}(h^*\xi)$. On the other hand,
\begin{eqnarray*}
 && \Tor(dh,v) + J\Tor(J\cdot dh,v) = \\
 && \nabla_{dh} v + J\cdot\nabla_{J\,dh} v - \nabla_v dh - J\cdot\nabla_v(J\cdot dh) - ([dh,v] + J[J\,dh,v]) = \\
 && \nabla_{dh} v + J\cdot\nabla_{J\,dh} v + J\cdot (L_v(J\,dh) - J\cdot L_v dh) = \\
 && \nabla_{dh} v + J\cdot\nabla_{J\,dh} v + J\cdot L_v J\cdot dh = \nabla_{dh} v + J\cdot\nabla_{J\,dh} v.
\end{eqnarray*}
>From $\nabla\partial_s =\nabla R =0$ it follows that $\Tor(dh,v_1) + J\Tor(J\cdot dh,v_1) = 0$, while for $v_2\in\xi$ 
we have $\nabla_{dh} v_2 + J\cdot\nabla_{J\,dh} v_2 \in \xi$, so that $D_{h} \cdot (v_1,v_2) = 
(\CR v_1, D_h^{\xi} v_2)$ with $D_h^{\xi}$ as in the lemma. Finally note that for the linearization 
of $\CR_J$ in the direction of $\IM_{0,n}$ there is obviously no variation in the $\xi$-direction,
\begin{equation*} 
D_j y= (i\cdot dh\cdot y, 0).\,\,\,\qed
\end{equation*}

Based on this result, the following lemmata describe kernel and cokernel of $D_{h,j}$. \\
\\
{\bf Proposition 2.2:} {\it The standard Cauchy-Riemann operator $\CR: H^{1,p,d}_{\cst}(\Si,\IC) \to L^{p,d}(T^*\Si\otimes_{j,i}\IC)$ is onto, 
so that} 
\begin{equation*} \coker D_{h,j} = \coker D_h^{\xi}. \end{equation*}
\\
{\it Proof:} The second part of the statement follows from the upper-triangle-form of $D_{h,j}$. 
For the first statement about $\CR$ fix a splitting
\begin{equation*}
 H^{1,p,d}_{\cst}(\Si,\IC) = H^{1,p,d}(\Si,\IC) \oplus \Gamma^n
\end{equation*}
where $\Gamma^n\subset C^{\infty}(\Si,\IC)$ is a $2n$-dimensional space of functions storing the constant shifts (see [BM]). 
Given a function $\varphi_d: \Si \to \IR$ 
with $(\varphi_d \circ \psi^{\pm}_k)(s,t) = e^{\pm d\cdot s}$, multiplication with $\varphi_d$ defines 
isomorphisms
\begin{eqnarray*}
H^{1,p,d}(\Si,\IC) &\stackrel{\cong}{\longrightarrow}& H^{1,p}(\Si,\IC), \\
L^{p,d}(T^*\Si\otimes_{j,i}\IC) &\stackrel{\cong}{\longrightarrow}& L^p(T^*\Si\otimes_{j,i}\IC),
\end{eqnarray*}
under which $\CR$ corresponds to a perturbed Cauchy-Riemann operator 
\begin{equation*}
 \CR_d = \CR + S_d: H^{1,p}(\Si,\IC) \to L^p(T^*\Si\otimes_{j,i}\IC).
\end{equation*}
With the asymptotic behaviour of $\varphi_d$ one computes
\begin{equation*}
 S^{\pm,k}_d(t) = (S_d \circ \psi^{\pm}_k)(\pm\infty,t) = \diag(\mp d,\mp d)
\end{equation*}
so that the Conley-Zehnder indices for the corresponding paths $\Psi^{\pm,k}: \IR \to \Sp(2)$ of 
symplectic matrices is $\mp 1$ for $d>0$ sufficiently small. Hence the index of $\CR: 
H^{1,p,d}_{\cst}(\Si,\IC) \to L^{p,d}(T^*\Si\otimes_{j,i}\IC)$ is given by
\begin{equation*}
\ind \CR = \dim\Gamma^n + \ind \CR_d = 2n - n + 1\cdot(2-n) = 2,
\end{equation*}
where the first summand is the dimension of $\Gamma^n$ and the second is the sum of the Conley-Zehnder 
indices. On the other hand it follows from Liouville's theorem that the kernel of $\CR$ consists 
precisely of the constant functions on $\Si$, so that $\dim \coker \CR = 0$. $\qed$ \\
\\
{\bf Proposition 2.3:} {\it The operator $D_h^{\xi}$ has a trivial kernel, so that} 
\begin{equation*}
 \ker D_{h,j} = T_{h,j} (\IR\times\IM).
\end{equation*}
\\
{\it Proof:} Note that here the second part of the statement follows from the first one using proposition 2.2 as follows: First it is clear that 
we have the inclusion $T_{h,j} (\IR\times\IM) \subset \ker D_{h,j}$, since $\IR\times\IM=\CR_J^{-1}(0)$. On the other hand, using the first part of the 
statement we know that the kernel of $D_{h,j}$ consists of all pairs $(\hb,y)\in H^{1,p,d}_{\cst}(\Si,\IC)\oplus T_j\IM_{0,n}$ satisfying 
$\CR\hb+D_j y = 0$. Since $\CR: H^{1,p,d}_{\cst}(\Si,\IC)\to L^{p,d}(T^*\Si\otimes_{j,i}\IC)$ is surjective, it follows that $\ker D_{h,j}$ projects 
surjectively onto $T_j\IM_{0,n}$, where the fibre can be identified with $\ker\CR=\IC$. In particular, we have that the dimension of $\ker D_{h,j}$ agrees 
with the dimension of $T_{h,j}(\RS\times\IM_{0,n})=T_{h,j}(\IR\times\IM)$, so that the inclusion must indeed be an equality. \\

The statement about the kernel of $D_h^{\xi}$ is the linearized version of lemma 5.4 in [BEHWZ]. For chosen $h=(h_1,h_2): 
(\Si,j)\to\RS\cong (\IR\times\gamma,J)$ and $v_2\in \ker D_h^{\xi} \subset H^{1,p}(h^*\xi)$ we can use the exponential map 
of some Riemannian metric on $\IR\times V$ to get for $r>0$ sufficiently small a family of curves 
\begin{equation*} \exp_h rv_2=(h_1,\exp_{h_2}rv_2): (\Si,j)\to\IR\times V. \end{equation*} 
Note that their $\omega$-energies $E_{\omega}(\exp_h rv_2) = \int_{\Si} (\exp_{h_2} rv_2)^*\omega$ by homological reasons agree with the $\omega$-energy 
of $h$ and hence vanish, 
\begin{equation*}
 E_{\omega}(\exp_h rv_2) = E_{\omega}(h) = 0,
\end{equation*}
since all curves in the family are asymptotically cylindrical over the same closed Reeb orbits near the punctures. 
Choosing an atlas $(U_{\alpha},\varphi_{\alpha})_{\alpha\in A}$ for the complex manifold $\Si$ with 
local holomorphic coordinates $(s_{\alpha},t_{\alpha})$ on $U_{\alpha}\subset \Si$, together with 
a subordinate partition of unity $(\psi_{\alpha})_{\alpha\in A}$, observe that the above integral can be rewritten as 
\begin{eqnarray*} 
 && \int_{\Si} (\exp_{h_2} rv_2)^*\omega \\
 && = \sum_{\alpha} \int_{U_{\alpha}} \psi_{\alpha} \cdot 
 \omega( \partial_{s_{\alpha}} \exp_{h_2}rv_2, \partial_{t_{\alpha}} \exp_{h_2}rv_2) ds_{\alpha}\wedge dt_{\alpha}\\
 && = \sum_{\alpha} \int_{U_{\alpha}} \psi_{\alpha} \cdot 
 \omega_{\xi}( \pi_{\xi}\partial_{s_{\alpha}} \exp_{h_2}rv_2, \pi_{\xi}\partial_{t_{\alpha}} \exp_{h_2}rv_2) 
 ds_{\alpha}\wedge dt_{\alpha},
\end{eqnarray*}
where $\pi_{\xi}$ denotes the projection $TV=\IC\oplus\xi\to\xi$ and the second equality follows from 
$R\in\ker\omega$. With the metric $\langle\cdot,\cdot\rangle_{\xi} = \omega_{\xi}(\cdot,J_{\xi}\cdot)$ on $\xi$ we get that 
the latter is equal to
\begin{eqnarray*} 
 && \sum_{\alpha} \int_{U_{\alpha}} \psi_{\alpha} \cdot 
 \langle \pi_{\xi}\partial_{s_{\alpha}} \exp_{h_2}rv_2, - \pi_{\xi}\, J\, \partial_{t_{\alpha}} \exp_{h_2}rv_2\rangle_{\xi} 
 \,\,ds_{\alpha}\wedge dt_{\alpha} = \\
 && \sum_{\alpha} \int_{U_{\alpha}} \psi_{\alpha} \cdot 
 \langle \pi_{\xi}\partial_{s_{\alpha}} \exp_{h_2}rv_2, \\
 && \pi_{\xi}\partial_{s_{\alpha}} \exp_{h_2}rv_2 -  \pi_{\xi}\CR_J \exp_{h_2}rv_2 \cdot\del_{s_{\alpha}}\rangle_{\xi} \,\,ds_{\alpha}\wedge dt_{\alpha}.  
\end{eqnarray*}
For $r=0$ observe that we have 
$\pi_{\xi}\CR_J \exp_{h_2}rv_2 = \pi_{\xi}\CR_J h_2 =0$ and $\pi_{\xi}\partial_{s_{\alpha}} \exp_{h_2}rv_2 = \pi_{\xi} \partial_{s_{\alpha}}h_2 =0$, 
where the latter uses that $h=(h_1,h_2)$ is a branched cover of a trivial cylinder, i.e., $h_2$ is contained in a trajectory of the 
Reeb vector field. Letting $\Phi_{h_2}rv_2$ denote parallel transport on $\xi$ starting from $h_2$ in the direction $v_2$ with respect to 
the complex connection $\nabla$ from before, where we additionally assume that it preserves $\omega_{\xi}$, i.e., the metric $\<\cdot,\cdot\>_{\xi}$, 
the Leibniz rule implies 
\begin{eqnarray*} 
 && \frac{d^2}{dr^2}|_{r=0} \bigl\langle \pi_{\xi}\partial_{s_{\alpha}} \exp_{h_2}rv_2, 
  \pi_{\xi}\partial_{s_{\alpha}} \exp_{h_2}rv_2 -  \pi_{\xi}\CR_J \exp_{h_2}rv_2\cdot\del_{s_{\alpha}}\bigr\rangle_{\xi} = \\
 && \frac{d^2}{dr^2}|_{r=0} \bigl\langle (\Phi_{h_2}rv_2)^{-1} \pi_{\xi}\partial_{s_{\alpha}} \exp_{h_2}rv_2, \\ 
 && (\Phi_{h_2}rv_2)^{-1}  \pi_{\xi}\partial_{s_{\alpha}} \exp_{h_2}rv_2 - (\Phi_{h_2}rv_2)^{-1} \pi_{\xi} \CR_J 
 \exp_{h_2}rv_2\cdot\del_{s_{\alpha}}\bigr\rangle_{\xi} = \\
 && \Bigl\langle\frac{d}{dr}|_{r=0} (\Phi_{h_2}rv_2)^{-1} \pi_{\xi}\partial_{s_{\alpha}} \exp_{h_2}rv_2, \\
 && \frac{d}{dr}|_{r=0}(\Phi_{h_2}rv_2)^{-1}  \pi_{\xi}\partial_{s_{\alpha}} \exp_{h_2}rv_2 \\
 && - \frac{d}{dr}|_{r=0}(\Phi_{h_2}rv_2)^{-1}  \pi_{\xi}\CR_J\exp_{h_2}rv_2\cdot\del_{s_{\alpha}}\Bigr\rangle_{\xi} = \\ 
 && \langle\nabla_{s_{\alpha}} v_2, \nabla_{s_{\alpha}} v_2 - D_h^{\xi} v_2 \cdot\del_{s_{\alpha}}\rangle_{\xi}.
 = |\nabla_{s_{\alpha}} v_2|_{\xi}^2  
\end{eqnarray*}
Hence, 
\begin{eqnarray*}
 0 &=& \frac{d^2}{dr^2}E_{\omega}(\exp_h rv_2) \\
   &=& \sum_{\alpha} \int_{U_{\alpha}} \psi_{\alpha} \cdot |\nabla_{s_{\alpha}} v_2|_{\xi}^2 \,\,ds_{\alpha}\wedge dt_{\alpha}, 
\end{eqnarray*}
so that $\nabla_{s_{\alpha}} v_2 = 0$. Since by the same arguments $\nabla_{t_{\alpha}} v_2 = 0$ we indeed have 
$\nabla v_2=0$ on $\Si$, which by $v_2 \in H^{1,p}(h^*\xi)$ implies $v_2=0$ as desired. $\qed$ \\    

Since the kernel of the linearized operator agrees with the tangent space to the moduli space of orbit curves, the dimension 
of the kernel of the linearization of $\CR_J$ is constant on the moduli space. Together with the constancy of the Fredholm index it
proves that the cokernel bundle is of constant rank over $\IM/(\IZ_{m^+}\times\IZ_{m^-}) = S^1\times\IM_{0,n}$. By the same 
arguments it follows that the cokernel bundle over the moduli space $\IM_{T,\LL}$ is of constant rank for any tree with level structure 
$(T,\LL)$. As in [MDSa] this rank constancy proves that $\Coker^{T,\LL}\CR_J$ is indeed a smooth vector bundle over the smooth manifold $\IM_{T,\LL}$:\\
\\
{\bf Corollary 2.4:} {\it $\Coker^{T,\LL}\CR_J$ naturally carries the structure of a smooth vector bundle over $\IM_{T,\LL}$.} \\

\subsection{Linear gluing}

This subsection is concerned with the following extension of the above result: \\
\\
{\bf Proposition 2.5:} {\it Using a linear gluing construction (relating the cokernel bundle over the moduli space with the cokernel bundles over the boundary strata) 
we can equip the cokernel bundle $\overline{\Coker}\CR_J$ over the compactified moduli space $\CM$ with the structure of a smooth vector bundle over a 
smooth manifold with corners.} \\

Recall that we have shown in proposition 1.4 that the compactified moduli space $\CM$ carries the structure of a smooth manifold with corners.
For the proof it suffices to establish linear gluing theorems for the cokernel bundle under gluing of 
the underlying moduli spaces of branched covers. For the gluing theorems we must distinguish the case of gluing of curves on different levels, i.e., 
gluing at punctures, and gluing of curves in the same level, which corresponds to gluing at a node. \\
$ $\\
{\it Gluing of moduli spaces:} \\
 
In order to describe gluing of the cokernel bundles, we must start with gluing of the underlying moduli spaces of branched covers. Although 
these moduli spaces are nonregular and we hence cannot apply the usual gluing theorems, the gluing can explicitly described as follows: \\

We will restrict our attention only to the case of gluing at a puncture and claim that the gluing at nodes follows along the same lines. Using the notation introduced in 1.2, let $(T,\LL)=(T,E,\Lambda^{\pm},\LL)$ denote the 
tree with level structure given by 
$T=\{1,2\}$, $1E2$ and $\LL(1)=1$, $\LL(2)=2$. Note that the moduli space $\IM_{T,\LL}$ is given by the fibre product 
$\IM_1\times_{\IZ_{m_{12}}}\IM_2$ where $\IM_1$, $\IM_2$ denote moduli spaces of connected branched covers without nodes. 
Let $(h,j,\theta,\mu^{\pm})\in\IM_{T,\LL}$ with $h=(h_1,h_2)$, $j=(j_1,j_2)$, 
$\theta=\theta_{12}\in\IZ_{m_{1,2}}$ and $\mu^{\pm}=(\mu^{\pm}_1,\mu^{\pm}_2)$ with $\mu^{\pm}_1=(\mu^{\pm}_k)_{k\in\Lambda^{\pm}_1}$, 
$\mu^{\pm}_2=(\mu^{\pm}_k)_{k\in\Lambda^{\pm}_2}$. Then the underlying punctured spheres are
$\Si_1=S^2-(Z_1^+\cup Z_1^-)$, $\Si_2=S^2-(Z_2^+\cup Z_2^-)$, where the connecting pair of punctures is $(z_{12},z_{21})$ with 
$z_{12}\in Z^-_1$ and $z_{21}\in Z^+_2$. We define the family of glued curves 
\begin{equation*} 
        (h^r,j^r,\mu^{\pm}) \;=\; \sharp_r (h,j,\theta,\mu^{\pm}) \;=\; (h_1,j_1,\mu_1) \sharp_{r,\theta} (h_2,j_2,\mu_2) 
\end{equation*}    
as follows, where $r=r_{12}\in\IR^+$ denotes the gluing parameter: \\

When $\psi_{12}:\IR^-\times S^1\to\Si_1$, $\psi_{21}:\IR^+\times S^1\to\Si_2$ denote the fixed cylindrical coordinates around $z_{12}\in\Si_1$, 
$z_{21}\in\Si_2$, let $\Si_1^r$, $\Si_2^r$ denote the punctured surfaces with boundary given by cutting out the half-cylinders $(-\infty,-r)\times 
S^1$, $(+r,+\infty)\times S^1$, respectively,
\begin{eqnarray*}
 \Si_1^r = \Si_1 - \psi_{12}((-\infty,-r)\times S^1),\;\;\Si_2^r = \Si_2-\psi_{21}((+r,+\infty)\times S^1).
\end{eqnarray*}
We introduce the punctured surface $\Si^r$ underlying $(h^r,j^r,\mu^{\pm})$ by gluing $\Si_1^r$ and $\Si_2^r$ along the boundary with the twist 
given by the maps $h_1$ and $h_2$ and the decoration $\theta_{12}$, 
\begin{equation*} 
  \Si^r \;=\; \Si_1^r \sharp_{\theta_{12}} \Si_2^r \;=\; \Si_1^r\coprod\Si_2^r/\{\psi_{12}(-r,t)\sim\psi_{21}(+r,t+\theta_{12})\}.
\end{equation*}
Note that here the decoration $\theta_{12}$ is viewed as an element in $S^1$ rather than in $\IZ_{m_{12}}$.
For this recall that the maps $h_1$, $h_2$ determine $m_{1,2}$ different asymptotic markers at $z_{12}\in\Si_1$ and $z_{21}\in\Si_2$, which determine 
$S^1$-coordinates in the cylindrical coordinates $\psi_{12}$ and $\psi_{21}$. Hence there are $m_{12}$ possible ways to glue 
$\Si_1^r$ and $\Si_2^r$ so that these $S^1$-coordinates match, and the element in $\IZ_{m_{12}}$ singles out the unique gluing twist. 
Note that $\Si^r$ is again diffeomorphic to a punctured sphere and the complex structures $j_1$ on $\Si_1$ and $j_2$ on $\Si_2$ determine 
a complex structure $j^r$ on $\Si^r$ since both agree with the standard complex structure on the embedded half-cylinders determined 
by $\psi_{12}$ and $\psi_{21}$. On the other hand, the branched covering map $h^r: (\Si^r,j^r)\to\RS$ is unique up to $\IR$-shift by the requirement that 
the asymptotic markers of $h^r$ match with those of the maps $h_1$ on $\Si_1^r$ and $h_2$ on $\Si_2^r$ and exists by the choice of the gluing twist 
$\theta_{12}\in S^1$, since it is chosen so that the $S^1$-shifts for $h_1$ and $h_2$ agree. Hence we found a natural gluing map for gluing 
at punctures
\begin{equation*} \sharp: (\IM_1\times_{\IZ_{m_{12}}}\IM_2) \times(0,+\infty)\hookrightarrow \IM,\;\; ((h,j,\theta,\mu^{\pm}),r)\mapsto (h^r,j^r,\mu^{\pm}). 
\end{equation*}

{\it Linear gluing of the cokernel bundle:} \\

We now start with the gluing of the cokernel bundles. It follows from proposition 2.2 in the last subsection that 
the fibres of the cokernel bundle over $(h,j)\in \IM$ are given by
\begin{equation*} 
 (\Coker\CR_J)_{(h,j)} \;=\; \coker D_{h,j} \;=\; \coker D^{\xi}_h \;=\; \ker (D^{\xi}_h)^*, 
\end{equation*}
where 
\begin{equation*} (D^{\xi}_h)^*: H^{1,q}(T^*\Si \otimes_{j,J_{\xi}} h^*\xi) \to L^q(h^*\xi),\;\; 1/p + 1/q = 1 \end{equation*}
denotes the formal adjoint of the linearization $D^{\xi}_h: H^{1,p}(h^*\xi)\to L^p(T^*\Si\otimes_{j,J_{\xi}} h^*\xi)$ of $\CR_J$ in the direction of the 
hyperplane distribution $\xi\subset TV$. Since by elliptic regularity all occuring kernels and hence cokernels are independent of the choice of 
$p\geq 2$, see [Sch], we set in the following $p=q=2$. Note that since $\ker D^{\xi}_h=\{0\}$ by proposition 2.3, the operators 
$(D^{\xi}_h)^*$ are surjective. \\

Recall that the cokernel bundle over 
$\IM_{T,\LL} = \IM_1\times_{\IZ_{m_{12}}}\IM_2$ is given as direct sum, 
\begin{equation*} 
 \Coker^{T,\LL}\CR_J \;=\; \pi_1^*\Coker^1\CR_J \oplus \pi_2^*\Coker^2\CR_J,
\end{equation*}
where $\Coker^1,\Coker^2\CR_J$ denote the cokernel bundles over $\IM_1,\IM_2$ and $\pi_1,\pi_2$ the projections 
from $\IM^{T,\LL}/(\IZ_{m^+}\times\IZ_{m^-})$ to $\IM_1/(\IZ_{m_1^+}\times\IZ_{m_1^-})$, $\IM_2/(\IZ_{m_2^+}\times\IZ_{m_2^-})$, 
respectively. Let $(h,j,\theta,\mu^{\pm})\in\IM_{T,\LL}=\IM_1\times_{\IZ_{m_{12}}}\IM_2$ with $h=(h_1,h_2)$, $j=(j_1,j_2)$. 
For 
\begin{equation*} 
  \eta=(\eta_1,\eta_2) \in (\Coker^{T,\LL}\CR_J)_{(h,j)} = (\Coker^1\CR_J)_{(h_1,j_1)}\oplus(\Coker^2\CR_J)_{(h_2,j_2)} 
\end{equation*}
with 
\begin{eqnarray*}
  &&\eta_1 \in (\Coker^1\CR_J)_{(h_1,j_1)} = \ker(D^{\xi}_{h_1})^* \subset H^{1,2}(T^*\Si_1\otimes_{j_1,J_{\xi}} h_1^*\xi), \\ 
  &&\eta_2 \in (\Coker^2\CR_J)_{(h_2,j_2)} = \ker(D^{\xi}_{h_2})^* \subset H^{1,2}(T^*\Si_2\otimes_{j_2,J_{\xi}} h_2^*\xi)
\end{eqnarray*}
we define a preglued section 
\begin{equation*}
\eta^r_0 = \sharp^0_r \eta \;=\; \eta_1 \sharp^0_r \eta_2 \in H^{1,2}(T^*\Si^r\otimes_{j^r,J_{\xi}}(h^r)^*\xi)
\end{equation*}
in the bundle of $j^r,J_{\xi}$-antiholomorphic one-forms over the glued surface $(\Si^r,j^r)$ with values in the pull-back bundle $(h^r)^*\xi$. Note that the 
integration measure for defining the $H^{1,2}$-norm agrees on the connecting cylindrical neck 
$\psi_{21}((0,+r]\times S^1)\sharp_{\theta_{12}}\psi_{12}([-r,0)\times S^1)$ with the standard measure $ds\wedge dt$ on the cylinder. \\

For $r>0$ let $\beta^r: [0,+r]\to [0,1]$ be a smooth cut-off function such that $\beta^r(s)=1$ for $0\leq s\leq r/4$ and $\beta^r(s)=0$ for $3r/4\leq s\leq r$ 
with $|\partial_s \beta^r| \leq 4/r$. Let 
\begin{equation*} \beta^r_1, \beta^r_2: \Si^r\to[0,1] \end{equation*}
be the two cut-off functions which are constant equal to zero on $\Si^r_2$, $\Si^r_1$, constant equal to one on $\Si^r_1-\psi_{12}([-r,0]\times S^1)$, 
$\Si^r_2-\psi_{21}([0,+r]\times S^1)$ and are on $\psi_{12}([-r,0)\times S^1)\subset \Si^r_1$, $\psi_{21}((0,+r]\times S^1)\subset \Si^r_2$ given by 
\begin{eqnarray*} \beta_1^r(\psi_{12}(s,t))=\beta^r(-s),\;\; \beta_2^r(\psi_{21}(s,t))=\beta^r(+s), \end{eqnarray*}
respectively. With this we define the preglued section 
$\eta_1\sharp^0_r\eta_2$ on $\Si^r=\Si^r_1\sharp\Si^r_2$ by 
\begin{equation*}
\eta^r_0 = \eta_1\sharp^0_r\eta_2 = \beta^r_1\eta_1 + \beta^r_2\eta_2. 
\end{equation*}
It follows that $\eta^r_0$ agrees with $\eta_1$, $\eta_2$ over $\Si^r_1-\psi_{12}([-r,0]\times S^1)$, $\Si^r_2-\psi_{21}([0,+r]\times S^1)$, 
respectively, while over the connecting neck  
we have 
\begin{eqnarray*}
 (\eta^r_0\circ\psi_{12})(s,t) = \beta^r(-s)\cdot(\eta_1\circ\psi_{12})(s,t),\\
 (\eta^r_0\circ\psi_{21})(s,t) = \beta^r(+s)\cdot(\eta_2\circ\psi_{21})(s,t).
\end{eqnarray*}
Observe that by $\beta^r(s)=0$ for $3r/4\leq s\leq r$ this indeed yields a well-defined section in $H^{1,2}(T^*\Si^r\otimes_{j^r,J_{\xi}}(h^r)^*\xi)$. \\

{\it The gluing lemma.} For $r\in (0,+\infty)$, $\theta\in S^1$ let 
\begin{eqnarray*} 
 \sharp^0_r (\Coker^{T,\LL}\CR_J)_{(h,j)} &=& \ker (D^{\xi}_{h_1})^* \sharp^0_{r,\theta} \ker (D^{\xi}_{h_2})^* \\
  &=& \{\sharp^0_r(\eta_1,\eta_2): \eta_i\in\ker (D^{\xi}_{h_i})^*, i=1,2\} \\
  &\subset& H^{1,2}(T^*\Si^r\otimes_{j^r,J_{\xi}}(h^r)^*\xi), \\
\end{eqnarray*}
denote the subspaces of preglued sections. With the orthogonal projections 
\begin{eqnarray*}
 \pi_r: H^{1,2}(T^*\Si^r\otimes_{j^r,J_{\xi}}(h^r)^*\xi) &\to& \coker D_{h^r,j^r} = \ker (D^{\xi}_{h^r})^* \\
\end{eqnarray*} 
we can state and prove the gluing lemma: \\
\\
{\bf Lemma 2.6:} {\it The projections from the spaces of preglued sections on the fibres of the cokernel bundles over the underlying glued branched covers, 
\begin{eqnarray*} 
  \pi_r: \sharp^0_r (\Coker^{T,\LL}\CR_J)_{(h,j)} \to (\Coker\CR_J)_{(h^r,j^r)},\;\; (h^r,j^r)=\sharp_r(h,j,\theta) 
\end{eqnarray*}
are isomorphisms for all $r>0$ sufficiently large.} \\
\\
{\it Proof:} For the proof we follow the proof of proposition 3.2.9 in [Sch]. However we emphasize that we cannot directly apply the linear gluing lemma 
in [Sch], since the linear operator $D^{\xi}_{h^r}$ over the glued surface does not agree with the glued operator 
$D^{\xi}_{h_1,j_1}\sharp_{r,\theta} D^{\xi}_{h_2,j_2}$ studied in [Sch]. \\

Observe that it suffices to find for every $r>0$ sufficiently large a constant $c>0$ such that 
$\|(D^{\xi}_{h^r,j^r})^*\eta\|_2 \geq c\|\eta\|_{1,2}$ for all $\eta \in (\sharp^0_r \Coker^{T,\LL}\CR_J)_{(h,j)}^{\perp} = 
(\ker (D^{\xi}_{h_1})^* \sharp^0_{r,\theta} \ker (D^{\xi}_{h_2})^*)^{\perp}$. 
Indeed, it then follows that 
\begin{equation*} \ker (D^{\xi}_{h^r})^* \cap (\ker (D^{\xi}_{h_1})^* \sharp^0_{r,\theta} \ker (D^{\xi}_{h_2})^*)^{\perp} =\{0\},  
\end{equation*}
which proves that the orthogonal projection is surjective. On the other hand, since $\dim \ker D^{\xi}_{h^r,j^r} = \dim \ker D^{\xi}_{h_1,j_1} 
= \dim \ker D^{\xi}_{h_2,j_2}=0$ by proposition 2.3 and the index of $D^{\xi}_{h^r,j^r}$ equals the sum of the indices of 
$D^{\xi}_{h_1,j_1}$ and $D^{\xi}_{h_2,j_2}$, it follows that 
\begin{equation*} \dim \ker (D^{\xi}_{h^r})^* \;=\; \dim \ker (D^{\xi}_{h_1})^* + \dim \ker (D^{\xi}_{h_2})^*. \end{equation*}
Since the latter agrees with the dimension of the space $\ker (D^{\xi}_{h_1})^* \sharp^0_{r,\theta} \ker (D^{\xi}_{h_2})^*$ of preglued sections, the 
surjectivity of the orthogonal projection directly implies that it is an isomorphism. \\

Assume to the contrary that there exists a sequence 
\begin{equation*} 
\eta_n\in (\ker (D^{\xi}_{h_1})^* \sharp_{r_n,\theta}^0 \ker (D^{\xi}_{h_2})^*)^{\perp},\;\;r_n\to\infty
\end{equation*}
with $\|\eta_n\|_{1,2}=1$ but $\|(D^{\xi}_{h^{r_n}})^*\eta_n\|_2\to 0$ as $n\to\infty$.  
Now observe that 
\begin{eqnarray*} 
 \|(D^{\xi}_{h^{r_n}})^*(\beta^{r_n}_1\eta_n)\|_2 &\leq& \|(D^{\xi}_{h^{r_n}})^*\eta_n\|_2 + c_1\|d\beta^{r_n}_1 \cdot \eta_n\|_2\\
 &\leq& \|(D^{\xi}_{h^{r_n}})^*\eta_n\|_2 + c_1 \|d\beta^{r_n}_1\|_{\infty} \cdot \|\eta_n\|_2
\end{eqnarray*}
for some $c_1>0$ with $\|d\beta^{r_n}_1\|_{\infty}\leq 4/r_n$ and $\|\eta_n\|_2\leq \|\eta_n\|_{1,2}=1$, so that 
$\|(D^{\xi}_{h^{r_n}})^*(\beta^{r_n}_1\eta_n)\|_2\to 0$ for $n\to\infty$. But since 
$(h^{r_n},j^{r_n})\to (h_1,j_1)$ on $\Si_1^{r_n} = \Si_1-\psi_{12}((-\infty,-r_n)\times S^1)$, this directly implies that 
\begin{equation*}
  \|(D^{\xi}_{h_1})^*(\beta^{r_n}_1\eta_n)\|_2 \to 0 
\end{equation*}
in the $L^2(\Si_1)$-sense and we can use the semi-Fredholm property of $(D^{\xi}_{h_1})^*$ and the boundedness of $(\eta_n)$ to deduce that, 
possibly after passing to a suitable subsequence, 
\begin{equation*}
\beta^{r_n}_1\eta_n \stackrel{H^{1,2}}{\to} \eta_1,\;\; \eta_1\in\ker (D^{\xi}_{h_1})^*.
\end{equation*}
Using the same arguments we deduce $\beta^{r_n}_2\eta_n \to \eta_2 \in\ker (D^{\xi}_{h_2})^*$. We use this to prove the desired contradiction 
by computing   
\begin{eqnarray*}
 1 = \lim_{n\to\infty}\|\eta_n\|_{1,2} &=& \lim_{n\to\infty} \<(\beta^{r_n}_1)^2\eta_n + (\beta^{r_n}_2)^2\eta_n,\eta_n\>_{1,2} \\
   &+& \lim_{n\to\infty} \<(1-(\beta^{r_n}_1)^2 - (\beta^{r_n}_2)^2)\cdot \eta_n,\eta_n\>_{1,2} \\
   &=& \lim_{n\to\infty} \<\beta^{r_n}_1 \eta_1 + \beta^{r_n}_2 \eta_2,\eta_n\>_{1,2} \\
   &=& \lim_{n\to\infty} \<\eta_1 \sharp_{r_n,\theta}^0 \eta_2,\eta_n\>_{1,2} = 0, 
\end{eqnarray*}
since $\eta_n\in (\ker (D^{\xi}_{h_1})^* \sharp_{r_n,\theta}^0 \ker (D^{\xi}_{h_2})^*)^{\perp}$, where it only remains to prove that 
\begin{equation*} \lim_{n\to\infty} \<(1-(\beta^{r_n}_1)^2 - (\beta^{r_n}_2)^2)\cdot \eta_n,\eta_n\>_{1,2} = 0. \end{equation*}
For this we use that $1-(\beta^{r_n}_1)^2 - (\beta^{r_n}_2)^2$ has only support in the middle part 
\begin{equation*}
\psi_{21}([+r_n/4,+r_n]\times S^1)\sharp_{\theta_{12}}\psi_{12}([-r_n,-r_n/4]\times S^1) \cong [-3r_n/4,+3r_n/4]\times S^1 
\end{equation*} 
of the cylindrical neck to prove that the $H^{1,2}$-norm of $(1-(\beta^{r_n}_1)^2 - (\beta^{r_n}_2)^2)\eta_n$ tends to zero as $n\to\infty$: \\

Choosing a unitary trivialization of the symplectic hyperplane bundle $\xi$ over the simple orbit $\gamma$, the restriction of the the differential 
operator $(D^{\xi}_{h^r})^*$ to $[-3r_n/4,+3r_n/4]\times S^1 \subset \Si^r$ is of the form 
\begin{eqnarray*}
 D_n = \partial_s + J_0\partial_t + S_n: 
 && H^{1,2}([-3r_n/4,+3r_n/4]\times S^1,\IR^{2m-2})\\
 &\to& L^2([-3r_n/4,+3r_n/4]\times S^1,\IR^{2m-2})
\end{eqnarray*}
with $S_n(s,t)\in\IR^{(2m-2)\times(2m-2)}$, which we extend to an operator on the full cylinder $\IR\times S^1$ by setting $S_n(+s,t)=S_n(+3r_n/4,t)$, 
$S_n(-s,t)=S_n(-3r_n/4,t)$ for $s>3r_n/4$. 
In order to study the operator $D_n$ let $h_n = h^{r_n}|_{[-3r_n/4,+3r_n/4]\times S^1}: [-3r_n/4,+3r_n/4]\times S^1 \to \RS$ and $x_n=h_n(0,\cdot): S^1\to\RS$. 
Since for $n\to\infty$ the length of the cylindrical neck goes to infinity, it follows that $h_n$ converges on each compact subinterval uniformly with all 
derivatives to the $\IR$-independent function $x_{\infty} = \lim_{n\to\infty}x_n: S^1\to\RS$ of the form $x_{\infty}(t)=(s_0,m_{12}t+t_0)$. From this 
it follows that $S_n(s,t)\to S_{\infty}(t)$, i.e., $D_n$ is converging in the operator norm 
to a translation invariant operator $D_{\infty}$. \\

Finishing the proof observe that from  
\begin{eqnarray*}
 && \|D_n (\eta_n-(\beta^{r_n}_1)^2\eta_n - (\beta^{r_n}_2)^2\eta_n)\|_2 \\
 && \leq \|D_n\eta_n\|_2 + c_2(\|d\beta^{r_n}_1\|_{\infty}\|\beta^{r_n}_1\eta_n\|_2+\|d\beta^{r_n}_2\|_{\infty}\|\beta^{r_n}_2\eta_n\|_2)\\
 && \leq \|D_n\eta_n\|_2 + c_2(\|d\beta^{r_n}_1\|_{\infty}+\|d\beta^{r_n}_2\|_{\infty})\|\eta_n\|_2
\end{eqnarray*}
and $\|d\beta^{r_n}_1\|_{\infty},\|d\beta^{r_n}_2\|_{\infty}\to 0$, $\|\eta_n\|_2=1$, $\|D_n\eta_n\|_2\to 0$ it follows that 
\begin{equation*}
\|D_{\infty} (\eta_n-(\beta^{r_n}_1)^2\eta_n - (\beta^{r_n}_2)^2\eta_n)\|_2 \to 0,\;\;n\to\infty.
\end{equation*}
But now we can use the fact that the operator $D_{\infty}: H^{1,2}(\RS,\IR^{2m-2})\to L^2(\RS,\IR^{2m-2})$ is an isomorphism ([Sch]) and hence   
\begin{equation*}
  \|(1-(\beta^{r_n}_1)^2 - (\beta^{r_n}_2)^2)\eta_n\|_{1,2} \leq c_3 \cdot \|D_{\infty} (\eta_n-(\beta^{r_n}_1)^2\eta_n - (\beta^{r_n}_2)^2\eta_n)\|_2 
\end{equation*}
for some $c_3>0$ to deduce that $\|(1-(\beta^{r_n}_1)^2 - (\beta^{r_n}_2)^2)\eta_n\|_{1,2}\to 0$ as $n$ goes to infinity. $\qed$\\

\subsection{Orientation}

In this section we show how the techniques by [BM] and [HWZ] for defining coherent orientations of the 
moduli spaces in symplectic field theory define an orientation of the cokernel bundle $\Coker\CR_J$ over the non-compactified moduli space 
$\IM$ and discuss the extension over the boundary strata. Although we have seen in the last section that the cokernel bundle $\Coker\CR_J$ 
naturally lives over the 
quotient $\IM/(\IZ_{m^+}\times\IZ_{m^-})$, obtained by forgetting the asymptotic markers $\mu^{\pm}\in\IZ_{m^{\pm}}$, 
we show in this section that in general we can orient $\Coker\CR_J$ only over the full moduli space $\IM$. 
For this we start with recalling the main points of the constructions of coherent orientations in [BM]: \\

Let $\Si=S^2-\{z^{\pm}_{1,0},...,z^{\pm}_{n^{\pm},0}\}$ denote as in 2.2 the punctured sphere underlying the moduli space $\IM$. For regular paths of 
symplectic matrices 
\begin{equation*}
 A^{\pm}_1,...,A^{\pm}_{n^{\pm}}: [0,1]\to\Sp(2m-2),\;\; \det(A^{\pm}_k(1)-A^{\pm}_k(0))\neq 0
\end{equation*} 
with $A^{\pm}_k(0)=\Id$, $\dot{A}^{\pm}_k(0)A^{\pm}_k(0)^{-1}=\dot{A}^{\pm}_k(1)A^{\pm}_k(1)^{-1}$ and where $2m-2$ is the rank of $\xi$, 
let $\IO((\Si,j,\mu^{\pm}),(A^{\pm}_k)_{k=1}^{n^{\pm}})$ denote the set of Cauchy-Riemann operators
\begin{eqnarray*}
 D: &&H^{1,p,d}_{\cst}(\Si,\IC)\oplus H^{1,p}(\Si,\IR^{2m-2})\\
    &&\to L^{p,d}(T^*\Si\otimes_{j,i}\IC) \oplus L^p(T^*\Si\otimes_{j,J_0}\IR^{2m-2})  \\
    &&D\cdot v = dv + J_0\cdot dv\cdot j + S\cdot v
\end{eqnarray*}
where $S:\Si\to\IR^{2m\times 2m}$ is a family of symmetric matrices such that the limit matrices are of the form 
\begin{equation*} 
 (S\circ\psi^{\pm}_k)(s,t+\mu^{\pm}_k)\stackrel{s\to\pm\infty}{\longrightarrow} (S\circ\psi^{\pm}_k)(\pm\infty,t+\mu^{\pm}_k) = 
 \binom{0\;\;\;\;\;\;\;\;0}{0\;\;S^{\pm}_k(t)},
\end{equation*}
and where $S^{\pm}_1,...,S^{\pm}_{n^{\pm}}: S^1\to\IR^{(2m-2)\times(2m-2)}$ are related to $A^{\pm}_1,...,A^{\pm}_{n^{\pm}}: [0,1]\to\Sp(2m-2)$ via 
\begin{equation*}
S^{\pm}_k(t)=-J_0\cdot \dot{A}^{\pm}_k(t)\cdot A^{\pm}_k(t)^{-1}
\end{equation*}
for all $k=1,...,n^{\pm}$.  \\

Since every operator $D\in \IO(\Si,(A^{\pm}_k)_{k=1}^{n^{\pm}})=\bigcup_{j,\mu}\IO((\Si,j,\mu^{\pm}),(A^{\pm}_k)_{k=1}^{n^{\pm}})$ 
is a Fredholm operator, we have the determinant line bundle $\Det(\Si,(A^{\pm}_k)_{k=1}^{n^{\pm}}))$ over $\IO$ with fibre  
\begin{equation*} 
  \Det(\Si,(A^{\pm}_k)_{k=1}^{n^{\pm}}))_D =\Det(D) = \Lambda^{\max} \ker D \otimes \Lambda^{\max} \coker D.
\end{equation*}
Since the space of Fredholm operators $\IO((\Si,j,\mu^{\pm}),(A^{\pm}_k)_{k=1}^{n^{\pm}})$ is contractible, it follows that the restriction  
$\Det((\Si,j,\mu^{\pm}),(A^{\pm}_k)_{k=1}^{n^{\pm}})$ of $\Det(\Si,(A^{\pm}_k)_{k=1}^{n^{\pm}}))$ to $\IO((\Si,j,\mu^{\pm}),(A^{\pm}_k)_{k=1}^{n^{\pm}})$ 
is trivial. On the other hand, it is shown in proposition 11 in [BM] that the determinant 
line bundle remains trivial when we allow the complex structure $j$ on the punctured sphere $\Si$ to vary. \\

In [BM] the authors describe a method to orient how the resulting bundles $\Det((\Si,\mu^{\pm}),(A^{\pm}_k)_{k=1}^{n^{\pm}})$ over 
$\IO(\Si,\mu^{\pm},(A^{\pm}_k)_{k=1}^{n^{\pm}})=\bigcup_j\IO((\Si,j,\mu^{\pm}),(A^{\pm}_k)_{k=1}^{n^{\pm}})$ for any number of punctures, 
directions $\mu^{\pm}$ and regular paths 
$A^{\pm}_1,...,A^{\pm}_{n^{\pm}}$ of symplectic matrices. The construction is based on arbitrarily fixing orientations for determinant bundles  
over the space $\IO((\IC^*,0),A)$ of Cauchy-Riemann operators on the holomorphic plane, constructing a gluing map for determinant bundles under 
gluing of Riemann surfaces and finally observing that we have a natural orientation of $\Det(S^2)$ induced by the complex orientation 
of the determinant line over the standard Cauchy-Riemann operator on $(S^2,i)=\CP$. Note that at this point the specification of the 
asymptotic markers $\mu=(\mu^+,\mu^-)$, $\mu^{\pm}=(\mu^{\pm}_k)_{k=1}^{n^{\pm}}$ becomes important, as they describe how to glue 
the holomorphic planes to the punctured sphere $\Si$ to obtain the closed sphere $S^2$. However it directly follows from the construction that 
the orientations on $\Det((\Si,\mu^{\pm}),(A^{\pm}_k)_{k=1}^{n^{\pm}})$ for different asymptotic markers $\mu$ fit together to give an orientation of 
the whole determinant bundle $\Det(\Si,(A^{\pm}_k)_{k=1}^{n^{\pm}}))$. \\  

Observe that the linearization of $\CR_J$ at some $(h,j,\mu^{\pm})\in\IM$,
\begin{eqnarray*}
 &D_{h,j}:& H^{1,p,d}_{\cst}(\Si,\IC)\oplus H^{1,p}(h^*\xi)\oplus T_j\IM_{0,n} \\
 &&\to L^{p,d}(T^*\Si\otimes_{j,i}\IC) \oplus L^p(T^*\Si\otimes_{j,J_{\xi}}h^*\xi),\\
\end{eqnarray*}
can be written as sum $D_{h,j}= D_h + D_j$ with
\begin{eqnarray*}
  &D_j:& T_j\IM_{0,n}\to L^{p,d}(T^*\Si\otimes_{j,i}\IC) \oplus L^p(T^*\Si\otimes_{j,J_{\xi}}h^*\xi)\\
  &D_h:& H^{1,p,d}_{\cst}(\Si,\IC)\oplus H^{1,p}(h^*\xi) \to L^{p,d}(T^*\Si\otimes_{j,i}\IC) \oplus L^p(T^*\Si\otimes_{j,J_{\xi}}h^*\xi)
\end{eqnarray*}
where $D_h$ is a Cauchy-Riemann operator. Using a unitary trivialization of the hyperplane bundle $\xi$ over the closed simple orbit $\gamma$, 
we get a unitary trivialization of $h^*\xi$ and a natural map 
\begin{equation*}
  \op: \IM\to\IO(\Si,(A^{\pm}_k)_{k=1}^{n^{\pm}}),\;\;(h,j,\mu^{\pm})\mapsto D_h,
\end{equation*}
where the regular paths of symplectic matrices $A^{\pm}_1,...,A^{\pm}_{n^{\pm}}$ are determined by the restriction to $\xi$ of the 
linearized Reeb flow along $\gamma$. Using the map $\op$ we can pull-back the determinant 
bundle $\Det=\Det(\Si,(A^{\pm}_k)_{k=1}^{n^{\pm}})$ to obtain the line bundle $\op^*\Det$ over $\IM$. On the other hand, 
following the arguments in [BM], we deduce from the fact that $D_{h,j}=D_j\oplus D_h$ is homotopic to the stabilization $0\oplus D_h$ 
with the complex vector space $T_j\IM_{0,n}$ that the determinant spaces of the linearization $D_{h,j}$ and the Cauchy-Riemann operator $D_h$ 
are canonically isomorphic, so that the pull-back of the determinant bundle over the space of Cauchy-Riemann operators is isomorphic to the 
determinant bundle of the fully linearized operator 
\begin{equation*}
\op^*\Det \cong \Lambda^{\max} \Ker\CR_J \otimes \Lambda^{\max} \Coker\CR_J
\end{equation*}
with fibre $\Lambda^{\max} \ker D_{h,j} \otimes \Lambda^{\max} \coker D_{h,j}$ over $(h,j,\mu^{\pm})\in\IM$.\\

Since $\Ker\CR_J$ and $\Coker\CR_J$ are bundles over $\IM/(\IZ_{m^+}\times\IZ_{m^-})$, it follows that the action of $\IZ_{m^+}\times\IZ_{m^-}$ 
lifts in an obvious way to an action on the vector bundle $\Lambda^{\max} \Ker\CR_J \otimes \Lambda^{\max} \Coker\CR_J$ which is trivial 
on the fibres. On the other hand, the fibres over $(h,j,\mu^{\pm}),(h,j,\mu'^{\pm})\in\IM$ do not neccessarily carry the same orientation. Indeed it 
is shown in theorem 3 in [BM] that this action is orientation-preserving if $\gamma$ is good, else, the action is orientation-preserving or 
-reversing if $\mu'-\mu\in\IZ_{m^+}\times\IZ_{m^-}$ is even or odd, respectively. In this case the even iterates $\gamma^{2k}$ of the simple orbit $\gamma$ are 
called bad. \\
\\
{\bf Proposition 2.7:} {\it For every tree with level structure $(T,\LL)$ with trees $T_1$,..., $T_L$, the choice of coherent orientations in [BM] 
equip the cokernel bundles $\Coker^{T_1}\CR_J$, ..., $\Coker^{T_L}\CR_J$ over $\IM_{T_1}$, ..., $\IM_{T_L}$ with orientations, which descend to 
an orientation of the cokernel bundle $\Coker^{T,\LL}\CR_J = \pi_1^*\Coker^{T_1}\CR_J\oplus...\oplus\pi_L^*\Coker^{T_L}\CR_J$ over 
$\IM_{T,\LL}=\IM_{T_1}\times ...\times \IM_{T_L}/\Delta$. The orientations of the cokernel bundles over the strata $\IM_{T,\LL}\subset\CM$ 
in general do} not {\it fit together to an orientation of the cokernel bundle $\overline{\Coker}\CR_J$ over the compactified moduli space $\CM$, but 
differ by a fixed sign due to reordering the punctures.} \\

We remark that the fact that the orientations of the cokernel bundles over the different strata differ by a fixed sign is not completely trivial, 
since the strata are in general not connected due to the possible choices for the asymptotic markers. Furthermore it directly follows from 
theorem 3 in [BM] that the cokernel bundle $\Coker\CR_J$ is orientable over the quotient $\IM/(\IZ_{m^+}\times\IZ_{m^-})$ only when all 
asymptotic orbits $\gamma^{m^{\pm}_1},...,\gamma^{m^{\pm}_{n^{\pm}}}$ are good. \\  
\\
{\it Proof:}  In the way described above the choice of coherent orientations in symplectic field theory following [BM] provides us with an orientation of the 
determinant bundles $\Lambda^{\max} \Ker\CR_J \otimes \Lambda^{\max} \Coker\CR_J$ of the Cauchy-Riemann operator $\CR_J$ over 
the moduli space of branched covers $\IM$.  But since by lemma 2.3 $\Ker\CR_J$ agrees with the tangent space to $\IR\times\IM$ and 
$\IR\times \IM = \RS\times \IM_{0,n} \times \IZ_{m^+}\times\IZ_{m^-}$ is a complex manifold, we always have a natural orientation
of $\Ker\CR_J$, which directly fixes an orientation on the cokernel bundle $\Coker\CR_J$ over $\IM$ by requiring that 
the orientations on $\Ker\CR_J$ and $\Coker\CR_J$ determine the orientation of the determinant bundle 
$\Lambda^{\max} \Ker\CR_J \otimes \Lambda^{\max} \Coker\CR_J$. \\

In order to see that the same arguments can be used to orient the cokernel bundles $\Coker^{T_{\ell}}\CR_J$ over the moduli spaces 
$\IM_{T_{\ell}}$ of nodal curves for $\ell=1,...,L$, observe that the constructions in [BM] immediately generalize to nodal curves 
in such a way that the orientation of the determinant bundle for the nodal surface fits with the orientation for the determinant 
bundle over the glued surface. 
Indeed this follows, using the gluing argument for the determinant line bundles, simply from the fact that also on closed surface with nodes we have a 
standard Cauchy-Riemann operator providing us with a natural orientation of the determinant line bundle over the space of Fredholm operators on a closed 
nodal surface, which clearly fits with the natural orientation of the determinant bundle over the space of Fredholm operators over the glued surface. 
In order to see that the orientations of $\Coker^{T_1}\CR_J$, ..., $\Coker^{T_L}\CR_J$ determine an orientation of the cokernel bundle over 
the stratum $\IM_{T,\LL}=\IM_{T_1}\times ...\times \IM_{T_L}/\Delta$, we must show that the lift of the action of $\Delta$ on 
$\IM_{T_1}\times ...\times \IM_{T_L}$ to the cokernel bundle $\Coker^{T,\LL}\CR_J= \pi_1^*\Coker^{T_1}\CR_J\oplus...\oplus\pi_L^*\Coker^{T_L}\CR_J$ is 
orientation-preserving: \\

For this recall that $\Delta=\prod_{\LL(\alpha)>\LL(\beta)}\Delta_{\alpha\beta}$, where $\Delta_{\alpha\beta}$ is the diagonal in 
$\IZ_{|m_{\alpha\beta}|} \times \IZ_{|m_{\beta\alpha}|}$ so that $\Delta_{\alpha\beta}$ acts on $\IM_{T_k}\times\IM_{T_{\ell}}$ for 
$k=\LL(\alpha)$, $\ell=\LL(\beta)$. Now it follows from theorem 3 in [BM] that the $\IZ_{|m_{\alpha\beta}|}$-actions on the cokernel bundles 
$\Coker^{T_k}\CR_J$ and $\Coker^{T_{\ell}}\CR_J$ are orientation-preserving if $\gamma^{|m_{\alpha\beta}|}$ is good, and simultaneuously 
orientation-preserving or -reversing for even or odd elements in $\IZ_{|m_{\alpha\beta}|}$ if $\gamma^{|m_{\alpha\beta}|}$ is bad. Hence 
the action on the direct sum $\pi_k^*\Coker^{T_k}\CR_J\oplus\pi_{\ell}^*\Coker^{T_{\ell}}\CR_J$ is orientation-preserving in all cases. \\

The statement about the behaviour of the orientations on the cokernel bundles under gluing directly follows from theorem 1 in [BM] which states 
that the gluing diffeomorphisms preserve the orientations up to a sign due to reordering of the punctures. This is however an immediate consequence of the 
behaviour of the orientation of moduli spaces under reordering the punctures. $\qed$ \\ 

\section{Perturbation theory and Euler numbers}

\subsection{Perturbed Cauchy-Riemann operator}

As outlined in the section about the linearized operator, the Cauchy-Riemann operator $\CR_J$ can be 
viewed as a smooth section in a Banach space bundle $\EE^{p,d}$ over a Banach manifold of maps $\BB^{p,d}$. Since for the 
contribution to the differential in contact homology and rational symplectic field theory we are interested in moduli spaces 
of branched covers $\IM$ of virtual dimension zero while the actual dimension is always strictly greater than zero, it follows 
that in the cases of interest the Cauchy-Riemann operator $\CR_J$ does not meet the zero section transversally. In other words, the 
image bundle $\Im\CR_J$ of $\CR_J$ over $\IM$ with fibre $(\Im\CR_J)_{h,j}=\im D_{h,j}$ is a true closed
subbundle of the Banach space bundle $\EE^{p,d}$ over the moduli space of branched covers $\IM =\CR_J^{-1}(0) \subset \BB^{p,d}$, where the 
closedness of $\im D_{h,j}$ in $\EE^{p,d}_{h,j}$ follows from the semi-Fredholm property of $D_{h,j}: T_{h,j}\BB^{p,d}\to\EE^{p,d}_{h,j}$. 
In particular, observe that we have a natural splitting 
\begin{equation*}
\EE^{p,d}|_{\CR_J^{-1}(0)} = \Im\CR_J \oplus \Coker\CR_J 
\end{equation*}
with the cokernel bundle $\Coker\CR_J$ introduced in section two. \\

For determining the contribution of $\IM$ to the differential in contact homology and rational symplectic field theory it follows 
that the Cauchy-Riemann operator $\CR_J$ has to be perturbed slightly to a transversal section in the Banach space bundle $\EE^{p,d}\to\BB^{p,d}$ 
in the sense that it meets the zero section transversally. This means that we have to add a compact perturbation $\nu$ to the Cauchy-Riemann operator 
to make it transversal and count elements in the regular moduli space $\IM^{\nu}$, 
\begin{equation*}
 \IM^{\nu}=(\CR^{\nu}_J)^{-1}(0)\subset \BB^{p,d},\;\;\CR^{\nu}_J=\CR_J+\nu.
\end{equation*}

We first prove the folk's theorem that it indeed suffices to study smooth sections in the cokernel bundle $\Coker\CR_J\subset \EE^{p,d}|_{\IM}$ over 
the moduli space $\IM$, i.e., the zero set of the Cauchy-Riemann operator $\CR_J$; for a different proof in the context of Gromov-Witten theory, see 
proposition 7.2.3 in [MDSa]. For this we extend a section in $\Coker\CR_J$ over $\IM$ to a smooth section in the Banach space bundle $\EE^{p,d}$ 
over the whole Banach manifold $\BB^{p,d}$ as follows: \\

Choosing a unitary trivialization of $(\xi,\omega_{\xi},J_{\xi})$ along the Reeb orbit $\gamma$, note that it can be extended to a unitary trivialization 
of $(\xi,\omega_{\xi},J_{\xi})$ over a sufficiently small neighborhood $N$ of $\gamma$ using parallel transport along geodesics with respect to a unitary connection 
$\nabla$. Further identifying $N$ with a neighborhood of the zero section in $\gamma^*\xi\cong S^1\times\IC^{m-1}$ we assume that $N\cong S^1\times B_{\eps}(0)$ with 
$B_{\eps}(0)=\{z\in\IC^{m-1}: |z|<\eps\}$. \\

Now observe that for a section $\nu$ in the cokernel bundle $\Coker\CR_J$ over $\IM$ we have 
$\nu(h,j)\in L^p(T^*\Si\otimes_{j,i}\IC^{m-1})$ for every $(h,j)\in\IM/(\IZ_{m^+}\times\IZ_{m^-})$, which for every tuple $(h,j,z)$, $z\in\Si$ defines an element 
$\nu(h,j,z)\in T_z^*\Si\otimes_{j,i}\IC^{m-1}$. Identifying for fixed complex structure $j$ the branched covering map $h$ with the direction $t\in S^1$ of the 
asymptotic marker, i.e., $(h,j)\equiv (t,j)\in S^1\times\IM_{0,n}\cong\IM/(\IZ_{m^+}\times\IZ_{m^-})$, note that this defines for every $(j,z)$ a smooth map 
$\nu_0(j,z): S^1\to T_z^*\Si\otimes_{j,i}\IC^{m-1}$. 
With the choice of a smooth cut off function $\varphi_{\eps}:[0,\eps]\to [0,1]$ with $\varphi_{\eps}(0)=1$ and $\varphi_{\eps}(\eps)=0$ we can extend $\nu_0(j,z)$ to 
a map starting from $N\cong S^1\times B_{\eps}(0)$ by setting $\nu_0(j,z)(t,v):=\varphi_{\eps}(|v|)\cdot \nu_0(j,z)(t)$ for $(t,v)\in S^1\times B_{\eps}(0)$. \\
  
Let $\UU=\UU^{p,d}$ denote the small neighborhood of $\IM$ in $\BB^{p,d}$ of all maps $u:\Si\to\IR\times V$ having image contained in the neighborhood 
$N$ of $\gamma$ in $V$. Writing $u=(h,v):\Si\to(\RS)\times B_{\eps}(0)\subset\IR\times V$ we can define an extension of $\nu$ from $\IM$ to $\UU^{p,d}$ by 
setting $\nu(u,j)(z):= \nu_0(j,z)(t(u),v(z))$ with $t(u)\in S^1$ denoting the direction of the asymptotic marker defined by the map $u$. Note that this indeed defines 
an extension and that $\nu(u,j)\in L^p(T^*\Si\otimes_{j,i}\IC^{m-1})$. In particular, $\nu$ defines a section in trivial bundle $\EE^{p,d}|_{\UU^{p,d}}$ with fibre 
$L^{p,d}(T^*\Si\otimes_{j,i}\IC)\oplus L^p(T^*\Si\otimes_{j,i}\IC^{m-1})$, which in turn after extending by zero defines a section in the Banach space bundle 
$\EE^{p,d}$ over the whole Banach manifold $\BB^{p,d}$. Then the following holds: \\
\\
{\bf Proposition 3.1:} {\it Let $\nu$ be a section in the cokernel bundle $\Coker\CR_J\subset \EE^{p,d}|_{\IM}$ over the moduli space 
$\IM=\CR_J^{-1}(0)\subset\BB^{p,d}$, which is extended to a section in $\EE^{p,d}$ as described above. Then it holds:} 
\begin{itemize} 
\item {\it The moduli space $\IM^{\nu}$ agrees with the zero set of $\nu$ in $\IM$,}   
\begin{equation*} \IM^{\nu}=\{(h,j)\in\IM: \nu(h,j)=0\}. \end{equation*}
\item {\it If $\nu$ is a transversal section in $\Coker\CR_J$, then $\CR^{\nu}_J$ is a transversal section in $\EE^{p,d}$, i.e., $\IM^{\nu}$ is regular.}
\item {\it The linearization of $\nu$ at every zero is a compact operator, so that the linearizations of $\CR_J$ and $\CR_J^{\nu}$ belong to 
the same class of Fredholm operators.}
\end{itemize}
$ $\\
{\it Proof:} First we find no zeroes of $\CR_J^{\nu}$ outside of the neighborhood $\UU$ of $\IM$ since there $\CR_J^{\nu}=\CR_J$. For every $(u,j)\in\UU$ 
with $u=(h,v):\Si\to(\RS)\times B_{\eps}(0)$ let $\pi_1$ denote the projection onto the first factor in 
$\EE^{p,d}_{u,j}=L^{p,d}(T^*\Si\otimes_{j,i}\IC)\oplus L^p(T^*\Si\otimes_{j,i}\IC^{m-1})$. Then we 
have by construction that $\pi_1\circ\nu(u,j)=0$ while $\pi_1\circ\CR_J(u)=\CR h$ with the standard Cauchy-Riemann operator $\CR: H^{1,p,d}_{\cst}(\Si,\IC) 
\to L^{p,d}(T^*\Si\otimes_{j,i}\IC)$. For $(u,j)\in \UU-\IM$ it follows that $\pi_1\circ \CR_J^{\nu}(u) = \CR h\neq 0$, so that we find no zeroes of $\CR_J^{\nu}$ 
in $\UU-\IM$. Finally, on $\IM$ we have $\CR_J^{\nu}=\nu$. \\ 

With respect to the splittings $T_{h,j}\BB^{p,d} = (\ker D_{h,j})^{\perp}\oplus \ker D_{h,j}$ and $\EE^{p,d}_{h,j} = \im D_{h,j} 
\oplus \coker D_{h,j}$ at $(h,j)\in\IM$, observe that the linearization $D_{h,j}\CR_J^{\nu}: T_{h,j}\BB^{p,d}\to\EE^{p,d}_{h,j}$ at a zero $\nu(h,j)=0$, 
$(h,j)\in\IM$, is of upper triangle form
\begin{equation*}
 D_{h,j}\CR_J^{\nu} = \binom{D_{h,j}\;\;\;\;\;\;\;0}{D^0_{h,j}\nu\;\;D^1_{h,j}\nu}, 
\end{equation*}
where $D^0,D^1$ denotes differentiation in the direction of $(\ker D_{h,j})^{\perp}$ and $\ker D_{h,j}=T_{h,j}\IM$, respectively. Since $\nu$ is a transversal 
section in $\Coker\CR_J$ over $\IM$ precisely when $D^1_{h,j}\nu: \ker D_{h,j}\to\coker D_{h,j}$ is surjective 
at every $\nu(h,j)=0$, the second statement follows from the fact that $D_{h,j}: (\ker D_{h,j})^{\perp}\to \im D_{h,j}$ is an isomorphism. \\

For the last statement it suffices to see that the linearization $D_{h,j}\nu: T_{h,j}B^{p,d}\to\coker D_{h,j}$ is an 
operator with finite-dimensional image. $\qed$\\

Since the cokernel bundle $\Coker\CR_J$ as well as its base space $\IM$ are oriented, it follows that the regular moduli space $\IM^{\nu}$ carries an 
orientation, which by the construction of orientations for $\Coker\CR_J$ agrees with the orientation of moduli spaces in symplectic field theory 
constructed in [BM]. 
The contribution of branched covers of orbit cylinders to the differential in rational symplectic field theory is given by the algebraic count of 
elements in $\IM^{\nu}$, 
which however might explicitly depend on the chosen perturbation $\nu$. \\

\subsection{Gluing compatibility}

In order to have transversality for all moduli spaces of connected branched covers without nodes we choose transversal sections 
$\nu=\nu_{\vec{m}}$ in the cokernel bundles over the moduli spaces 
\begin{equation*}
 \IM=\IM_{\vec{m}}= \IM_{0,0}(\gamma^{m^+_1},...,\gamma^{m^+_{n^+}};\gamma^{m^-_1},...,\gamma^{m^-_{n^-}})/\IR
\end{equation*}
for all tuples $\vec{m}=(\vec{m}^+,\vec{m}^-)$, $\vec{m}^{\pm}=(m^{\pm}_1,...,m^{\pm}_{n^{\pm}})$ with 
\begin{equation*}
 |\vec{m}|:= m^+_1+...+m^+_{n^+} = m^-_1+...+m^-_{n^-},
\end{equation*}
i.e., for which $\IM_{\vec{m}}\neq\emptyset$. \\

To be precise we choose transversal sections $\nu$ in the cokernel bundles over the quotient $\IM/(\IZ_{m^+}\times\IZ_{m^-})$, where we forget the position of 
the asymptotic markers at the positive, respectively negative punctures. In this way we ensure that the chosen abstract perturbation and hence the contribution 
of curves to the differential does not depend on the choice of asymptotic markers, which is implicit in the algebraic setup of symplectic field theory. 
At this point recall that although the cokernel bundle $\Coker\CR_J$ naturally lives over the quotient $\IM/(\IZ_{m^+}\times\IZ_{m^-})$, it is in general 
only orientable over the complete moduli space $\IM$, since the orientation of a fibre in general depends on the choice of the asymptotic markers at the 
punctures. \\

In order to have the compactness and gluing results for the resulting regular moduli spaces which are implicit in the definition of 
algebraic invariants in symplectic field theory we consider only sets of cokernel sections $(\nu_{\vec{m}})_{\vec{m}}$, which are compatible 
with compactness and gluing in symplectic field theory in the following sense: \\

Let $(h^q,j^q)$, $q\in\IN$ be a sequence of curves in the regular moduli space $\IM^{\nu}=\IM^{\nu_{\vec{m}}}_{\vec{m}}$, which converges for $q\to\infty$ 
to a level branched covering $(h,j)\in\IM_{T,\LL}\subset \CM$ with 
\begin{equation*} \IM_{T,\LL}= \IM_{T_1}\times ...\times\IM_{T_L}/\Delta \end{equation*}
and $\IM_{T_{\ell}}=\IM_{T_{\ell,1}}\times ...\times \IM_{T_{\ell,N_{\ell}}}\times\IR^{N_{\ell}-1}$. Then all components 
$(h_{\ell,k},j_{\ell,k})\in\IM_{T_{\ell,k}}$, $\ell=1,...,L$, $k=1,...,N_{\ell}$ again satisfy a perturbed Cauchy-Riemann equation. When the moduli space $\IM_{T,\LL}$ 
is made up of curves with no nodes, i.e., for which the trees $T_{1,1},...,T_{L,N_L}$ are trivial, the moduli spaces $\IM_{\ell,k}=\IM_{T_{\ell,k}}$ are 
again moduli spaces of connected branched covers without nodes, 
\begin{equation*} \IM_{\ell,k}=\IM_{\vec{m}_{\ell,k}} \end{equation*}
for new tuples $\vec{m}_{\ell,k}=(\vec{m}^+_{\ell,k},\vec{m}^-_{\ell,k})$, $\vec{m}^{\pm}_{\ell,k}=(m^{\pm}_{\ell,k,1},...,m^{\pm}_{\ell,k,n^{\pm}_{\ell,k}})$. 
Assuming that the abstract perturbations $\nu_{\ell,k}=\nu_{\vec{m}_{\ell,k}}$ for the moduli spaces $\IM_{\ell,k}$ are already chosen, 
compatibility with gluing in symplectic field theory now means that the abstract perturbation $\nu=\nu_{\vec{m}}$ is chosen in such a way that 
$(h_{\ell,k},j_{\ell,k})$ 
satisfies the Cauchy-Riemann equation with perturbation $\nu_{k,\ell}$, i.e., is an element in the regular moduli space 
$\IM^{\nu_{\ell,k}}_{\ell,k}\subset \IM_{\ell,k}$. Observe that for every level $\ell=1,...,L$ we have 
\begin{equation*} |\vec{m}_{\ell,1}| + ... + |\vec{m}_{\ell,N_{\ell}}| = |\vec{m}|, \end{equation*} 
while for the number of punctures $\sharp\vec{m}=n=n^++n^-$ and $\sharp\vec{m}_{\ell,k}=n_{\ell,k}=n^+_{\ell,k}+n^-_{\ell,k}$ we have  
\begin{equation*} \sharp\vec{m}_{\ell,k} < \sharp\vec{m}. \end{equation*}
It follows that the choice of the abstract perturbation $\nu_{\vec{m}}$ depends only on abstract perturbations $\nu_{\vec{m'}}$ with 
$\sharp\vec{m'}<\sharp\vec{m}$ and $|\vec{m'}|\leq |\vec{m}|$. \\

The correct setup for constructing perturbations $\nu=\nu_{\vec{m}}$ with the desired properties is to study smooth transversal sections $\bar{\nu}$ in the 
cokernel bundle $\overline{\Coker}\CR_J$ over the compactification $\CM$ of the moduli space $\IM=\IM_{\vec{m}}$. 
More precisely, we study smooth transversal sections in the cokernel bundle over the quotient $\CM/(\IZ_{m^+}\times\IZ_{m^-})$, i.e., we again forget the 
positions of the asymptotic markers. Then the 
abstract perturbation $\nu=\nu_{\vec{m}}$ for the moduli space $\IM$ is given by the restriction $\nu=\bar{\nu}|_{\IM}$ of the section to the 
interior, while the abstract perturbations $\nu_{\vec{m'}}$ 
for the moduli spaces for tuples $\vec{m'}$ with $\sharp\vec{m'}<\sharp\vec{m}$ and $|\vec{m'}|\leq |\vec{m}|$ determine $\bar{\nu}$ on the boundary 
$\del\CM=\CM-\IM$ as follows: \\
   
Let $(T,\LL)$ be a tree with level structure which represents curves with no nodes, i.e., for which all trees $T_{\ell,k}$, $\ell=1,...,L$, 
$k=1,...,N_{\ell}$ are trivial, and denote again $\IM_{\ell,k}=\IM_{T_{\ell,k}}$ the corresponding moduli spaces of branched covers. Let us further 
assume that $\IM_{T,\LL}$ is indeed a boundary stratum, i.e., does not agree with the top stratum $\IM$. Denoting by 
$\Coker^{\ell,k}\CR_J$ the cokernel bundle over $\IM_{\ell,k}/(\IZ_{m^+_{\ell,k}}\times\IZ_{m^-_{\ell,k}})$ with the sets 
$\IZ_{m^{\pm}_{\ell,k}}=\IZ_{m^{\pm}_{\ell,k,1}}\times ... \times\IZ_{m^{\pm}_{\ell,k,n^{\pm}_{\ell,k}}}$ of asymptotic markers at the 
positive, respectively negative punctures, recall that the cokernel bundle over $\IM_{T,\LL}/(\IZ_{m^+}\times\IZ_{m^-})$ is given as sum of pullback bundles
\begin{equation*}
  \Coker^{T,\LL}\CR_J = \pi_{1,1}^*\Coker^{1,1}\CR_J\oplus ... \oplus \pi_{L,N_L}^*\Coker^{L,N_L}\CR_J
\end{equation*}
under the projections 
\begin{equation*}
 \pi_{\ell,k}: \IM_{T,\LL}/(\IZ_{m^+}\times\IZ_{m^-}) \to \IM_{\ell,k}/(\IZ_{m^+_{\ell,k}}\times\IZ_{m^-_{\ell,k}}).
\end{equation*}
For the section $\bar{\nu}$ in the cokernel bundle over $\CM/(\IZ_{m^+}\times\IZ_{m^-})$ we now require that the restriction $\nu_{T,\LL}$ to 
$\IM_{T,\LL}/(\IZ_{m^+}\times\IZ_{m^-})$ is given by  
\begin{equation*}
 \nu_{T,\LL}(h,j)=(\nu_{1,1}(h_{1,1},j_{1,1}),...,\nu_{L,N_L}(h_{L,N_L},j_{L,N_L}))
\end{equation*}
for $(h,j)\in\IM_{T,\LL}$ with $h=(h_{1,1},...,h_{L,N_L})$, $j=(j_{1,1},...,j_{L,N_L})$. In other words, $\bar{\nu}$ is over $\IM_{T,\LL}$ given as 
sum of pullback sections 
\begin{equation*}
   \bar{\nu}|_{\IM_{T,\LL}} = \pi_{1,1}^*\nu_{1,1}\oplus ... \oplus \pi_{L,N_L}^*\nu_{L,N_L}
\end{equation*}
with $\nu_{1,1}$, ..., $\nu_{L,N_L}$ chosen before. 
Note that this makes sense, since all sections $\nu_{\ell,k}$ indeed live in the cokernel bundle $\Coker^{\ell,k}\CR_J$ over the quotient 
$\IM_{\ell,k}/(\IZ_{m^+_{\ell,k}}\times\IZ_{m^-_{\ell,k}})$. \\

We define $\CM^{\bar{\nu}}\subset\CM$ by pulling back the section $\bar{\nu}$ from the cokernel bundle over $\CM/(\IZ_{m^+}\times\IZ_{m^-})$ to 
the cokernel bundle over $\CM$ and setting 
\begin{equation*} \CM^{\bar{\nu}}=\bar{\nu}^{-1}(0). \end{equation*}
Recall that we have seen in proposition 2.5 that the cokernel bundle $\overline{\Coker}\CR_J$ can be equipped with the structure of a smooth vector bundle over 
the compactified moduli space $\CM$, which by proposition 1.4 is a smooth manifold with corners. Since $\bar{\nu}$ is assumed to be smooth and 
transversal to the zero section, it follows from a version of the implicit function theorem that $\CM^{\bar{\nu}}$ is a smooth submanifold with corners 
of $\CM$, which is furthermore neat in the sense that 
\begin{equation*} \del\CM^{\bar{\nu}}=\CM^{\bar{\nu}}\cap \del\CM. \end{equation*}
More precisely it follows that $\CM^{\bar{\nu}}$ is again a stratified space with strata 
\begin{equation*} \IM^{\bar{\nu}}_{T,\LL}=\CM^{\bar{\nu}}\cap\IM_{T,\LL}. \end{equation*}
Since $\bar{\nu}$ is independent of the directions of the asymptotic markers at the punctures, it follows that $\IZ_{m^+}\times\IZ_{m^-}$ 
still acts on the regular moduli space $\CM^{\bar{\nu}}$. Furthermore the conditions on the section $\bar{\nu}$ imply that for  
$\IM^{\bar{\nu}}_{T,\LL}$ with trivial trees $T_{1,1},...,T_{L,N_L}$ we have
\begin{equation*} 
     \IM^{\bar{\nu}}_{T,\LL}= \prod_{\ell=1}^L\IR^{N_{\ell}-1}\times \IM^{\nu_{1,1}}_{1,1}\times ...\times\IM^{\nu_{L,N_L}}_{L,N_L}/\Delta. 
\end{equation*}
This motivates the following definition: \\
\\
{\bf Definition 3.2:} {\it A section $\bar{\nu}$ in the cokernel bundle $\overline{\Coker}\CR_J$ over the compactified moduli space $\CM$ 
is called} coherent {\it if it is the pullback of a section in the cokernel bundle over the quotient $\CM/(\IZ_{m^+}\times\IZ_{m^-})$, 
which over each boundary stratum $\IM_{T,\LL}/(\IZ_{m^+}\times\IZ_{m^-})$ for trees with level structure $(T,\LL)$ with trivial trees $T_{1,1},...,T_{L,N_L}$ 
is given as sum} 
\begin{equation*}
       \bar{\nu}|_{\IM_{T,\LL}} = \pi_{1,1}^*\nu_{1,1}\oplus ... \oplus \pi_{L,N_L}^*\nu_{L,N_L}
\end{equation*}
{\it of pullbacks of sections in the cokernel bundles $\Coker^{T_{1,1}}\CR_J$,...,$\Coker^{T_{L,N_L}}\CR_J$ under the projections} 
\begin{equation*}
 \pi_{\ell,k}: \IM_{T,\LL}/(\IZ_{m^+}\times\IZ_{m^-}) \to \IM_{T_{\ell,k}}/(\IZ_{m^+_{\ell,k}}\times\IZ_{m^-_{\ell,k}}).
\end{equation*}

We emphasize that our notion of coherency is weaker than the usual definition: While we just require that the abstract perturbations are of a special form over 
each boundary stratum, one usually additionally requires that for every moduli space one chooses a unique abstract perturbation in the sense that if a moduli space 
appears in two different boundary strata the two perturbations agree. However it follows from our proof of theorem 3.3 that our weaker assumption indeed suffices 
to prove our desired result. \\
  
Let $\IM_1\times_{\IZ_{m_{12}}}\IM_2\subset\del\CM$ be an arbitrary codimension one boundary stratum. Recall from subsection 2.1 that the restriction 
of $\overline{\Coker}\CR_J$ to $\CM_1\times_{\IZ_{m_{12}}}\CM_2 = \overline{\IM_1\times_{\IZ_{m_{12}}}\IM_2}\subset\CM$ is 
given by the sum of pullback bundles  
\begin{equation*} 
 \overline{\Coker}\CR_J|_{\CM_1\times_{\IZ_{m_{1,2}}}\CM_2} \;=\; 
 \pi_1^*\overline{\Coker}^1\CR_J \oplus \pi_2^*\overline{\Coker}^2\CR_J. 
\end{equation*}
It directly follows from the definition that any coherent (and transversal) section $\bar{\nu}$ in $\overline{\Coker}\CR_J$ is given over 
$\CM_1\times_{\IZ_{m_{12}}}\CM_2$ by 
\begin{equation*} 
     \bar{\nu}|_{\CM_1\times_{\IZ_{m_{12}}}\CM_2}=\pi_1^*\bar{\nu}_1\oplus\pi_2^*\bar{\nu}_2 
\end{equation*}
with coherent (and transversal) sections in $\overline{\Coker}^1\CR_J$, $\overline{\Coker}^2\CR_J$, respectively. Furthermore 
\begin{equation*}
    \overline{\IM_1\times_{\IZ_{m_{12}}}\IM_2}^{\bar{\nu}}= \CM_1^{\bar{\nu}_1}\times_{\IZ_{m_{12}}}\CM_2^{\bar{\nu}_2}.
\end{equation*} 
  
\subsection{Euler numbers for Fredholm problems}

We have seen that the perturbation chosen for a moduli space explicitly depends on the perturbations chosen 
for the moduli spaces forming the boundary of its compactification. However, in this last section we prove that for 
any coherent and transversal section $\bar{\nu}$ in $\overline{\Coker}\CR_J$ in the sense of definition 3.2 the algebraic count of elements 
in the regular compactified moduli space $\CM^{\bar{\nu}}$ is zero, independent of all choices. Together with the discussion in 0.3 
we have then shown that branched covers over trivial orbit cylinders do not contribute to the differential of rational 
symplectic field theory, i.e., we have proven the main theorem. \\
\\
{\bf Theorem 3.3:} {\it For the cokernel bundle $\overline{\Coker}\CR_J$ over the compactification $\CM$ of every moduli space 
of branched covers over a trivial cylinder with $\dim_{\virt}\IM=\dim\IM-\rank\Coker\CR_J=0$ the following holds:}
\begin{itemize} 
\item {\it For every pair $\bar{\nu}^0$, $\bar{\nu}^1$ of coherent and transversal sections in $\overline{\Coker}\CR_J$ the algebraic count of zeroes 
      of $\bar{\nu}^0$ and $\bar{\nu}^1$ are finite and agree, so that we can define an Euler number 
      $\chi(\overline{\Coker}\CR_J)$ for coherent sections in $\overline{\Coker}\CR_J$ by} 
      \begin{equation*} \chi(\overline{\Coker}\CR_J) \,:=\, \sharp (\bar{\nu}^0)^{-1}(0) \,=\, \sharp (\bar{\nu}^1)^{-1}(0). 
      \end{equation*} 
\item {\it This Euler number is $\chi(\overline{\Coker}\CR_J) = 0$.} 
\end{itemize}
$ $ \\
{\it Proof:} We prove this statement for all moduli spaces of orbit curves by induction on the number of punctures $n\geq 3$. \\

Let $\bar{\nu}$ be a coherent and transversal section in $\overline{\Coker}\CR_J$. In order to see that the zeroes of $\bar{\nu}$ 
can be counted to give a finite number, observe that it follows from $\dim\IM-\rank\overline{\Coker}\CR_J=0$ and the implicit function theorem that 
$\bar{\nu}^{-1}(0)$ is a neat zero-dimensional submanifold of $\CM$, i.e., a discrete set of points in $\IM\subset\CM$, which is compact as 
a closed subset of a compact set. \\

Now let $\bar{\nu}^0$ and $\bar{\nu}^1$ be two coherent and transversal sections in $\overline{\Coker}\CR_J$. In order to see that the numbers of zeroes 
$\sharp (\bar{\nu}^0)^{-1}(0)$ and $\sharp (\bar{\nu}^1)^{-1}(0)$ indeed agree, let $\bar{\nu}^{01}$ be a section in the cokernel bundle 
$\overline{\Coker_0}\CR_J$ over $\overline{\IM^0}/(\IZ_{m^+}\times\IZ_{m^-})$, which is coherent and compatible with $\bar{\nu}^0$ and $\bar{\nu}^1$ in the 
sense that over each stratum 
\begin{eqnarray*} 
&& \IM^0_{T,\LL,\ell_0}/(\IZ_{m^+}\times\IZ_{m^-}) \\
&& = \frac{\IM_{T_1}\times ... \times \IM_{T_{\ell_0-1}}\times \IM^0_{T_{\ell_0}}\times \IM_{T_{\ell_0+1}}\times ... \times \IM_{T_L}}
{\Delta\times \IZ_{m^+}\times\IZ_{m^-}}
\end{eqnarray*}
the restriction $\nu^{01}_{T,\LL,\ell_0} = \bar{\nu}^{01}|_{\IM^0_{T,\LL,\ell_0}}$ is of the form  
\begin{eqnarray*}
  \nu^{01}_{T,\LL,\ell_0} &=& \pi_1^*\nu^0_{T_1}\oplus ... 
       \oplus \pi_{\ell_0-1}^*\nu^0_{T_{\ell_0-1}}\oplus \pi_{\ell_0}^*\nu^{01}_{T_{\ell_0}} \\
       && \oplus \pi_{\ell_0+1}^*\nu^1_{T_{\ell_0+1}}\oplus ... \oplus \pi_L^*\nu^1_{T_L},
\end{eqnarray*}
with the projections 
\begin{equation*}
 \pi_{\ell,k}: \IM^0_{T,\LL,\ell_0}/(\IZ_{m^+}\times\IZ_{m^-}) \to \IM^{(0)}_{T_{\ell,k}}/(\IZ_{m^+_{\ell,k}}\times\IZ_{m^-_{\ell,k}}),
\end{equation*}
where $\nu^0_{T_1},...,\nu^0_{T_{\ell_0-1}}$ and $\nu^1_{T_{\ell_0+1}},...,\nu^1_{T_L}$ are given by $\bar{\nu}^0$ and $\bar{\nu}^1$, respectively.
Note that this implies that  
\begin{eqnarray*}
   \bar{\nu}^{01}|_{\overline{\IM^0_1}\times_{\IZ_{m_{12}}}\CM_2} &=& \pi_1^*\bar{\nu}^{01}_1\oplus\pi_2^*\bar{\nu}^1_2, \\
   \bar{\nu}^{01}|_{\CM_1\times_{\IZ_{m_{12}}}\overline{\IM^0_2}} &=& \pi_1^*\bar{\nu}^0_1\oplus\pi_2^*\bar{\nu}^{01}_2, \\
   \bar{\nu}^{01}|_{\{\point\}\times\CM} = \bar{\nu}^1 &\textrm{and}& \bar{\nu}^{01}|_{\CM\times\{\point\}} = \bar{\nu}^0
\end{eqnarray*}
and that we can always find a section $\bar{\nu}^{01}$ in $\overline{\Coker_0}\CR_J$ with the above properties by iteratively extending as in 3.2 the sections 
from the boundary of $\IM^0$ to the interior of the moduli space. In particular, observe that by proposition 1.2 the numbers of punctures of the curves in 
$\IM^0_1$ and $\IM^0_2$ in the codimension one boundary of $\IM^0$ are strictly smaller than the number of punctures of the curves in $\IM^0$. \\
 
Note that by the propositions 1.4 and 2.5 the cokernel bundle $\overline{\Coker_0}\CR_J$ over $\overline{\IM^0}$ can also be equipped with the structure of a 
smooth vector bundle over a manifold with corners. With this we further again assume that $\bar{\nu}^{01}$ is a smooth and transversal section in 
$\overline{\Coker_0}\CR_J$, which in turn implies that 
for each stratum $\IM^0_{T,\LL,\ell_0}$ the underlying sections $\nu^0_{T_1}$, ..., $\nu^0_{T_{\ell_0-1}}$, $\nu^{01}_{T_{\ell_0}}$, $\nu^1_{T_{\ell_0+1}}$, 
..., $\nu^1_{T_{\ell_0}}$ of the cokernel bundles $\Coker^{T_1}\CR_J$, ..., $\Coker^{T_{\ell_0}}_0\CR_J$, ..., $\Coker^{T_L}\CR_J$ are again smooth and transversal. 
Now it follows from 
\begin{equation*} 
 \dim\IM^0-\rank\overline{\Coker_0}\CR_J = 1+\dim\IM-\rank\overline{\Coker}\CR_J = 1
\end{equation*}
that the resulting regular moduli space 
\begin{equation*} \overline{\IM^0}^{\bar{\nu}^{01}} = (\bar{\nu}^{01})^{-1}(0) \subset \overline{\IM^0} \end{equation*}
is a neat one-dimensional submanifold of $\overline{\IM^0}$. In other words, we have that $\overline{\IM^0}^{\bar{\nu}^{01}}$ is a one-dimensional manifold with boundary 
given by 
\begin{equation*}
   \del\overline{\IM^0}^{\bar{\nu}^{01}} = \overline{\IM^0}^{\bar{\nu}^{01}} \cap \del\overline{\IM^0}. 
\end{equation*}
 
In order to determine the boundary of $\overline{\IM^0}^{\bar{\nu}^{01}}$ observe that after setting 
\begin{eqnarray*}
   \bigl(\overline{\IM^0_1}\times_{\IZ_{m_{12}}}\CM_2\bigr)^{\bar{\nu}^{01}} &:=& 
   \overline{\IM^0}^{\bar{\nu}^{01}} \cap \bigl(\overline{\IM^0_1}\times_{\IZ_{m_{12}}}\CM_2\bigr), \\
   \bigl(\CM_1\times_{\IZ_{m_{12}}}\overline{\IM^0_2}\bigr)^{\bar{\nu}^{01}} &:=& 
   \overline{\IM^0}^{\bar{\nu}^{01}} \cap \bigl(\CM_1\times_{\IZ_{m_{12}}}\overline{\IM^0_2}\bigr), \\
   \bigl(\{\point\}\times \CM\bigr)^{\bar{\nu}^{01}} &:=& \overline{\IM^0}^{\bar{\nu}^{01}} \cap \bigl(\{\point\}\times\CM\bigr), \\
   \textrm{and}\;\;\;\; \bigl(\CM\times\{\point\}\bigr)^{\bar{\nu}^{01}} &:=& \overline{\IM^0}^{\bar{\nu}^{01}} \cap \bigl(\CM\times\{\point\}\bigr),
\end{eqnarray*}
the boundary conditions for $\bar{\nu}^{01}$ yield 
\begin{eqnarray*}
   \bigl(\overline{\IM^0_1}\times_{\IZ_{m_{12}}}\CM_2\bigr)^{\bar{\nu}^{01}}
   &=& \overline{\IM^0_1}^{\bar{\nu}^{01}_1}\times_{\IZ_{m_{12}}}\CM_2^{\bar{\nu}^1_2}, \\
   \bigl(\CM_1\times_{\IZ_{m_{12}}}\overline{\IM^0_2}\bigr)^{\bar{\nu}^{01}} 
   &=& \CM_1^{\bar{\nu}^0_1}\times_{\IZ_{m_{12}}}\overline{\IM^0_2}^{\bar{\nu}^{01}_2}, \\
   \bigl(\{\point\}\times\CM\bigr)^{\bar{\nu}^{01}} &=& \{\point\}\times\CM^{\bar{\nu}^1}, \\
   \textrm{and}\;\;\;\; \bigl(\CM\times\{\point\}\bigr)^{\bar{\nu}^{01}} &=& \CM^{\bar{\nu}^0}\times\{\point\}. 
\end{eqnarray*}
All together it follows that the boundary of $\overline{\IM^0}^{\bar{\nu}^{01}}$ is given by 
\begin{eqnarray*}
 \del\overline{\IM^0}^{\bar{\nu}^{01}} &=& (\CM^{\bar{\nu}^0}\times\{\point\}) \;\cup\; (\{\point\}\times\CM^{\bar{\nu}^1}) \\
 &\cup& \bigcup_{2<n_1,n_2<n} \Bigl((\CM_1^{\bar{\nu}^0_1}\times_{\IZ_{m_{12}}}\overline{\IM^0_2}^{\bar{\nu}^{01}_2})
 \;\cup\; (\overline{\IM^0_1}^{\bar{\nu}^{01}_1}\times_{\IZ_{m_{12}}}\CM_2^{\bar{\nu}^1_2})\Bigr),
\end{eqnarray*}   
where we take the union over all those codimension one boundary strata of $\overline{\IM^0}$ where the number of punctures 
$n_1,n_2$ for $\overline{\IM^0_1}$, $\CM_2$ (and $\CM_1$, $\overline{\IM^0_2}$) is strictly between two and the number of punctures $n$ for 
$\overline{\IM_0}$, i.e., $\overline{\IM^0_1}$, $\overline{\IM^0_2}\neq\{\point\}$. \\

Now since $\del\overline{\IM^0}^{\bar{\nu}^{01}}$ is 
the boundary of a one-dimensional manifold and taking into account the orientation of the codimension one boundary of the base space 
$\overline{\IM^0}$ it follows that 
\begin{eqnarray*}
 0 &=& \#(\CM^{\bar{\nu}^0}\times\{\point\}) \;-\; \#(\{\point\}\times\CM^{\bar{\nu}^1}) \\
 &+& \sum_{2<n_1,n_2<n} \Bigl(\#(\CM_1^{\bar{\nu}^0_1}\times_{\IZ_{m_{12}}}\overline{\IM^0_2}^{\bar{\nu}^{01}_2})
 \;-\; \#(\overline{\IM^0_1}^{\bar{\nu}^{01}_1}\times_{\IZ_{m_{12}}}\CM_2^{\bar{\nu}^1_2})\Bigr).
\end{eqnarray*}   
Note that here $\#$ refers to the orientation as boundary of $(\IM^0)^{\bar{\nu}^{01}}$, which itself is induced by the orientation 
of the cokernel bundle $\Coker_0\CR_J$ over $\IM^0$.
In order to show that 
\begin{equation*}
\#(\bar{\nu}^0)^{-1}(0)=\#\CM^{\bar{\nu}^0}=\#\CM^{\bar{\nu}^1}=\#(\bar{\nu}^1)^{-1}(0), 
\end{equation*}
i.e., to prove the existence of the Euler number $\chi(\overline{\Coker}\CR_J)$, it hence suffices to show that 
\begin{eqnarray*}
             &&\#(\CM_1^{\bar{\nu}^0_1}\times_{\IZ_{m_{12}}}\overline{\IM^0_2}^{\bar{\nu}^{01}_2})= 0 \\
\textrm{and} &&\#(\overline{\IM^0_1}^{\bar{\nu}^{01}_1}\times_{\IZ_{m_{12}}}\CM_2^{\bar{\nu}^1_2})= 0
\end{eqnarray*}
for every other boundary stratum: \\

For this observe that in order to have 
\begin{equation*}\#(\CM_1^{\bar{\nu}^0_1}\times_{\IZ_{m_{12}}}\overline{\IM^0_2}^{\bar{\nu}^{01}_2})\neq 0, \end{equation*} 
we in particular must have 
\begin{equation*} \CM_1^{\bar{\nu}^0_1}\times_{\IZ_{m_{12}}}\overline{\IM^0_2}^{\bar{\nu}^{01}_2}\neq\emptyset, \end{equation*} 
which is equivalent to 
\begin{equation*} \CM_1^{\bar{\nu}^0_1}\neq\emptyset \;\;\;\;\;\;\textrm{and}\;\;\;\;\;\; \overline{\IM^0_2}^{\bar{\nu}^{01}_2}\neq\emptyset. \end{equation*}
Now since both $\bar{\nu}^0_1$ and $\bar{\nu}^{01}_2$ are transversal, i.e., have zero as a regular value, it follows that 
\begin{eqnarray*} 
 \dim\IM_1 - \rank\overline{\Coker}^1\CR_J = \dim \IM_1^{\nu^0_1} &\geq& 0 \\
 \textrm{and}\;\;\;\;\;\; \dim\IM^0_2 - \rank\overline{\Coker_0}^2\CR_J = \dim (\IM^0_2)^{\nu^{01}_2} &\geq& 0.
\end{eqnarray*}
On the other hand, since 
\begin{eqnarray*} 
 1 &=& \dim\IM^0 - \rank\overline{\Coker_0}\CR_J \\
   &=& 1 + \dim \IM_1 + \dim \IM^0_2 - \rank\overline{\Coker}^1\CR_J - \rank\overline{\Coker_0}^2 \CR_J
\end{eqnarray*}
it follows that we indeed must have equality, i.e., 
\begin{eqnarray*} 
 \dim\IM_1 - \rank\overline{\Coker}^1\CR_J = \dim \IM_1^{\nu^0_1} &=& 0 \\
 \textrm{and}\;\;\;\;\;\; \dim\IM^0_2 - \rank\overline{\Coker_0}^2\CR_J = \dim (\IM^0_2)^{\nu^{01}_2} &=& 0.
\end{eqnarray*}

In other words, we can immediately forget about all boundary components $\CM_1^{\bar{\nu}^0_1}\times_{\IZ_{m_{12}}}\overline{\IM^0_2}^{\bar{\nu}^{01}_2}$ 
where the virtual dimension of $\IM_1$ (and of $\IM^0_2$) is not equal to zero, i.e., the underlying Fredholm index of 
$\CR_J$ is not equal to one, where a corresponding statement clearly also holds for the boundary strata 
$\overline{\IM^0_1}^{\bar{\nu}^{01}_1}\times_{\IZ_{m_{12}}}\CM_2^{\bar{\nu}^1_2}$. In particular, observe that this would directly prove the existence 
of the desired Euler number $\chi(\overline{\Coker}\CR_J)$ if we were able to show that none of the moduli spaces $\IM_1$ or $\IM_2$ appearing 
in the codimension one boundary has virtual dimension zero. While this is typically the case when the compactification is not ``too large'' , 
note that here there is no way to exclude the latter from happening. However, at this point, we can now make use of the induction hypothesis as follows: \\

Since the number of punctures for the moduli space $\CM_1$ and $\CM_2$ is strictly smaller than the number of punctures for the original 
moduli space $\CM$, it follows that we do not only have Euler numbers $\chi(\overline{\Coker}^1\CR_J)$ and $\chi(\overline{\Coker}^2\CR_J)$ 
for coherent and transversal sections in the cokernel bundles $\overline{\Coker}^1\CR_J$ and $\overline{\Coker}^2\CR_J$, but by assumption 
further know that they are zero. In other words, we already know that 
\begin{equation*} \#_1\CM_1^{\bar{\nu}^0_1} \;=\;\chi(\overline{\Coker}^1\CR_J) \;=\; 0,\;\;\;\; 
                  \#_2\CM_2^{\bar{\nu}^1_2} \;=\;\chi(\overline{\Coker}^2\CR_J) \;=\; 0, \end{equation*}
where $\#_1$, $\#_2$ refers to the algebraic count with respect to the orientation on the cokernel bundle $\Coker^1\CR_J$, $\Coker^2\CR_J$ over $\IM_1$, 
$\IM_2$, respectively. Denoting by $\#_1$, $\#_2$ further the algebraic count with respect to the induced orientation on $\Coker_0^1\CR_J$, $\Coker_0^2\CR_J$ 
it follows that 
\begin{eqnarray*} 
   &&\#_{12} (\CM_1^{\bar{\nu}^0_1}\times_{\IZ_{m_{12}}}\overline{\IM^0_2}^{\bar{\nu}^{01}_2}) \;=\; 
     \frac{1}{m_{12}} \cdot \#_1 \CM_1^{\bar{\nu}^0_1}\cdot \#_2\overline{\IM^0_2}^{\bar{\nu}^{01}_2} \;=\; 0, \\  
   &&\#_{12} (\overline{\IM^0_1}^{\bar{\nu}^{01}_1}\times_{\IZ_{m_{12}}}\CM_2^{\bar{\nu}^1_2}) \;=\; 
     \frac{1}{m_{12}} \cdot \#_1\overline{\IM^0_1}^{\bar{\nu}^{01}_1} \cdot \#_2 \CM_2^{\bar{\nu}^1_2} \;=\; 0,  
\end{eqnarray*} 
where $\#_{12}$ refers to the induced orientations on $\pi_1^*\Coker^1\CR_J\oplus\pi_2^*\Coker_0^2\CR_J$. But since the algebraic counts  
$\#$ and $\#_{12}$ differ only by sign by proposition 2.7, it follows that   
\begin{equation*}  \# (\CM_1^{\bar{\nu}^0_1}\times_{\IZ_{m_{12}}}\overline{\IM^0_2}^{\bar{\nu}^{01}_2}) = 0,\;\;\;
                   \# (\overline{\IM^0_1}^{\bar{\nu}^{01}_1}\times_{\IZ_{m_{12}}}\CM_2^{\bar{\nu}^1_2}) = 0, \end{equation*}
which proves the first part of the theorem. \\ 
  
It remains to prove $\chi(\overline{\Coker}\CR_J)=0$: \\

But for this we must only observe that  the rank of $\overline{\Coker}\CR_J$ is always 
odd, since it agrees with the dimension of $\IM$, which itself is the product 
of a one-dimensional manifold with a complex manifold. Indeed, we have 
\begin{equation*} \rank \overline{\Coker}\CR_J = \dim \IM \,=\, \dim(S^1\times\IM_{0,n}) = 2(n-3)+1 \equiv 1 \mod 2.
\end{equation*} 

Following the idea of proving the vanishing of the Euler characteristic for odd-dimensional closed manifolds, 
observe that for any coherent and transversal section $\bar{\nu}$ in $\overline{\Coker}\CR_J$ the section $-\bar{\nu}$ has 
the same property and we have 
\begin{equation*} 
 \chi(\overline{\Coker}\CR_J) \,=\, \sharp (-\bar{\nu})^{-1}(0) 
 \,=\,\, -\,\sharp \bar{\nu}^{-1}(0) = - \chi(\overline{\Coker}\CR_J),
\end{equation*}
\\
implying $\chi(\overline{\Coker}\CR_J) = 0$. $\qed$

\section{Consequences}
\subsection{Action filtration on rational symplectic field theory} 

In this section we want to discuss the implications of our main theorem on rational symplectic field theory. While we have seen that the problem 
of achieving regularity for moduli spaces already appears in the case of orbit curves, which we however have settled above using obstruction bundles, 
note that our method does not allow us to solve the problem for the other moduli spaces studied in rational symplectic field theory. Beside the fact 
that we cannot assume the nonregular moduli spaces to be manifolds in general, we further cannot assume that the cokernels of the linearizations 
of the Cauchy-Riemann operator fit together to give a vector bundle of the right rank over the nonregular moduli space. For this recall that we have proven 
the latter by a linearized energy argument in proposition 2.3 which is not available in the general case. In order to settle the transversality problem in 
symplectic field theory H. Hofer, K. Wysocki and E. Zehnder invented the theory of polyfolds, which however at the moment of writing this paper is 
still on its way of being completed. While our result about orbit curves in rational symplectic field theory is itself independent of the methods 
used to achieve regularity in the general case, let us outline how our result can be embedded in the general story: \\

While the most natural way consists in using our obstruction bundle perturbations for the moduli spaces of orbit curves and extending them via the 
polyfold theory to abstract perturbations for all other moduli spaces, we claim that the statement of the main theorem is true independent of the 
method used to define the coherent compact perturbations. In particular it should hold for the 
abstract perturbations constructed using the polyfold theory of [HWZ] as well as the domain-dependent Hamiltonian perturbations used in [F1]. Since 
the analytical foundations of symplectic field theory are not yet established, we cannot make the above statement rigorous in full detail. However, 
let us point out the important consequences of our result to symplectic field theory of which we are confident that they can be shown once the analytical 
tools from polyfold theory are available. Despite the fact that we can not make them rigorous by the aforementioned reasons, we decided to state them as 
propositions with proofs as it is common in recent papers on symplectic field theory, see e.g. [B] and [EGH]. \\
\\
{\bf Proposition 4.1:} {\it For all choices of coherent compact perturbations $\nu$ which make the perturbed Cauchy-Riemann operator $\CR_J^{\nu}=\CR_J+\nu$ 
transversal to the zero section in an appropriate Banach space bundle (or polyfold) setup, the algebraic count of elements in the resulting regular 
moduli space $\CM^{\nu}=(\CR_J^{\nu})^{-1}(0)$ is zero. 
It follows that branched covers over orbit cylinders do not contribute to the algebraic invariants of rational symplectic field theory.} \\
\\
{\it Proof:} Here we proceed as in the proof of theorem 3.3 and prove the statement by induction on the number of punctures. 
For every moduli space of orbit curves $\IM$ assume we are given an arbitrary coherent perturbations 
$\bar{\nu}^0$ and $\bar{\nu}^1$, constructed e.g. using the polyfold theory of [HWZ], which, after being added to $\CR_J$, make all strata of the 
compactification $\CM$ regular. Using polyfold theory we can construct a compact perturbation $\bar{\nu}^{01}$ of $\overline{\IM^0}$ so that, 
in the notation from before, the codimension one boundary strata of the resulting regular moduli space $\overline{\IM^0}^{\bar{\nu}^{01}}$ are again given by
\begin{eqnarray*}
   \bigl(\overline{\IM^0_1}\times_{\IZ_{m_{12}}}\CM_2\bigr)^{\bar{\nu}^{01}}
   &=& \overline{\IM^0_1}^{\bar{\nu}^{01}_1}\times_{\IZ_{m_{12}}}\CM_2^{\bar{\nu}^1_2}, \\
   \bigl(\CM_1\times_{\IZ_{m_{12}}}\overline{\IM^0_2}\bigr)^{\bar{\nu}^{01}} 
   &=& \CM_1^{\bar{\nu}^0_1}\times_{\IZ_{m_{12}}}\overline{\IM^0_2}^{\bar{\nu}^{01}_2}, \\
   \bigl(\{\point\}\times\CM\bigr)^{\bar{\nu}^{01}} &=& \{\point\}\times\CM^{\bar{\nu}^1}, \\
   \textrm{and}\;\;\;\; \bigl(\CM\times\{\point\}\bigr)^{\bar{\nu}^{01}} &=& \CM^{\bar{\nu}^0}\times\{\point\}. 
\end{eqnarray*}
In particular we again have  
\begin{eqnarray*}
 && \#\CM^{\bar{\nu}^0}\;-\; \#\CM^{\bar{\nu}^1} \\
 && = \;\sum_{2<n_1,n_2<n} \Bigl(\#(\overline{\IM^0_1}^{\bar{\nu}^{01}_1}\times_{\IZ_{m_{12}}}\CM_2^{\bar{\nu}^1_2})
 \;-\; \#(\CM_1^{\bar{\nu}^0_1}\times_{\IZ_{m_{12}}}\overline{\IM^0_2}^{\bar{\nu}^{01}_2})\Bigr).
\end{eqnarray*}   
Using the induction hypothesis it follows as before that the right hand side of the equation is equal to zero, so that 
\begin{equation*}
 \#\CM^{\bar{\nu}^0} \;\;=\;\; \#\CM^{\bar{\nu}^1},
\end{equation*}
i.e., the number of elements in the regular moduli space is independent of {\it any} choice of coherent compact perturbations. 
Assuming in particular that $\bar{\nu}^0$ is a coherent compact perturbation resulting from a section in the cokernel bundle $\overline{\Coker}\CR_J$ as studied 
before, it follows that this number is zero. $\qed$ \\

Like in Gromov-Witten theory and symplectic Floer homology the orbit curves in symplectic field theory can be characterized by the fact 
that they carry no energy in a certain sense, which, as in Floer homology, can be expressed as difference of actions assigned to the asymptotic 
periodic orbits. More precisely, we can introduce a  natural action filtration on rational symplectic field theory as follows: \\

The action 
\begin{equation*} S(\gamma) = \int f_{\gamma}^*\omega, \end{equation*} 
which we defined in 0.2 using the spanning surface $f_{\gamma}$ for every closed Reeb orbit $\gamma$, naturally defines an action filtration 
$\FF$ on the chain algebras $\IA$ and $\IP$ underlying contact homology 
and rational symplectic field theory. For this observe that over the group ring over $H_2(V)$ $\IA$ and $\IP$ are generated by the formal variables 
$q_{\gamma}$ (and $p_{\gamma}$) assigned to every good orbit $\gamma$ in the sense of [BM], so that for every monomial we can define
\begin{equation*}
 \FF(q_{\gamma^-_1}...q_{\gamma^-_{n^-}} p_{\gamma^+_1}...p_{\gamma^+_{n^+}}e^A) := 
 \sum_{k=1}^{n^-} S(\gamma^-_k) - \sum_{\ell=1}^{n^+} S(\gamma^+_{\ell}) + \omega(A).
\end{equation*}
Note that in the contact case, i.e., where the one-form $\lambda$ of the Hamiltonian structure on $V$ is contact 
and $\omega=d\lambda$, we have $\omega(A)=0$ and the action for the periodic orbits $\gamma$, i.e., the closed Reeb orbits, 
is given by integrating the one-form $\lambda$ along $\gamma$. \\ 
\\
{\bf Corollary 4.2:} {\it Like in cylindrical contact homology the differential in contact homology and rational symplectic 
field theory is strictly decreasing with respect to the action filtration.} \\
\\ 
{\it Proof:} Since the differential $d=d^{\hh}=\{\hh,\cdot\}:\IP\to\IP$ in rational symplectic field theory, given by the generating function $\hh\in\IP$ counting 
holomorphic curves in $\IR\times V$, satisfies a graded Leibniz rule, it is strictly decreasing with respect 
to $\FF$ precisely when for every orbit $\gamma$, 
\begin{eqnarray*} 
                     \<d p_{\gamma},p^{\Gamma^+}q^{\Gamma^-}e^A\>\neq 0 &\textrm{implies}& \FF(p_{\gamma})>\FF(p^{\Gamma^+}q^{\Gamma^-}e^A) \\
\textrm{and} \;\;\;\;\<d q_{\gamma},p^{\Gamma^+}q^{\Gamma^-}e^A\>\neq 0 &\textrm{implies}& \FF(q_{\gamma})>\FF(p^{\Gamma^+}q^{\Gamma^-}e^A), 
\end{eqnarray*}
where $\<d p_{\gamma},p^{\Gamma^+}q^{\Gamma^-}e^A\>$ and $\<d q_{\gamma},p^{\Gamma^+}q^{\Gamma^-}e^A\>$
denote the coefficients of 
\begin{equation*}p^{\Gamma^+}q^{\Gamma^-}e^A=p_{\gamma_1^+}...p_{\gamma_{n^+}^+}q_{\gamma_1^-}...q_{\gamma_{n^-}^-}e^A \end{equation*} 
in the series expansion of $d p_{\gamma}$ and $d q_{\gamma}$, respectively. On the other hand it follows from the definition of $d$ that  
\begin{eqnarray*} 
   \<d p_{\gamma},p^{\Gamma^+}q^{\Gamma^-}e^A\> &=&\<\{\Ih,p_{\gamma}\},p^{\Gamma^+}q^{\Gamma^-}e^A\>\\
   &=&\kappa_{\gamma}\;\<\frac{\del\Ih}{\del q_{\gamma}},p^{\Gamma^+}q^{\Gamma^-}e^A\>\\&=&\pm\kappa_{\gamma}\;\<\Ih,p^{\Gamma^+}(q^{\Gamma^-}q_{\gamma})e^A\> 
\end{eqnarray*}
with the Hamiltonian $\Ih\in\IP$ of rational symplectic field theory, and similar for $d q_{\gamma}$, so that the requirement on $d$ is 
equivalent to requiring that 
\begin{equation*}\<\Ih,p^{\Gamma^+}q^{\Gamma^-}e^A\>\neq 0\;\;\;\;\textrm{implies}\;\;\;\;\FF(p^{\Gamma^+}q^{\Gamma^-}e^A)>0. 
\end{equation*}
Note that here we use $\FF(q_{\gamma})=-\FF(p_{\gamma})$. In order to see how this follows 
from the above proposition, recall that $\<\Ih,p^{\Gamma^+}q^{\Gamma^-}e^A\>$ is given by the algebraic count of elements in the moduli space described 
by the monomial $p^{\Gamma^+}q^{\Gamma^-}e^A$, which consists of the curves which are 
asymptotically cylindrical over the orbits $\gamma_1^{\pm},...,\gamma_{n^{\pm}}^{\pm}$ at the positive, respectively negative punctures and represent 
the homology class $A\in H_2(V)$. On the other hand recall from 0.2 that the $\omega$-energy of a holomorphic curve $u$ in the moduli space can be expressed 
in terms of the actions of the closed orbits $\gamma_1^{\pm},...,\gamma_{n^{\pm}}^{\pm}$ and the integral of $\omega$ over the homology class $A\in H_2(V)$ by 
\begin{equation*} 
 E_{\omega}(u) = \sum_{k=1}^{n^+} S(\gamma_k^+) - \sum_{\ell=1}^{n^-} S(\gamma_{\ell}^-) + \omega(A),   
\end{equation*}
i.e., $E_{\omega}(u)= \FF(q_{\gamma^-_1}...q_{\gamma^-_{n^-}} p_{\gamma^+_1}...p_{\gamma^+_{n^+}}e^A)$. But since the algebraic count of curves in moduli spaces 
of curves with $E_{\omega}(u)=0$ is zero by proposition 4.1, we get the desired result. $\qed$ \\

Recall that this statement is trivial in the case of cylindrical contact homology and symplectic Floer homology since the only orbit curves 
in these cases are trivial cylinders. 
   
\subsection{Marked points, differential forms and the spectral sequence for filtered complexes}

Since orbit curves are characterized by the fact that they have trivial $\omega$-energy and this quantity is preserved under taking boundaries and 
gluing of moduli spaces, it follows that every algebraic invariant of rational symplectic field theory has a natural analog defined by counting 
only those orbit curves. More precisely, observe that the generating function $\hh\in\IP$ counting holomorphic curves in $(\IR\times V,J)$ can be 
written as a sum $\hh=\hh_0+\hh_{>0}$ where $\hh_0\in\IP$ is the generating function for the curves with trivial $\omega$-energy and $\hh_{>0}$ the 
one for the curves with strictly positive $\omega$-energy, which in turn immediately implies that also the differential $d=d^{\hh}:\IP\to\IP$ is given as a sum 
$d=d_0+d_{>0}$ with $d_0=d^{\hh_0}$, $d_{>0}=d^{\hh_{>0}}$. \\

In the same way as we use the study of the boundaries of one-dimensional moduli spaces 
(after quotiening out the $\IR$-action) to deduce the fundamental identity $\{\hh,\hh\}=0$ implying $d^2=0$, it follows from the aforementioned 
fact that the $\omega$-energy is preserved under taking boundaries and gluing of moduli spaces that we already have $\{\hh_0,\hh_0\}=0$ and 
therefore $d_0^2=0$. Even further it is clear that we already have $\{\hh_{0,\gamma},\hh_{0,\gamma}\}=0$ where $\hh_{0,\gamma}$ is 
the generating function counting all orbit curves over the closed Reeb orbit $\gamma\in P(V)$ so that $\hh_0=\sum_{\gamma\in P(V)} \hh_{0,\gamma}$. \\
 
Denoting by $\IP_{\gamma}$ the graded Poisson subalgebra of $\IP$ generated only by the variables $p_{\gamma^k}$, $q_{\gamma^k}$ assigned to multiple 
covers of the chosen Reeb orbit $\gamma$, observe that we have $\hh_{0,\gamma}\in\IP_{\gamma}$ so that $d_{0,\gamma}=\{\hh_{0,\gamma},\cdot\}$ defines a 
differential on $\IP_{\gamma}$. We call its homology $H_*(\IP_{\gamma},d_{0,\gamma})$ the {\it rational symplectic field theory of $\gamma$}. \\

While it follows from our main theorem that $h_{0,\gamma}=0$ and therefore $H_*(\IP_{\gamma},d_{0,\gamma})=\IP_{\gamma}$ when no differential forms 
are chosen, let us spend the remaining time studying what can be said about the general case described in [EGH] when a string of closed differential forms is introduced: \\

To this end, let $\Theta=(\theta_1,...,\theta_N)\in(\Omega^*(V))^N$ be a string of closed differential forms. Abbreviating $\gamma^{\vec{m}^{\pm}}=
(\gamma^{m_1^{\pm}},...,\gamma^{m_{n^{\pm}}^{\pm}})$, note that on every moduli space $\IM_{0,0,r}(\gamma^{\vec{m}^+},\gamma^{\vec{m}^-})$ of orbit curves 
with additional $r$ marked points $\uw=(w_1,...,w_r)\in\Si^r$ we have $r$ evaluation maps 
\begin{equation*} \ev_i: \;\IM_{0,0,r}(\gamma^{\vec{m}^+},\gamma^{\vec{m}^-})/\IR \;\to\; V,\;\; i=1,...,r \end{equation*}
given by mapping the tuple $(h,j,\mu,\uw)\in\IM_{0,0,r}(\gamma^{\vec{m}^+},\gamma^{\vec{m}^-})/\IR$ to $h(w_i)\in V$, which extend to the compactified moduli space 
$\overline{\IM_{0,0,r}(\gamma^{\vec{m}^+},\gamma^{\vec{m}^-})/\IR}$. Since we still cannot expect the moduli space $\IM_{0,0,r}(\gamma^{\vec{m}^+},\gamma^{\vec{m}^-})$ 
to be transversally cut out by the Cauchy-Riemann operator, we must proceed as before and choose 
coherent sections $\bar{\nu}$ in the cokernel bundles $\overline{\Coker}\CR_J$ over the compactified moduli spaces 
$\CM=\overline{\IM_{0,0,r}(\gamma^{\vec{m}^+},\gamma^{\vec{m}^-})/\IR}$ to obtain the regular moduli spaces 
\begin{equation*}
    \overline{\IM_{0,0,r}(\gamma^{\vec{m}^+},\gamma^{\vec{m}^-})/\IR}^{\bar{\nu}} =\bar{\nu}^{-1}(0) \subset 
    \overline{\IM_{0,0,r}(\gamma^{\vec{m}^+},\gamma^{\vec{m}^-})/\IR}. 
\end{equation*}
$ $\\

Assigning to each chosen differential form $\theta_i\in\Omega^*(V)$ a graded formal variable $t_i$ with $\deg t_i= \deg \theta_i - 2$ 
and abbreviating $p_m = p_{\gamma^m}$ and $q_m = q_{\gamma^m}$ we let $\IP_{\gamma}$ be the graded Poisson algebra of formal power series in the variables 
$p_m$ with coefficients which are polynomials in the $q_m$'s and formal power series of the $t_i$'s. Following [EGH] we define the 
generating function $h_{0,\gamma}\in\IP_{\gamma}$ by 
\begin{equation*} 
 \hh_{0,\gamma} = \sum_{\vec{m}^{\pm},\vec{i}} \frac{1}{n^+! n^-! r!} \int_{\overline{\IM_{0,0,r}(\gamma^{\vec{m}^+},\gamma^{\vec{m}^-})/\IR}^{\bar{\nu}}} 
 \ev_1^*\theta_{i_1}\wedge ... \wedge \ev_r^*\theta_{i_r}\;\;\; p_{\vec{m}^+} q_{\vec{m}^-} t_{\vec{i}}.
\end{equation*}
$ $\\
{\bf Theorem 4.3:} {\it For a chosen string of closed differential forms $\Theta=(\theta_1,...,\theta_N)\in(\Omega^*(V))^N$ 
the generating function $\hh_{0,\gamma}\in\IP_{\gamma}$ 
is given by}
\begin{equation*} \hh_{0,\gamma} = \sum_{i: \deg\theta_i=1} \sum_{m\in\IN} m \int_{\gamma} \theta_i  \;\cdot\; p_m q_m t_i. \end{equation*}
$ $\\
{\it Proof:} Since the positions of the marked points are not fixed, it follows that the dimension of the regular moduli space 
$\IM_{0,0,r}(\gamma^{\vec{m}^+},\gamma^{\vec{m}^-})^{\nu}$ is given by $2r$ plus the dimension of the underlying regular moduli space 
$\IM_{0,0}(\gamma^{\vec{m}^+},\gamma^{\vec{m}^-})^{\nu}$ with no additional marked points. In the case of true branched covers, i.e., 
$\IM_{0,0}(\gamma^{\vec{m}^+},\gamma^{\vec{m}^-})\neq \IM_{0,0}(\gamma^m,\gamma^m)$ it follows that $\IM_{0,0,r}(\gamma^{\vec{m}^+},\gamma^{\vec{m}^-})^{\nu}$ 
has dimension greater or equal to $2r+1$. In other words, the top stratum of $\overline{\IM_{0,0,r}(\gamma^{\vec{m}^+},\gamma^{\vec{m}^-})/\IR}^{\bar{\nu}}$ has 
dimension greater or equal to $2r$, which in turn must agree with the degree of the differential form $\ev_1^*\theta_{i_1}\wedge ... \wedge \ev_r^*\theta_{i_r}$ 
in order to get a nonzero contribution to $\hh_{0,\gamma}$. In particular, at least one differential form $\theta_{i_k}$, $k\in\{1,...,r\}$ must have degree greater 
or equal to two. On the other hand, observing that the image of the evaluation map $\ev_k$ from $\overline{\IM_{0,0,r}(\gamma^{\vec{m}^+},\gamma^{\vec{m}^-})/\IR}$ 
to $V$ is clearly 
contained in the closed Reeb orbit $\gamma$ and that the pullback of a form on $V$ under the inclusion map $\gamma \hookrightarrow V$ 
is nonzero only for forms of degree zero or one, it follows that $\ev_k^*\theta_{i_k}=0$. So, while we have shown in this paper that moduli spaces of true branched covers 
without additional marked points do not contribute to the generating function $\hh_{0,\gamma}$, it follows from the last observation that this remains true 
when we introduce additional marked points and differential forms by simple topological reasons. Finally, observe that for moduli spaces of trivial cylinders the 
top stratum of $\overline{\IM_{0,0,r}(\gamma^m,\gamma^m)/\IR}$ has dimension $2r-1$, so that here we might get nonzero contributions from moduli spaces 
with one additional marked point if the corresponding differential form has degree one. Since the moduli spaces of trivial cylinders are automatically regular, 
it is easily seen that this contribution is given by integrating the one-form along the closed Reeb orbit $\gamma$. $\qed$ \\
  
Observe that the generating function is in general no longer equal to zero when a string of 
closed differential forms is chosen, which implies that the differential in rational symplectic field theory and contact homology is no longer strictly decreasing 
with respect to the action filtration, where we have set $\FF(t_i)=0$ for each formal variable $t_i$. 
However, in order to show how theorem 2.4.3 can be used to compute SFT invariants, we follow [FOOO] in employing   
the spectral sequence for filtered complexes, where for simplicity we restrict our attention only to the computation of the contact homology for contact manifolds 
and symplectic mapping tori. Recall from the introduction that in both cases the contact homology is indeed well-defined.  \\
\\
{\bf Corollary 4.4:} {\it Let $V$ be a contact manifold or a symplectic mapping torus. Then there exists a spectral sequence $(E^r,d^r)$ computing the 
contact homology, $E^{\infty}=H_*(\IA,\del)$, where the $E^2$-page is given by the graded commutative algebra $\IA_0$ which, in contrast to $\IA$, is now only 
freely generated by the formal variables $q_{\gamma}$ with $\int_{\gamma}\theta_i=0$ for all $i=1,...,N$.} \\
\\
{\it Proof:} First observe that it follows from the theorem of Arzela-Ascoli that for any given maximal period $T>0$ the set of closed Reeb orbits of period $\leq T$ 
is compact. Together with the assumption that the contact one-form $\lambda$ is chosen generically in the sense that every closed orbit is nondegenerate and hence isolated, 
it follows that the number of closed orbits with period less or equal $T$ is finite for every $T>0$, so that, in particular, the action spectrum 
$\{\int_{\gamma}\lambda: \gamma\in P(V)\}$ is a discrete subset of $\IR^+$. Note that this automatically implies that the set of action values $\FF(q^{\Gamma})\in\IR^+$, 
$\Gamma\subset P(V)$ is discrete, and hence can be identified with the discrete set $\{a_1,a_2,...\}\subset\IR^+$ with $a_k \leq a_{k+1}$. \\

Using this we equip the chain complex $(\IA,\del)$ underlying contact homology with a filtration $(\FF^k\IA)_{k\in\IN}$ by requiring that $\FF^k\IA$ is spanned by 
monomials $q^{\Gamma}$ with $\FF(q^{\Gamma})\leq a_k$. Note that it follows from the fact all curves have nonnegative contact area that the 
differential is indeed respecting the filtration, $\del: \FF^k\IA\to\FF^k\IA$. Now we can use as in [FOOO] the spectral sequence $(E^r,d^r)$ for filtered complexes to 
compute the homology of $(\IA,\del)$. In order to see how the theorem implies the corollary it suffices 
to observe that the differential $d^1:E^1_{k,\ell}\to E^1_{k,\ell-1}$ agrees with the part $\del_0$ of the differential $\del:\IA\to\IA$, which is counting 
only curves with zero contact area, i.e., orbit curves. Hence $E^2=H_*(\IA,\del_0)$ and it is easily deduced from the fact that $\del_0$ satisfies the Leibniz 
rule that the latter agrees with $\IA_0$ as defined above. \\     

On the other hand, for symplectic mapping tori one can use the splitting of the chain complex with respect to the total period and again use the compactness of the 
set of closed orbits of bounded period to get discreteness of the action spectrum. $\qed$


\begin{thebibliography}{10000000}

\bibitem[BEHWZ]{BEHWZ} Bourgeois, F., Eliashberg, Y., Hofer, H., Wysocki, K. and Zehnder, E.: {\it Compactness results in 
      symplectic field theory.} Geom. and Top. {\bf 7}, 2003. 

\bibitem[CM]{CM} Cieliebak, K. and Mohnke, K.:{\it Compactness for punctured holomorphic curves.} J. Symp. Geom. {\bf 3}(4), 2005.

\bibitem[B]{B} Bourgeois, F.: {\it A Morse-Bott approach to contact homology.} Ph.D. thesis, Stanford University, 2002.

\bibitem[BM]{BM} Bourgeois, F. and  Mohnke, K.: {\it Coherent orientations in symplectic field theory.} Math. Z. {\bf 248}, 2003.

\bibitem[EGH]{EGH} Eliashberg, Y., Givental, A. and Hofer, H.: {\it Introduction to symplectic field theory.} 
      Geom. Funct. Anal., Special Volume, Part II, 2000.      

\bibitem[EKP]{EKP} Eliashberg, Y., Kim, S. and Polterovich, L.: {\it Geometry of contact transformations and domains: 
      orderability vs. squeezing.} Geom. and Top. {\bf 10}, 2006.

\bibitem[F]{F} Fabert, O.: {\it Contact homology of Hamiltonian mapping tori.} ArXiv preprint (math.SG/0609406), 2006.

\bibitem[FO]{FO} Fukaya, K. and Ono, K.: {\it Arnold conjecture and Gromov-Witten invariants for general symplectic manifolds.} 
      Fields Inst. Comm. {\bf 24}, 1999.

\bibitem[FOOO]{FOOO} Fukaya, K., Oh, Y., Ohta, H. and Ono, K.: {\it Lagrangian intersection Floer theory - anomaly and obstruction.} preprint, 2000.
 
\bibitem[HT1]{HT1} Hutchings, M. and Taubes, C.: {\it Gluing pseudoholomorphic curves along branched covered cylinders I.} 
      ArXiv preprint (math.SG/0701300), 2007.
 
\bibitem[HT2]{HT2} Hutchings, M. and Taubes, C.: {\it Gluing pseudoholomorphic curves along branched covered cylinders II.} 
      ArXiv preprint (math.SG/0705.2074), 2007.

\bibitem[HWZ]{HWZ} Hofer, H., Wysocki, K. and Zehnder, E.: {\it A general Fredholm theory I: A splicing-based differential geometry.} 
      ArXiv preprint (math.FA/0612604), 2006. 

\bibitem[L]{L} Long, Y.: {\it Index theory for symplectic paths with applications.} Progress in Mathematics {\bf 207}, 
      Birkhauser, 2002. 

\bibitem[LT]{LT} Li, J. and Tian, G.: {\it Virtual moduli cycles and Gromov-Witten invariants of general symplectic manifolds.} 
      First Int. Press Lect. Ser. {\bf I}, 1998.
 
\bibitem[LiuT]{LiuT} Liu, G. and Tian, G.: {\it Floer homology and Arnold conjecture.} J. Diff. Geom. {\bf 49}, 1998.

\bibitem[MD]{MD} McDuff, D.: {\it The virtual moduli cycle.} AMS Transl. {\bf 196}(2), 1999. 

\bibitem[MDSa]{MDSa} McDuff, D. and Salamon, D.A.: {\it $J$-holomorphic curves and symplectic topology.} AMS Colloquium Publications, 
      Providence RI, 2004.

\bibitem[Sch]{Sch} Schwarz, M.: {\it Cohomology operations from $S^1$-cobordisms in Floer homology.} Ph.D. thesis, Swiss Federal Inst. of 
        Techn. Zurich, Diss. ETH No. 11182, 1995.

\end{thebibliography}
\end{document}
